\documentclass[11pt]{amsart}
\usepackage{amsmath,amssymb,amsthm}
\usepackage{mathrsfs}
\usepackage{microtype}
\usepackage[margin=1in]{geometry}
\usepackage[colorlinks=true,linkcolor=blue,citecolor=blue,urlcolor=blue]{hyperref}

\newtheorem{theorem}{Theorem}[section]
\numberwithin{equation}{section}
\newtheorem{proposition}[theorem]{Proposition}
\newtheorem{lemma}[theorem]{Lemma}
\newtheorem{corollary}[theorem]{Corollary}
\theoremstyle{definition}
\newtheorem{definition}[theorem]{Definition}
\newtheorem{example}[theorem]{Example}
\theoremstyle{remark}
\newtheorem{remark}[theorem]{Remark}
\newtheorem{problem}[theorem]{Problem}
\newtheorem{conjecture}[theorem]{Conjecture}

\newcommand{\C}{\mathbb{C}}
\newcommand{\R}{\mathbb{R}}
\newcommand{\Z}{\mathbb{Z}}
\newcommand{\T}{\mathbb{T}}
\newcommand{\A}{\mathscr{A}}
\newcommand{\Log}{\operatorname{Log}}
\newcommand{\Arg}{\operatorname{Arg}}
\newcommand{\MV}{\operatorname{MV}}
\newcommand{\Vol}{\operatorname{Vol}}
\newcommand{\MVol}{\operatorname{MultiVol}}
\newcommand{\conj}{\operatorname{conj}}
\newcommand{\Hm}{\mathcal{H}}
\newcommand{\G}{\mathcal{G}}
\newcommand{\Sg}{\mathcal{S}}
\newcommand{\Mv}{\mathcal{M}}
\newcommand{\Vc}{\mathcal{V}}
\newcommand{\Fix}{\operatorname{Fix}}
\newcommand{\re}{\operatorname{Re}}

\newcommand{\rk}{\operatorname{rank}}
\newcommand{\Trop}{\operatorname{Trop}}

\title{The topology of the amoeba map and maximal volume} 
\author{Mounir Nisse}

 \address{Mounir Nisse\\
Department of Mathematics, Xiamen University Malaysia, Jalan Sunsuria, Bandar Sunsuria, 43900, Sepang, Selangor, Malaysia.
}
\email{\href{mailto:mounir.nisse@gmail.com}{mounir.nisse@gmail.com}, {mounir.nisse@xmu.edu.my}}
\urladdr{https://sites.google.com/site/mounirnisse1/}

\thanks{This research is supported in part by Xiamen University Malaysia Research Fund (Grant no. XMUMRF/2025-C15/IMAT/0035).}
 
\subjclass[2020]{Primary 14T20, 14P25, 32A60, 52A39, 14P10} 
\keywords{Amoebas, coamoebas, multivolume, mixed volume, Hausdorff measure, real algebraic varieties, Harnack curves, linear spaces, symmetric powers, Cauchy determinants, transversality}

\begin{document}

\begin{abstract}
For a smooth complete intersection $V\subset(\C^*)^n$ of dimension $k$, $n\ge2k$, we bound the volume of its amoeba by $\pi^{2k}$ times a sum of mixed volumes, divided by the order of any finite group $\Gamma$ of symmetries commuting with the logarithmic map. For $n=2k\ge4$ almost every fibre of the amoeba map has a number of points divisible by $2|\Gamma|$. For $n=2k$ we characterise the varieties of maximal volume, show that for $k\ge2$ they are not real, and prove that no plane in general position in $(\C^*)^4$ has maximal volume. For $n>2k$ the bound is never attained, varieties of maximal volume have amoeba maps injective off a null set, and injectivity is generic. The amoeba map of almost every real $k$-plane is two-to-one, and that of the $k$-th symmetric power of a real line exactly $2^k$-to-one; conversely, real $k$-planes in general position in $(\C^*)^{2k}$ with $2^k$-to-one amoeba maps are symmetric powers. We also classify the real planes of $(\C^*)^4$ and study lines in $(\C^*)^3$.
\end{abstract}

\maketitle


\section{Introduction}

Let $V\subset(\C^*)^n$ be an algebraic variety of pure complex dimension $k$. Its amoeba, introduced by Gelfand, Kapranov and Zelevinsky \cite[Chapter 6]{GKZ94}, is the image $\A(V)=\Log(V)\subset\R^n$ of $V$ under the map $\Log(z_1,\dots,z_n)=(\log|z_1|,\dots,\log|z_n|)$, and its coamoeba is the image of $V$ under the argument map $\Arg(z_1,\dots,z_n)=(\arg z_1,\dots,\arg z_n)\in\T^n=(\R/2\pi\Z)^n$, see \cite{Mik04,NP17}. Since $V$ has real dimension $2k$, the amoeba has dimension at most $\min\{2k,n\}$, and for a generic linear space the dimension is exactly $\min\{2k,n\}$ by Nisse and Passare \cite{NP17}. When $n\ge 2k$ the natural measure of the size of the amoeba is therefore its $2k$-dimensional Hausdorff measure $\Hm^{2k}(\A(V))$, which we simply call its volume.

The first result on this volume is due to Passare and Rullg{\aa}rd \cite{PR04}, who proved that the area of the amoeba of a plane curve $\{f=0\}\subset(\C^*)^2$ is at most $\pi^2\Vol(\Delta_f)$, where $\Delta_f$ is the Newton polygon of $f$; Mikhalkin and Rullg{\aa}rd \cite{MR01} showed that equality holds exactly for the simple Harnack curves of \cite{Mik00}, up to a multiplicative translation. Madani and Nisse \cite{MN13} proved that the $2k$-volume of the amoeba of a $k$-dimensional algebraic variety in $(\C^*)^n$ with $n\ge 2k$ is finite; for rational curves they gave the explicit bound $\pi^2\sum_{i<j}n_in_j$, where $n_j$ is the number of zeros and poles of the $j$-th coordinate of a parametrisation, while for $k\ge2$ their bound involves constants which are not computed in terms of the variety, and in \cite{MN15} it is shown that for generic analytic $k$-dimensional varieties finiteness of the amoeba volume characterises algebraicity. In the half-dimensional case $n=2k$, Mikhalkin \cite{Mik17} proved that $\Vol(\A(V))\le\frac{\pi^{2k}}2\alpha(V)$, where $\alpha(V)$ is the toric intersection number of $V$ with its complex conjugate, and that for a complete intersection with Newton polytopes $\Delta_1,\dots,\Delta_k\subset\R^{2k}$ one has $\alpha(V)=\MV(\Delta_1,\dots,\Delta_k,\Delta_1,\dots,\Delta_k)$. The expected value of this volume for random complete intersections was studied by Bayraktar and Ki\c{s}isel \cite{BK23}, and the amoebas and coamoebas of linear spaces were described by Nisse and Passare \cite{NP17}, who computed in particular the volume $\pi^{2k}$ of the coamoeba of a generic $k$-plane and stated that the amoeba of a real generic $k$-plane in $(\C^*)^{2k}$ has volume $\pi^{2k}/2^k$ \cite[Theorem 1.2(ii)]{NP17}, deducing it from the assertion that its amoeba map is exactly $2^k$-to-one \cite[Proposition 5.3]{NP17}. The value is correct for $k\le2$. For $k\ge3$ we prove that the assertion of \cite[Proposition 5.3]{NP17} fails for every real $k$-plane in general position which is not a symmetric power of a line (Corollary \ref{cor:fibre4}). Since the present author is a coauthor of \cite{NP17}, these results should be read as a correction of Proposition 5.3 of that paper; the point of its proof which fails is identified in Section \ref{sec:why}.

Mikhalkin's bound, and the results of the present paper on the counting function of the amoeba map and on the varieties of maximal volume, were conjectured by the author in his talk at the memorial conference for Mikael Passare in summer 2013. Bounds of this type are attributed to a conjecture of the author in \cite[Remark 1.6]{Mik17}, where Mikhalkin's bound is proved. The present paper proves the statements on the counting function, in the form of the divisibility of the fibres of the amoeba map by twice the order of the symmetry group of the variety (Theorems \ref{thm:half} and \ref{thm:proj}), and on the varieties of maximal volume, in the form of their characterisation by three conditions and of the nonexistence of planes of maximal volume in general position in $(\C^*)^4$ (Theorems \ref{thm:char} and \ref{thm:planesgp}), and it extends them to the case $n>2k$ (Theorems \ref{thm:strict}, \ref{thm:emb}, \ref{thm:max} and \ref{thm:geninj}).

This paper treats the case $n>2k$, in which the amoeba is a set of positive codimension and Mikhalkin's topological argument, based on the properness of the amoeba map, is no longer available. Throughout, $V$ is a smooth complete intersection: $V=\{f_1=\dots=f_{n-k}=0\}$ for Laurent polynomials $f_m$ with Newton polytopes $\Delta_m\subset\R^n$, whose differentials are linearly independent at every point of $V$; note that a complete intersection of dimension $k$ needs $n-k$ equations, so that a line in $(\C^*)^3$ is cut out by two equations and a $k$-plane in $(\C^*)^{2k+1}$ by $k+1$ equations.

Our first result, Theorem \ref{thm:main}, is the estimate

$$\Hm^{2k}(\A(V))\le\sum_{|I|=2k}\MVol_I(V)\le\pi^{2k}\sum_{|I|=2k}\alpha_I(\Delta_1,\dots,\Delta_{n-k}),$$

where $I$ runs over the subsets of $\{1,\dots,n\}$ of cardinality $2k$, $\MVol_I(V)$ is the volume of the projected amoeba $p_I(\A(V))\subset\R^{2k}$ counted with multiplicity, and $\alpha_I$ is the mixed volume in $\R^{2n-2k}$ of the polytopes $\Delta_m$ placed in two different ways in the coordinates $(a_I,a_{\hat I},0)$ and $(a_I,0,a_{\hat I})$, with $\hat I$ the complement of $I$. The proof combines the Cauchy--Binet identity, the area formula of geometric measure theory \cite{Fed69}, the equality of the Jacobians of the logarithmic and argument maps on a complex manifold, and a count of the points of the fibres of the rolled coamoeba map by Bernstein's theorem \cite{Ber75} applied to a doubled system, in the spirit of \cite{Mik17}. For polynomials of degrees $d_1,\dots,d_{n-k}$ it gives $\Hm^{2k}(\A(V))\le\binom n{2k}\pi^{2k}(d_1\cdots d_{n-k})^2$, and it is sharp in each projection for linear spaces in general position. For a line in $(\C^*)^3$ it gives $3\pi^2$, while the area lies strictly between $\sqrt3\,\pi^2$ and $3\pi^2$ (Theorem \ref{thm:line}). For $k\ge2$ it is, to our knowledge, the first bound with explicit constants.

Our second result concerns the factor $\frac12$ of Mikhalkin's theorem. Let $\G$ be the group of maps $z\mapsto\theta z$ and $z\mapsto\theta\bar z$ with $\theta$ in the compact torus; all of them commute with $\Log$ in the sense that $\Log\circ g=\Log$. If $\Gamma$ is a finite group of such maps preserving $V$, then the fibres of the amoeba map are unions of free $\Gamma$-orbits off a negligible set, and Theorem \ref{thm:sym} gives $\Hm^{2k}(\A(V))\le|\Gamma|^{-1}\MVol(V)\le\pi^{2k}|\Gamma|^{-1}\sum_I\alpha_I$. Since $V$ is real up to a multiplicative translation exactly when $\Gamma$ can be chosen to contain an antiholomorphic involution, the factor $\frac12$ holds for such varieties (Corollary \ref{cor:half}). This factor cannot hold in general for $n>2k$, because the amoeba map of a nonreal line in $(\C^*)^3$ is injective and its area exceeds $\frac{\pi^2}2\sum_I\alpha_I=\frac32\pi^2$ (Theorem \ref{thm:line}). In the half-dimensional case we prove in Theorem \ref{thm:half} that, if $k\ge2$ or $\Gamma$ consists of holomorphic maps, almost every fibre of the amoeba map has a number of points divisible by $2|\Gamma|$. The proof removes the image under $\Log$ of the fixed points of $\Gamma$, which has codimension at least $k\ge2$ and therefore does not disconnect $\R^{2k}$ \cite{HW41}, and uses the vanishing of the mod $2$ degree \cite{Mil65} of the proper map induced by $\Log$ on the quotient by $\Gamma$. For real varieties of dimension $k\ge2$ this gives fibres divisible by $4$ and the factor $\frac14$ (Corollary \ref{cor:real4}), attained by real generic planes in $(\C^*)^4$, for which it gives a new proof that the amoeba map is four-to-one (Corollary \ref{cor:sizes}). Theorem \ref{thm:proj} extends the divisibility to the coordinate projections when $n>2k$: the number of points of $V$ over a point of the projection $p_I(\A(V))$ is divisible by $2|\Gamma|$, which bounds the volume of the projected amoeba by $\MVol_I(V)/2|\Gamma|$, while the projected part of the volume of $\A(V)$ is at most $\MVol_I(V)/|\Gamma|$, with equality exactly when almost every fibre of the amoeba map is a single $\Gamma$-orbit. We show in Theorem \ref{thm:twotoone} that the amoeba map of almost every real $k$-plane in $(\C^*)^n$, $n>2k$, is two-to-one, and that of almost every complex $k$-plane injective, almost everywhere on the amoeba. Hence the projected parts of the volume of almost every real $k$-plane in $(\C^*)^{2k+1}$ equal $\pi^{2k}/2$ (Corollary \ref{cor:kplanes2}). For real planes in $(\C^*)^5$, Theorem \ref{thm:conic} gives an explicit condition: if the six lines cut on the projective closure of $V$ by the coordinate hyperplanes are not tangent to a common conic, which amounts to the linear independence of the three vectors $(1/\gamma_m,1/\delta_{m1},1/\delta_{m2})$ formed from the coefficients $z_m=\gamma_m+\delta_{m1}z_1+\delta_{m2}z_2$ of $V$ in graph form, then the amoeba map is two-to-one outside the image of the real points. If they are tangent to a smooth conic, $V$ is the coordinatewise square $\ell\odot\ell$ of a real line $\ell$ and its amoeba is the Minkowski sum $\A(\ell)+\A(\ell)$ (Corollary \ref{cor:hadamard}). In Theorem \ref{thm:square} we prove that the amoeba map of $\ell\odot\ell\subset(\C^*)^n$, $n\ge4$, is exactly four-to-one: two pairs $\{q_1,q_2\}$, $\{q_3,q_4\}$ of parameters have the same image if and only if $\{q_3,q_4\}$ is obtained from $\{q_1,q_2\}$ by conjugating some of its entries, because a real polynomial of degree at most $3$ vanishing at $n\ge4$ points is zero. Hence $\Hm^4(\A(\ell)+\A(\ell))=\frac14\MVol(\ell\odot\ell)$, which Proposition \ref{prop:squarevol} writes as an explicit integral whose integrand is built from Cauchy determinants. Each projected part equals $\pi^4/4$ (Corollary \ref{cor:squareproj}), the lower bound of Corollary \ref{cor:plane5}. Since five lines in general position in a real projective plane are always tangent to a conic, every real plane in general position in $(\C^*)^4$ is such a square, which gives an elementary proof that its amoeba map is four-to-one and its amoeba has volume $\pi^4/4$ (Corollary \ref{cor:square4}).

Section \ref{sec:squares} treats more generally the $k$-th symmetric power $\ell^{\odot k}=\{\ell(q_1)\odot\dots\odot\ell(q_k)\}$ of a line $\ell=\{(q-\varpi_m)_m\}$, an affine $k$-plane whose amoeba is the Minkowski sum $\A(\ell)^{[k]}$ of $k$ copies of $\A(\ell)$. For a real line and $n\ge2k$ the amoeba map of $\ell^{\odot k}$ is exactly $2^k$-to-one, its fibres being the orbits of the conjugations of the single factors (Theorem \ref{thm:square}). The volume is an explicit integral of Cauchy determinants (Proposition \ref{prop:squarevol}), the projected multivolumes equal $\pi^{2k}$ by a direct computation of the rolled coamoeba, which amounts to the identity of Proposition \ref{prop:coamsym} extending the integral used by Passare \cite{Pas08}, and for $n=2k$ the amoeba has volume exactly $\pi^{2k}/2^k$ (Corollary \ref{cor:squareproj}). Thus the value of \cite[Theorem 1.2(ii)]{NP17} holds on the family of symmetric powers, of codimension $2$ among real $3$-planes of $(\C^*)^6$. By contrast, for $k\ge3$ the $2^k$-to-one property on which it rests fails for every other real $k$-plane in general position (Corollary \ref{cor:fibre4}). For squares of real lines in $(\C^*)^5$, the volume of $\A(\ell)+\A(\ell)$ is a function of three real moduli; we compute its limit when two points of $\ell$ collide in terms of the multivolume of a line of $(\C^*)^3$ (Theorem \ref{thm:collide}), which shows that its supremum is at least $\frac{1+2\sqrt2}4\pi^4$. For nonreal lines the picture changes: the amoeba map of $\ell^{\odot k}$ is injective almost everywhere for almost every nonreal $\ell$ when $n\ge2k+1$, and everywhere when $n\ge k^2+2k$ (so $n\ge8$ for $k=2$) and the points $\varpi_m$ lie on no curve of an explicit space $\mathscr P_k$ of dimension $k^2+2k$, which for $k=1$ consists of lines and recovers the condition $c\notin\R$, and for $k=2$ consists of circular cubics (Theorem \ref{thm:nonrealsq}).

Section \ref{sec:products} asks whether the products of symmetric powers of real lines in disjoint sets of coordinates, which have $2^k$-to-one amoeba maps (Proposition \ref{prop:prodfib}), are the only real $k$-planes with this property. Without any condition the answer is negative: the fibred extensions of Theorem \ref{thm:fibred}, in which lines of $(\C^*)^2$ are attached to a product of symmetric powers through functions whose moduli are invariant under the partial conjugations, have $2^k$-to-one amoeba maps and are not products, already for $k=2$. In general position, which excludes products with more than one factor, the answer is positive in a strong form (Theorem \ref{thm:classif}): a real $k$-plane in general position in $\C^n$, $n\ge2k$, whose amoeba map has at least $2^k$ points in almost every fibre is a symmetric power of a real line. The proof is an induction on $k$ combining three ingredients: near the tentacle of the amoeba in the direction of a coordinate hyperplane, the fibres of the amoeba map have twice as many points as those of the hyperplane section (Lemma \ref{lem:tentacle}), coordinate projections transport lower bounds for the fibres from $n>2k$ to $n=2k$ (Lemma \ref{lem:projfib}), and $2k+1$ points of $\mathbb P^k$ whose projections from two of them lie on rational normal curves lie on a rational normal curve (Lemma \ref{lem:cones}); the case of planes is Theorem \ref{thm:conic}. Consequently every real $k$-plane in general position in $(\C^*)^{2k}$, $k\ge3$, which is not a symmetric power has fibres of fewer than $2^k$ points on a set of positive measure, exactly $4$ points when $k=3$ (Corollary \ref{cor:fibre4}). This proves that \cite[Proposition 5.6]{NP17} fails for every $k\ge3$. It does not by itself decide the volume formula \cite[Theorem 1.2(ii)]{NP17}, since fibres of more than $2^k$ points also occur. Outside general position we classify completely the real planes of $(\C^*)^4$ (Theorem \ref{thm:planes4}): up to monomial maps exchanging the hyperplane at infinity with the coordinate hyperplanes (Lemma \ref{lem:monomial}), every real plane with a nondegenerate amoeba is a square, a fibred extension or a product of two lines, its amoeba map is four-to-one and its amoeba has volume $\pi^4/4$. For real $3$-planes of $(\C^*)^6$ whose dual configuration has a collinear triple, the amoeba map is eight-to-one exactly when the line of the triple meets two opposite edges of the tetrahedron of the remaining points, and then the plane is a fibred extension (Theorem \ref{thm:triple}); this supports Conjecture \ref{conj:k3}, which lists all eight-to-one real $3$-planes.

Our third result is a detailed study of the lines of $(\C^*)^3$, which up to multiplicative translation form the family $L_c=\{(s,1-s,1-cs)\}$, $c\in\C\setminus\{0,1\}$. The three projected multivolumes all equal $\pi^2$, the amoeba map is injective for $c\notin\R$ and two-to-one for $c\in\R$, and this gives the two-sided estimate $\sqrt3\,\pi^2<\Hm^2(\A(L_c))<3\pi^2$ for nonreal lines. We then study the total multivolume $\Mv(c)=\MVol(L_c)$, which equals the area for nonreal $c$ and twice the area for real $c$. It depends only on the similarity class of the triangle with vertices $0,1,c$ (Proposition \ref{prop:sym12}). It is continuous, and real analytic off the real axis (Proposition \ref{prop:reg}). It tends to $(1+\sqrt2)\pi^2$ when two vertices of the triangle collide (Theorem \ref{thm:degen}), with the correction $-\kappa(\arg c)\sqrt{|c|}+O(|c|^{2/3})$ of Proposition \ref{prop:rate}, uniformly in $\arg c$, where $\kappa>0$ involves the constant $4\pi^{3/2}/\Gamma(\frac14)^2$, so that $\Mv(c)<(1+\sqrt2)\pi^2$ near the degenerate values. Its critical points include the equilateral and the harmonic configurations (Proposition \ref{prop:crit}).

\subsection*{Varieties of maximal volume}
Let $\Delta_1,\dots,\Delta_{n-k}\subset\R^n$ be lattice polytopes with $n\ge2k$, and let $\Vc(\Delta)$ be the set of smooth complete intersections $V=\{f_1=\dots=f_{n-k}=0\}\subset(\C^*)^n$ of dimension $k$, the $f_m$ being Laurent polynomials with Newton polytopes $\Delta(f_m)=\Delta_m$ whose differentials are linearly independent at every point of $V$. We regard $\Vc(\Delta)$ as a subset of the space of coefficient vectors (see before Theorem \ref{thm:geninj}). By Theorem \ref{thm:main}, the area formula \cite[Theorem 3.2.3]{Fed69} and \eqref{eq:mvsum},
\begin{equation}\label{eq:chain}
\Hm^{2k}\big(\A(V)\big)\le\MVol(V)\le\sum_{|I|=2k}\MVol_I(V)\le\pi^{2k}\sum_{|I|=2k}\alpha_I,
\end{equation}
where the multivolumes $\MVol(V)$ and $\MVol_I(V)$ are defined in Section \ref{sec:notation} and $\alpha_I$ is the mixed volume of Definition \ref{def:alpha}. When $n=2k$, $\alpha_I=M:=\MV(\Delta_1,\dots,\Delta_k,\Delta_1,\dots,\Delta_k)$ and Mikhalkin's sharper bound $\Vol(\A(V))\le\frac{\pi^{2k}}2M$ holds \cite[Theorem 1]{Mik17}.

We characterise the varieties realising $\sup_{\Vc(\Delta)}\Hm^{2k}(\A(V))$ through the antiholomorphic involutions $\sigma_\varepsilon(z)=\varepsilon\bar z$, $\varepsilon\in\T^n=(S^1)^n$, and $\sigma'_x(z)=e^{2x}/\bar z$, $x\in\R^n$, taken coordinatewise; the involution $\sigma'_x$ is denoted $\conj'_x$ in Section \ref{sec:half}. Since $\Fix(\sigma'_x)=\Log^{-1}(x)$, the fibre of the amoeba map over $x$ is $V\cap\Fix(\sigma'_x)$, and since $\Fix(\sigma_\varepsilon)=\varepsilon^{1/2}(\R^*)^n$, for $n=2k$ the fibres of the rolled coamoeba map $\Arg_\pi|_V$ are the sets $V\cap\Fix(\sigma_\varepsilon)$. With the counting function $N(x)=\#(V\cap\Log^{-1}(x))$, the area formula gives $\MVol(V)=\int_{\R^n}N\,d\Hm^{2k}$ and
\begin{equation}\label{eq:area}
\Hm^{2k}(\Log(S))\le\int_SJ(\Log|_V)\,d\mathrm{vol}\quad\text{for every Borel set }S\subset V.
\end{equation}

Let $n=2k$. We say that $V$ has maximal volume if $\Vol(\A(V))=\frac{\pi^{2k}}2M$ and $M>0$. Let $\alpha(V)$ be Mikhalkin's toric intersection number of $V$ with $\conj(V)$, that is, the number of points, counted with multiplicity, of $V\cap\varepsilon\conj(V)$ for generic $\varepsilon\in(\C^*)^{2k}$ \cite[Definition 1.3]{Mik17}.

\begin{theorem}\label{thm:char}
Let $n=2k$, $V\in\Vc(\Delta)$ and $M>0$. Then $V$ has maximal volume if and only if the following three conditions hold. $(\mathrm A)$: $\alpha(V)=M$; this holds whenever $V$ is nondegenerate in the sense of Khovanskii. $(\mathrm B)$: for almost every $\varepsilon\in\T^{2k}$, every point of $V\cap\sigma_\varepsilon(V)$ is fixed by $\sigma_\varepsilon$; equivalently, $\Arg_\pi|_V$ has exactly $\alpha(V)$ points in almost every fibre, that is, $V$ is multiHarnack in the sense of \cite[Definition 3.2]{Mik17}. $(\mathrm C)$: for almost every $x\in\A(V)$, exactly two points of $V\cap\sigma'_x(V)$ are fixed by $\sigma'_x$; equivalently, $\Log|_V$ is two-to-one almost everywhere on $\A(V)$, and $\Log$ restricted to the preimage of the complement of the critical values in $\A(V)$ is a two-sheeted covering.
\end{theorem}

For curves, $(\mathrm A)$, $(\mathrm B)$, $(\mathrm C)$ characterise the simple Harnack curves (Section \ref{sec:harnack}). The inequalities $\Vol(\A(V))\le\frac12\MVol(V)\le\frac{\pi^{2k}}2\alpha(V)$ are due to Mikhalkin \cite[Lemma 6 and Corollary 7]{Mik17}, and $(\mathrm B)$ is his multiHarnack property \cite[Definition 3.2]{Mik17}; what is new in Theorem \ref{thm:char} is the identification of the three equality cases, the sufficient condition for $(\mathrm A)$ (Proposition \ref{prop:A}), and the covering structure of $(\mathrm C)$.

The existence of varieties of maximal volume for $k\ge2$ is open \cite[Problems 3.1 and 3.3]{Mik17}. We show that for $k\ge2$ they have trivial stabiliser and are not real (Corollary \ref{cor:nomaxreal}), that products of simple Harnack curves satisfy $(\mathrm A)$ and $(\mathrm B)$ but not $(\mathrm C)$ (Example \ref{ex:products}), and that no plane in general position in $(\C^*)^4$ has maximal volume (Theorem \ref{thm:planesgp}).

\smallskip\noindent\emph{The case $n>2k$.} The parity which forces $N\ge2$ when $n=2k$ disappears, and the amoeba map of a variety of large volume must be injective.

\begin{theorem}\label{thm:strict}
Let $n>2k$ and $V\in\Vc(\Delta)$. The first inequality of \eqref{eq:chain} is an equality if and only if $N=1$ almost everywhere on $\A(V)$ for $\Hm^{2k}$. If a connected component of $V$ on which $J(\Log|_V)\not\equiv0$ lies in no translate $T^I_t=\{z:z_j=t_j\text{ for }j\notin I\}$ of a $2k$-dimensional coordinate subtorus, the second inequality is strict. The third inequality is an equality if and only if condition $(\mathrm B_I)$ of Proposition \ref{prop:BI} holds for every $I$. Finally, if $\sum_I\alpha_I>0$, then $\Hm^{2k}(\A(V))<\pi^{2k}\sum_{|I|=2k}\alpha_I$ for every $V\in\Vc(\Delta)$.
\end{theorem}

Condition $(\mathrm B_I)$ is the analogue of $(\mathrm B)$ for the projection to $(\C^*)^I$; the part of \eqref{eq:chain} which detects the maximum is the first inequality.

\begin{theorem}\label{thm:emb}
Let $n>2k$ and $V\in\Vc(\Delta)$. The following conditions are equivalent. $(\mathrm I_1)$: $\Hm^{2k}(\A(V))=\MVol(V)$. $(\mathrm I_2)$: for $\Hm^{2k}$-almost every $x\in\A(V)$, exactly one point of $V$ is fixed by $\sigma'_x$. $(\mathrm I_3)$: there is a closed set $Z\subset\A(V)$ with $\Hm^{2k}(Z)=0$ such that $\A(V)\setminus Z$ is a closed embedded submanifold of dimension $2k$ of $\R^n\setminus Z$ and $\Log$ is a diffeomorphism from $V\setminus\Log^{-1}(Z)$ onto $\A(V)\setminus Z$. If they hold and $\MVol(V)>0$, then the stabiliser $\G_V$ of $V$ in the group $\G$ of the maps $z\mapsto\theta z$, $z\mapsto\theta\bar z$, $\theta\in\T^n$, is trivial, $V$ is not real up to a multiplicative translation, and no connected component of $V$ with $J(\Log|_V)\not\equiv0$ lies in a translate of a $2k$-dimensional coordinate subtorus.
\end{theorem}

Put $m^*(\Delta)=\sup_{V\in\Vc(\Delta)}\MVol(V)$, which is finite by \eqref{eq:chain}. We say that $\Vc(\Delta)$ has the injectivity property if the set of $V\in\Vc(\Delta)$ satisfying the conditions of Theorem \ref{thm:emb} is dense in $\Vc(\Delta)$.

\begin{theorem}\label{thm:max}
Let $n>2k$. Then $\sup_{\Vc(\Delta)}\Hm^{2k}(\A(V))\le m^*(\Delta)$, and a variety $V\in\Vc(\Delta)$ satisfies $\Hm^{2k}(\A(V))=m^*(\Delta)$ if and only if its amoeba map satisfies the conditions of Theorem \ref{thm:emb} and $\MVol(V)=m^*(\Delta)$. If $\Vc(\Delta)$ has the injectivity property, then $\sup_{\Vc(\Delta)}\Hm^{2k}(\A(V))=m^*(\Delta)$. In that case the varieties of maximal $2k$-volume are exactly the varieties with an amoeba map which is an embedding off a null set and with maximal multivolume. They are not real up to translation, have trivial stabiliser in $\G$, and satisfy $\MVol(V)<\sum_I\MVol_I(V)$.
\end{theorem}

Theorem \ref{thm:max} describes the supremum, which need not be attained. For the lines of $(\C^*)^3$, Theorem \ref{thm:degen} and Proposition \ref{prop:rate} show that the area of the amoeba of a nonreal line tends to $(1+\sqrt2)\pi^2$ when the line degenerates, and stays strictly below this value for all lines close enough to the degenerate ones.

When $n>2k$, pairs of points of $V$ with the same image under $\Log$ form a set of expected dimension $4k-n<2k$, and we turn this into a theorem. Let $A_m=\Delta_m\cap\Z^n$, let $\C^A=\prod_m\C^{A_m}$ be the space of coefficient vectors $c=(c_{m,a})$, $f_{c,m}(z)=\sum_{a\in A_m}c_{m,a}z^a$, and let $V_c$ be the zero set of $f_{c,1},\dots,f_{c,n-k}$. From now on we regard $\Vc(\Delta)$ as the set of $c\in\C^A$ such that every vertex coefficient of every $f_{c,m}$ is nonzero and such that $df_{c,1},\dots,df_{c,n-k}$ are linearly independent along $V_c$; it contains a nonempty Zariski open subset of $\C^A$ (Lemma \ref{lem:U}). For $m=1,\dots,n-k$ let $L_m\subset\Z^n$ be the subgroup generated by the differences $a-b$, $a,b\in A_m$, and for $S\subset\{1,\dots,n-k\}$ put $L_S=\sum_{m\in S}L_m$, with $L_\emptyset=0$.

\begin{theorem}\label{thm:geninj}
Let $n>2k$ and assume that the polytopes satisfy the following two conditions. $(\mathrm L_1)$: $L_{\{1,\dots,n-k\}}=\Z^n$. $(\mathrm L_2)$: for every subset $S\subsetneq\{1,\dots,n-k\}$, $\rk L_S\ge2|S|+2k-n+1$. Then there is a proper real algebraic subset $\mathcal B\subset\C^A$ such that, for every $c\in\Vc(\Delta)\setminus\mathcal B$, the amoeba map of $V_c$ is injective at $\Hm^{2k}$-almost every point of $\A(V_c)$. In particular $\Vc(\Delta)$ has the injectivity property, $\sup_{\Vc(\Delta)}\Hm^{2k}(\A(V_c))=m^*(\Delta)$, and $V_c$ has maximal $2k$-dimensional volume if and only if its amoeba map is an embedding off a closed null set and its multivolume is maximal.
\end{theorem}

For $S=\emptyset$ condition $(\mathrm L_2)$ is exactly $n>2k$. Both conditions hold when every $\Delta_m$ contains $0,e_1,\dots,e_n$, in particular for linear spaces (Corollary \ref{cor:planes}) and for complete intersections of general hypersurfaces.

\subsection*{Organisation of the paper}
Section \ref{sec:notation} fixes notation, proves the two linear algebra facts on which everything rests and collects some preliminaries. Section \ref{sec:counting} counts the fibres of the rolled coamoeba map and Section \ref{sec:main} proves the main estimate. Section \ref{sec:sym} introduces the symmetries and proves the factor $\frac12$, the divisibility of the fibres in the half-dimensional case, the explanation of the failure of the count $2^k$ by partial conjugations, the extension of the divisibility to projections, and the two-to-one property of the amoeba maps of real linear spaces with the explicit conic condition for planes. Section \ref{sec:maximal} characterises the varieties of maximal volume in the half-dimensional case, shows that for $k\ge2$ they are not real, and proves that no plane in general position in $(\C^*)^4$ has maximal volume. Section \ref{sec:squares} studies symmetric powers of lines, real and nonreal, and the volume of $\A(\ell)+\A(\ell)$, and Section \ref{sec:products} studies products of symmetric powers, their fibred extensions, and the classification of real $k$-planes in general position with $2^k$-to-one amoeba maps. Section \ref{sec:lines} treats lines in $(\C^*)^3$ and Section \ref{sec:cross} studies their multivolume as a function of the cross-ratio. Section \ref{sec:kplanes} collects the consequences for $k$-planes in $(\C^*)^{2k+1}$. Sections \ref{sec:chain}, \ref{sec:inj} and \ref{sec:geninj} treat the case $n>2k$: the equality cases of the chain of inequalities \eqref{eq:chain}, the varieties with injective amoeba maps and maximal volume, and the generic injectivity of the amoeba map. Section \ref{sec:curves} rewrites the estimate for curves in terms of eliminants and mixed fibre polytopes, Section \ref{sec:open} collects open questions, and Appendix \ref{app:transv} contains the transversality lemmas used in Section \ref{sec:geninj}.

\bigskip

\noindent{\it Acknowledgements.} The author is  gratefully acknowledges the use the AI system Claude (Anthropic), which used 
 Python, NumPy, SciPy, and LaTeX for the exact computations, and numerical experiments in this work.  


\section{Notation and two linear algebra facts}\label{sec:notation}

Throughout, $(\C^*)^n$ carries the flat metric pulled back by the covering $\C^n\to(\C^*)^n$, $w\mapsto e^w$, that is the metric $\sum_j|dz_j|^2/|z_j|^2$, and every complex submanifold carries the induced Riemannian metric. Writing $w_j=x_j+i\theta_j$ with $x_j=\log|z_j|$ and $\theta_j=\arg z_j$, the maps $\Log$ and $\Arg$ are respectively the real and imaginary parts of the multivalued holomorphic map $w=\log z$, and both are Riemannian submersions onto $\R^n$ and $\T^n$ with their flat metrics. Following Mikhalkin \cite{Mik17} we also use the rolled argument map $\Arg_\pi$, with values in $\T_\pi^n=(\R/\pi\Z)^n$, the composition of $\Arg$ with the $2^n$-sheeted covering $\T^n\to\T_\pi^n$.

For a subset $I\subset\{1,\dots,n\}$ with $|I|=2k$ we denote by $\hat I$ its complement, by $p_I:\R^n\to\R^{I}\simeq\R^{2k}$ the coordinate projection, by $\pi_I:(\C^*)^n\to(\C^*)^{I}$ the corresponding monomial projection, and we put $\Log_I=p_I\circ\Log=\Log\circ\pi_I$, and $\Arg_I$, $\Arg_{\pi,I}$ analogously. Clearly $p_I(\A(V))=\A(\pi_I(V))$. For a $k$-dimensional complex submanifold $V$ we write $\MVol(V)=\int_VJ(\Log|_V)\,d\mathrm{vol},\quad \MVol_I(V)=\int_VJ(\Log_I|_V)\,d\mathrm{vol},$ where $J$ denotes the $2k$-dimensional Jacobian. By the area formula \cite[Theorem 3.2.3]{Fed69}, $\MVol(V)=\int_{\R^n}\#(V\cap\Log^{-1}(x))\,d\Hm^{2k}(x)$ and $\MVol_I(V)=\int_{\R^{2k}}\#(V\cap\Log_I^{-1}(y))\,dy$, so these are the volumes of the amoeba and of its projections counted with the multiplicity of the respective maps; we call them multivolumes.

The first linear algebra fact compares the volume of a set with the volumes of its coordinate projections. Let $T\subset\R^n$ be a real $2k$-dimensional linear subspace and let $M$ be an $n\times 2k$ matrix whose columns form an orthonormal basis of $T$. For $|I|=2k$, the Jacobian of $p_I|_T:T\to\R^I$ is $J_I(T)=|\det M_I|$, where $M_I$ is the submatrix of rows indexed by $I$. The Cauchy--Binet formula gives $\sum_{|I|=2k}(\det M_I)^2=\det(M^{t}M)=1$; this identity is the one used by Madani and Nisse \cite[Section 3]{MN13} to write the volume element of an amoeba. Consequently

\begin{equation}\label{eq:cb}
1=\Big(\sum_{|I|=2k}J_I(T)^2\Big)^{1/2}\le\sum_{|I|=2k} J_I(T)\le \binom{n}{2k}^{1/2},
\end{equation}

the second inequality being the Cauchy--Schwarz inequality. Equality on the left holds only when $T$ is a coordinate $2k$-plane, since it means that at most one $J_I(T)$ is nonzero, and then $J_I(T)=1$ forces $T=\R^I$. More generally, for real $n\times p$ matrices $A$ and $B$ with $p\le n$ the Cauchy--Binet formula reads $\det(A^tB)=\sum_{|I|=p}\det A_I\,\det B_I$, and writing a linear map $L:\R^p\to\R^n$ as $L=MA'$, with $M$ as above and $A'$ square, gives
\begin{equation}\label{eq:CB0}
J(L)^2=\sum_{|I|=p}J(p_I\circ L)^2,
\end{equation}
both sides vanishing if $L$ has rank smaller than $p$; in this paper $p=2k$ always. Applied to the tangent planes of the image of $\Log|_V$, the same identity gives $J(\Log|_V)^2=\sum_IJ(\Log_I|_V)^2$, hence

\begin{equation}\label{eq:mvsum}
\MVol(V)\le\sum_{|I|=2k}\MVol_I(V).
\end{equation}

We shall use repeatedly the following identity, due in substance to Passare \cite{Pas08}. For distinct $q,q'\in\C$,
\begin{equation}\label{eq:passare}
\int_\C\frac{|\sin\theta(s)|}{|s-q|\,|s-q'|}\,dA(s)=\pi^2,
\end{equation}
where $dA$ is the Lebesgue measure of $\C$ and $\theta(s)$ is the angle at $s$ of the triangle with vertices $q,s,q'$. Indeed, the measure $|\sin\theta(s)|\,|s-q|^{-1}|s-q'|^{-1}\,dA(s)$ is invariant under the similarities of $\C$, so we may take $q=0$ and $q'=1$; then the integrand is $|\operatorname{Im}(\bar s^{-1}(s-1)^{-1})|$, which is the Jacobian of the amoeba map of the line $\{z_1+z_2=1\}$ of $(\C^*)^2$ parametrised by $z_1=s$, so that the left-hand side is the multivolume of this line, that is twice the area $\pi^2/2$ of its amoeba \cite{PR04}, the amoeba map of a real line being two-to-one. Written out, \eqref{eq:passare} is the identity $\int_\C|\operatorname{Im}s|\,|s|^{-2}|1-s|^{-2}\,dA(s)=\pi^2$, which Passare used to prove $\zeta(2)=\pi^2/6$ \cite{Pas08}.

The second fact is the following lemma, which in the half-dimensional case is \cite[Lemma 4 and Corollary 5]{Mik17}; for curves it reduces to the classical observation that $dx_1\wedge dx_2=d\theta_1\wedge d\theta_2$ on a holomorphic curve (see \cite{PR04,Mik00}).

\begin{lemma}\label{lem:jac}
Let $V\subset(\C^*)^n$ be a $k$-dimensional complex submanifold, $z\in V$ and $|I|=2k$. Then $J(\Log_I|_V)(z)=J(\Arg_I|_V)(z)=J(\Arg_{\pi,I}|_V)(z)$. In particular $J(\Log|_V)(z)=J(\Arg|_V)(z)$, where these are the $2k$-dimensional Jacobians of maps into $\R^n$ and $\T^n$.
\end{lemma}

\begin{proof}
Choose a unitary basis $u_1,\dots,u_k$ of $T_zV$, so that $u_1,iu_1,\dots,u_k,iu_k$ is a real orthonormal basis of $T_zV$. The differential of $\pi_I\circ\log$ restricted to $T_zV$ is a complex linear map $A:\C^k\to\C^{I}$ whose matrix in these bases is $P+iQ$, with $P,Q$ real $2k\times k$ matrices. The vector $u_j$ is sent to $Pe_j+iQe_j$ and $iu_j$ to $-Qe_j+iPe_j$. Taking real parts, the matrix of $d\Log_I$ in the basis $u_1,\dots,u_k,iu_1,\dots,iu_k$ is $[P\,|\,-Q]$, and taking imaginary parts, the matrix of $d\Arg_I$ is $[Q\,|\,P]$. Now $\det[P\,|\,-Q]=(-1)^k\det[P\,|\,Q]$ and $\det[Q\,|\,P]=(-1)^{k^2}\det[P\,|\,Q]=(-1)^k\det[P\,|\,Q]$, so the absolute values agree. The equality with $\Arg_{\pi,I}$ holds because $\T^{I}\to\T^{I}_\pi$ is a local isometry, and the last assertion follows from the Cauchy--Binet identity applied to both maps.
\end{proof}

Several statements below concern $k$-planes ``in general position'' or ``generic'' $k$-planes of $(\C^*)^{2k}$, and use the fact that their rolled coamoeba is covered once \cite[Proposition 8]{Mik17}. We make the genericity condition explicit. Let $V$ be the intersection with $(\C^*)^{2k}$ of an affine $k$-plane $\{z\in\C^{2k}:Lz=b\}$, with $L$ a complex $k\times2k$ matrix of rank $k$ and $b\in\C^k$. For $\varphi\in\T_\pi^{2k}$ write $D_\varphi=\operatorname{diag}(e^{i\varphi_1},\dots,e^{i\varphi_{2k}})$ and let $M(\varphi)$ be the real $2k\times2k$ matrix whose first $k$ rows are $\operatorname{Re}(LD_\varphi)$ and whose last $k$ rows are $\operatorname{Im}(LD_\varphi)$; the matrix $D_\varphi$ is defined up to the signs of its entries, which change $M(\varphi)$ by the multiplication of columns by $-1$ only. A point $z=D_\varphi r$ with $r\in(\R^*)^{2k}$ lies in $V$ if and only if $M(\varphi)r=(\operatorname{Re}b,\operatorname{Im}b)$. We say that $V$ satisfies condition (G) if $\det M(\varphi)$ is not identically zero as a function of $\varphi$ and if no coordinate of the solution $r(\varphi)=M(\varphi)^{-1}(\operatorname{Re}b,\operatorname{Im}b)$ vanishes identically.

\begin{lemma}\label{lem:G}
Condition (G) holds outside a proper real algebraic subset of the space of pairs $(L,b)$. If $V$ satisfies (G), then for almost every $\varphi\in\T_\pi^{2k}$ the fibre $V\cap\Arg_\pi^{-1}(\varphi)$ consists of exactly one point, and $\MVol(V)=\pi^{2k}$. If $V$ does not satisfy (G), then $\MVol(V)=0$. In particular $\MVol(V)\in\{0,\pi^{2k}\}$, and (G) holds if and only if $\Vol_{2k}(\A(V))>0$.
\end{lemma}

\begin{proof}
The function $\det M(\varphi)$ is a trigonometric polynomial in $\varphi$ whose coefficients are real polynomials in the real and imaginary parts of the entries of $L$, and by Cramer's rule $\det M(\varphi)\,r_j(\varphi)$ is a trigonometric polynomial whose coefficients are real polynomials in the entries of $(L,b)$. Hence each of the conditions $\det M\equiv0$ and $\det M\cdot r_j\equiv0$ is the vanishing of finitely many real polynomials in $(L,b)$, and (G) fails on a real algebraic subset. This subset is proper: for $L=(I_k\ \,iI_k)$ and $b_m=1+i$ the system at $\varphi=0$ reads $r_m=1$, $r_{k+m}=1$, so $\det M(0)\ne0$ and $r(0)=(1,\dots,1)$. If (G) holds, the zero sets of the nonzero real analytic functions $\det M$ and $\det M\cdot r_j$ on $\T_\pi^{2k}$ are null, and outside their union the system has exactly one solution, all of whose coordinates are nonzero. Hence $\MVol(V)=\pi^{2k}$ by \eqref{eq:rolled} below, applied with $I=\{1,\dots,2k\}$. If (G) fails, either $\det M\equiv0$, and then for every $\varphi$ the solution set of the system is empty or an affine subspace of $\R^{2k}$ of positive dimension, so that the fibre of $\Arg_\pi|_V$ over $\varphi$ is empty or infinite, hence empty for almost every $\varphi$ by Proposition \ref{prop:fibre}; or some $r_j$ vanishes identically, and then the fibre is empty for almost every $\varphi$. In both cases $\MVol(V)=0$ by \eqref{eq:rolled}. Finally, by the area formula $\MVol(V)=\int_{\R^{2k}}N(x)\,dx\ge\Vol_{2k}(\A(V))$, and $\MVol(V)>0$ forces $N>0$ on a set of positive measure, that is $\Vol_{2k}(\A(V))>0$.
\end{proof}

Throughout the paper, a generic $k$-plane of $(\C^*)^{2k}$ is one satisfying (G), equivalently one whose amoeba has positive volume, and when a $k$-plane $V$ of $(\C^*)^n$ with $n>2k$ is called generic, we require in addition that $\pi_I|_V$ be injective and that its image $\pi_I(V)$ satisfy (G) for every $I$ with $|I|=2k$; by Lemma \ref{lem:G} this excludes a proper real algebraic subset of the space of $k$-planes, and it gives $\MVol_I(V)=\pi^{2k}$ for every $I$. This is the content of \cite[Proposition 8]{Mik17} which we use.

\subsection{Further preliminaries}\label{sec:prelim}
For $z\in V$ and $|I|=2k$ we write $J(z)=J(\Log|_V)(z)$ and $J_I(z)=J(\Log_I|_V)(z)$. By \eqref{eq:CB0}, $J(z)^2=\sum_IJ_I(z)^2$, hence
\begin{equation}\label{eq:CBJ}J(z)\le\textstyle\sum_{|I|=2k}J_I(z),\end{equation}
with equality at $z$ if and only if at most one $J_I(z)$ is nonzero. The functions $J^2$ and $J_I^2$ are real analytic on $V$. In the real coordinates of $\C^n$ they are rational functions of the coordinates, because $T_zV=\ker df(z)$ and $d\Log=\re\,d\log$, so they are semialgebraic.

\begin{lemma}\label{lem:lsc}
The function $V\mapsto\MVol(V)$ is lower semicontinuous on $\Vc(\Delta)$.
\end{lemma}

\begin{proof}
Let $c_0\in\Vc(\Delta)$ with variety $V_0$ and let $K\subset V_0$ be compact. For $p\in K$ let $N_p$ be the Euclidean orthogonal complement of $T_pV_0$ in $\C^n$, and choose $\rho>0$ such that $(p,v)\mapsto p+v$, $v\in N_p$, $|v|<\rho$, is injective with values in $(\C^*)^n$. Since $df_{c_0}(p)$ has rank $n-k$ and kernel $T_pV_0$, its restriction to $N_p$ is invertible, and the implicit function theorem gives, for $c$ near $c_0$, a unique small $v(c,p)\in N_p$ with $f_c(p+v(c,p))=0$, smooth in $(c,p)$ and vanishing at $c_0$. Hence $\psi_c(p)=p+v(c,p)$ embeds $K$ into $V_c$ and tends to the identity in $C^1(K)$, so that $\MVol(V_c)\ge\int_KJ(\Log|_{V_c})(\psi_c(p))\,\mathrm{Jac}(\psi_c)(p)\,d\mathrm{vol}(p)\to\int_KJ\,d\mathrm{vol}$. Exhausting $V_0$ by compact sets gives $\liminf_{c\to c_0}\MVol(V_c)\ge\MVol(V_0)$.
\end{proof}

\begin{lemma}\label{lem:parity}
Let $W\subset(\C^*)^{2k}$ be a closed complex submanifold of dimension $k$ whose amoeba has finite Lebesgue measure. Then the counting function of $\Log|_W$ is even at every regular value, and at least $2$ at every regular value in $\A(W)$.
\end{lemma}

\begin{proof}
The map $\Log|_W:W\to\R^{2k}$ is proper, since $W$ is closed in $(\C^*)^{2k}$ and $\Log$ is proper, and it is a map between manifolds of the same dimension $2k$ without boundary. It therefore has a degree mod $2$, equal to the number mod $2$ of preimages of any regular value \cite[Section 4]{Mil65}; the proof given there for compact manifolds applies verbatim to proper maps, since the preimage of a small ball around a regular value is then a finite disjoint union of balls mapped diffeomorphically. An amoeba of finite measure is not all of $\R^{2k}$, so some regular value has no preimage, and the degree mod $2$ vanishes. At a regular value in $\A(W)$ the number of preimages is then even and positive.
\end{proof}

\section{Counting the fibres of the rolled coamoeba map}\label{sec:counting}

By Lemma \ref{lem:jac} and the area formula applied to $\Arg_{\pi,I}|_V$,

\begin{equation}\label{eq:rolled}
\MVol_I(V)=\int_{\T_\pi^{2k}}\#\big(V\cap\Arg_{\pi,I}^{-1}(\varphi)\big)\,d\varphi.
\end{equation}

This is the identity which in \cite[Corollary 5]{Mik17} is stated for $n=2k$; its advantage is that the right-hand side is an integral over a compact torus of volume $\pi^{2k}$, so that any uniform bound for the cardinality of the fibres of $\Arg_{\pi,I}|_V$ yields a bound for $\MVol_I(V)$.

Fix $\varphi\in\T_\pi^{I}$. A point $z\in(\C^*)^n$ lies in $\Arg_{\pi,I}^{-1}(\varphi)$ if and only if $z_i=e^{i\varphi_i}r_i$ with $r_i\in\R^*$ for $i\in I$, the coordinates $z_j$, $j\in\hat I$, being arbitrary. Hence the fibre $V\cap\Arg^{-1}_{\pi,I}(\varphi)$ is in bijection with the set of $(r_I,z_{\hat I})\in(\R^*)^{2k}\times(\C^*)^{n-2k}$ such that $F_m(r_I,z_{\hat I}):=f_m(e^{i\varphi_I}r_I,z_{\hat I})=0$ for $m=1,\dots,n-k$. This is a system of $2n-2k$ real equations in $2k+2(n-2k)=2n-2k$ real unknowns. To count its solutions we complexify it, as in \cite[Section 2.1]{Mik17}, treating $z_{\hat I}$ and its complex conjugate as independent variables. Let $\bar f_m$ be the polynomial obtained from $f_m$ by conjugating its coefficients and consider the square system in the variables $(r_I,z_{\hat I},w_{\hat I})\in(\C^*)^{2n-2k}$

\begin{equation}\label{eq:doubled}
F_m(r_I,z_{\hat I})=0,\qquad G_m(r_I,w_{\hat I}):=\bar f_m\big(e^{-i\varphi_I}r_I,\,w_{\hat I}\big)=0,\qquad m=1,\dots,n-k.
\end{equation}

If $(r_I,z_{\hat I})$ solves the real system with $r_I$ real, then $(r_I,z_{\hat I},\bar z_{\hat I})$ solves \eqref{eq:doubled}, because $G_m(r_I,\bar z_{\hat I})=\overline{F_m(r_I,z_{\hat I})}$. The Newton polytope of $F_m$ in the variables $(r_I,z_{\hat I},w_{\hat I})$ is the image of $\Delta_m$ under the embedding $\iota_1^I:\R^n\to\R^{2n-2k}$, $(a_I,a_{\hat I})\mapsto(a_I,a_{\hat I},0)$, and that of $G_m$ is the image of $\Delta_m$ under $\iota_2^I:(a_I,a_{\hat I})\mapsto(a_I,0,a_{\hat I})$; multiplying a variable by a nonzero constant does not change Newton polytopes.

\begin{definition}\label{def:alpha}
For $|I|=2k$ we put

$$\alpha_I(\Delta_1,\dots,\Delta_{n-k})=\MV\big(\iota_1^I\Delta_1,\dots,\iota_1^I\Delta_{n-k},\ \iota_2^I\Delta_1,\dots,\iota_2^I\Delta_{n-k}\big),$$

the mixed volume in $\R^{2n-2k}$ normalised as in Bernstein's theorem \cite{Ber75}, so that $\MV(P,\dots,P)=(2n-2k)!\Vol(P)$.
\end{definition}

When $n=2k$ we have $\hat I=\emptyset$, both embeddings are the identity, and $\alpha_I$ is Mikhalkin's conj-degree $\MV(\Delta_1,\dots,\Delta_k,\Delta_1,\dots,\Delta_k)$ of \cite[Proposition 2]{Mik17}. In general $\alpha_I$ is the number of points, for a generic multiplicative translation, of the intersection of $V\times\conj(V)\subset(\C^*)^n\times(\C^*)^n$ with the subtorus $\{(z,w):z_I=w_I\}$, so it is the toric intersection number of \cite[Definition 1.3]{Mik17} relative to the coordinates in $I$.

\begin{proposition}\label{prop:fibre}
For almost every $\varphi\in\T_\pi^{2k}$ the fibre $V\cap\Arg_{\pi,I}^{-1}(\varphi)$ is finite and has at most $\alpha_I(\Delta_1,\dots,\Delta_{n-k})$ points. Consequently $\MVol_I(V)\le \pi^{2k}\,\alpha_I(\Delta_1,\dots,\Delta_{n-k})$.
\end{proposition}

\begin{proof}
By Sard's theorem almost every $\varphi$ is a regular value of $\Arg_{\pi,I}|_V:V\to\T_\pi^{2k}$, a smooth map between manifolds of dimension $2k$. Fix such a $\varphi$ and put $E=\Arg^{-1}_{\pi,I}(\varphi)$, a submanifold of dimension $2n-2k$ which is a level set of the submersion $\Arg_{\pi,I}$. At $z\in V\cap E$ the kernel of $d(\Arg_{\pi,I}|_V)$ is $T_zV\cap T_zE$, so regularity means $T_zV\cap T_zE=0$. Let $\iota:(\R^*)^{2k}\times(\C^*)^{n-2k}\to E$ be the parametrisation $(r_I,z_{\hat I})\mapsto(e^{i\varphi_I}r_I,z_{\hat I})$ and $\Phi=(F_1,\dots,F_{n-k})=f\circ\iota$. Then $\ker d\Phi=d\iota^{-1}(\ker df)=d\iota^{-1}(T_zV\cap T_zE)=0$, since $\ker df=T_zV$ by the independence of the differentials of the $f_m$ along $V$. Source and target of $d\Phi$ have real dimension $2n-2k$, so $d\Phi$ is invertible at every point of the fibre.

The complex Jacobian matrix of \eqref{eq:doubled} at $(r_I,z_{\hat I},\bar z_{\hat I})$ is obtained from the real Jacobian matrix of $(\Phi,\bar\Phi)$ in the coordinates $(r_I,\operatorname{Re}z_{\hat I},\operatorname{Im}z_{\hat I})$ by the fixed invertible change of coordinates $(\operatorname{Re}z_{\hat I},\operatorname{Im}z_{\hat I})\leftrightarrow(z_{\hat I},\bar z_{\hat I})$; hence it is invertible, and every point of the fibre gives a nondegenerate, hence isolated, solution of \eqref{eq:doubled} in $(\C^*)^{2n-2k}$. Distinct points of the fibre give distinct solutions. By the theorem of Bernstein \cite{Ber75} (see also Kushnirenko \cite{Kus76} and Fulton \cite[Section 5.5]{Ful93}), for any coefficients the number of isolated solutions in the torus of a square system of Laurent polynomials, counted with multiplicity, does not exceed the mixed volume of their Newton polytopes. Hence the fibre has at most $\alpha_I$ points, and the bound for $\MVol_I(V)$ follows from \eqref{eq:rolled}.
\end{proof}

No genericity of the coefficients of the $f_m$ is needed in Proposition \ref{prop:fibre}: Bernstein's bound for the number of isolated solutions holds for arbitrary coefficients, and the nondegeneracy of the solutions coming from points of the fibre is a consequence of the regularity of $\varphi$, which holds for almost every $\varphi$ by Sard's theorem.

\begin{remark}
Working with $\Arg_I$ and $\T^{2k}=(\R/2\pi\Z)^{2k}$ instead of the rolled map would give the weaker bound $(2\pi)^{2k}\alpha_I$, because the fibres of $\Arg_I$ correspond to positive $r_I$ only. The rolled coamoeba, used systematically by Mikhalkin \cite{Mik17}, gains the factor $4^k$ and makes the estimate sharp for linear spaces (Sections \ref{sec:lines} and \ref{sec:kplanes}).
\end{remark}

\section{The main estimate}\label{sec:main}

\begin{theorem}\label{thm:main}
Let $V\subset(\C^*)^n$ be a smooth complete intersection of dimension $k$, $n\ge 2k$, with Newton polytopes $\Delta_1,\dots,\Delta_{n-k}\subset\R^n$. Then

$$\Hm^{2k}\big(\A(V)\big)\ \le\ \sum_{|I|=2k}\MVol_I(V)\ \le\ \pi^{2k}\sum_{|I|=2k}\alpha_I(\Delta_1,\dots,\Delta_{n-k}),$$

and $\Hm^{2k}(\A(V))\ge\big(\sum_{|I|=2k}\Vol_{2k}(\A(\pi_I(V)))^2\big)^{1/2}$.
\end{theorem}

\begin{proof}
Put $S=\A(V)$. Since $V$ is a countable union of compact pieces and $\Log$ is smooth, $S$ is countably $2k$-rectifiable and has an approximate tangent $2k$-plane $T_xS$ at $\Hm^{2k}$-almost every point \cite[Theorem 3.2.19]{Fed69}. Write $J_I^S(x)=J_I(T_xS)$. The area formula for rectifiable sets \cite[Theorem 3.2.22]{Fed69} gives $\int_S J^S_I\,d\Hm^{2k}=\int_{\R^{2k}}\#(S\cap p_I^{-1}(y))\,dy$. Since $\Log$ maps $V\cap\Log_I^{-1}(y)$ onto $S\cap p_I^{-1}(y)$, the integrand is bounded by $\#(V\cap\Log_I^{-1}(y))$, whose integral is $\MVol_I(V)$. By \eqref{eq:cb}, $1\le\sum_IJ^S_I$ almost everywhere on $S$, so $\Hm^{2k}(S)\le\sum_I\int_SJ^S_I\,d\Hm^{2k}\le\sum_I\MVol_I(V)$, and Proposition \ref{prop:fibre} gives the second inequality. For the lower bound, the vector $(J^S_I(x))_I$ has Euclidean norm one, hence, by the triangle inequality for integrals of vector valued functions, $\Hm^{2k}(S)=\int_S\|(J^S_I)_I\|\,d\Hm^{2k}\ge\|(\int_SJ^S_I\,d\Hm^{2k})_I\|$, and $\int_SJ^S_I\,d\Hm^{2k}\ge\Vol_{2k}(p_I(S))=\Vol_{2k}(\A(\pi_I(V)))$.
\end{proof}

\begin{corollary}\label{cor:degree}
If $\Delta_m\subset d_m\Sigma_n$ for every $m$, where $\Sigma_n$ is the standard simplex of $\R^n$, then $\Hm^{2k}(\A(V))\le\binom{n}{2k}\pi^{2k}(d_1\cdots d_{n-k})^2$.
\end{corollary}

\begin{proof}
The embeddings $\iota_1^I,\iota_2^I$ map $\Sigma_n$ into the standard simplex of $\R^{2n-2k}$, and mixed volume is monotone and multilinear, so $\alpha_I\le\MV(d_1\Sigma,\dots,d_{n-k}\Sigma,d_1\Sigma,\dots,d_{n-k}\Sigma)=(d_1\cdots d_{n-k})^2$ by B\'ezout's theorem; there are $\binom{n}{2k}$ subsets $I$.
\end{proof}

For $n=2k$ the corollary gives $\pi^{2k}(d_1\cdots d_k)^2$, twice Mikhalkin's bound \cite[Corollary 3]{Mik17}. The missing factor $\frac12$ comes, in the half-dimensional case, from the fact that $\Log|_V$ is proper and has degree zero \cite[Corollary 7]{Mik17}, so that almost every point of the amoeba is covered at least twice. When $n>2k$ the maps $\Log_I|_V$ are not proper in general, since a point of $V$ may escape to the toric boundary along the coordinates in $\hat I$ while $\Log_I$ stays bounded. Section \ref{sec:sym} shows that a factor $\frac12$ nevertheless holds for varieties which are real up to translation, and Section \ref{sec:lines} shows that it fails for nonreal lines.

\section{Symmetries of the amoeba map and the real case}\label{sec:sym}

\subsection{\texorpdfstring{The group $\G$ and the stabiliser of $V$}{The group G and the stabiliser of V}}
Let $\T=(S^1)^n\subset(\C^*)^n$ be the compact torus and let $\G$ be the group of the maps $z\mapsto\theta z$ and $z\mapsto\theta\bar z$, $\theta\in\T$, under composition; one checks that $(\theta\,\bar\cdot)\circ(\theta'\,\bar\cdot)$ is multiplication by $\theta\bar\theta'$, so that $\G$ is the semidirect product of $\T$ by the complex conjugation. Every $g\in\G$ satisfies $\Log\circ g=\Log$ and is an isometry of the flat metric. The elements of $\T$ act holomorphically and the elements of the coset $\T\cdot\conj$ antiholomorphically; each of the latter is an involution, since $\theta\overline{\theta\bar z}=\theta\bar\theta z=z$. The stabiliser of $V$ is $\G_V=\{g\in\G:g(V)=V\}$, and for a subgroup $\Gamma\subset\G_V$ we denote by $\Gamma^+$ its subgroup of elements which preserve the complex orientation of $V$. A holomorphic map preserves this orientation, and an antiholomorphic diffeomorphism of a complex manifold of dimension $k$ multiplies it by $(-1)^k$, since its differential is a linear isomorphism composed with the conjugation of $\C^k$, whose real determinant is $(-1)^k$. Hence $\Gamma^+=\Gamma$ if $k$ is even and $\Gamma^+=\Gamma\cap\T$ if $k$ is odd.

\begin{lemma}\label{lem:real}
The variety $V$ is real up to a multiplicative translation, that is $V=tW$ with $t\in(\C^*)^n$ and $\bar W=W$, if and only if $\G_V$ contains an antiholomorphic element.
\end{lemma}

\begin{proof}
If $V=tW$ with $\bar W=W$, then for $z=tw\in V$ we have $(t/\bar t)\bar z=t\bar w\in V$, so $z\mapsto\theta\bar z$ with $\theta=t/\bar t\in\T$ preserves $V$. Conversely, if $z\mapsto\theta\bar z$ preserves $V$, choose $\eta\in\T$ with $\eta^2=\theta$ and put $W=\eta^{-1}V$; for $w=\eta^{-1}z\in W$ we get $\bar w=\eta\bar z=\eta^{-1}(\theta\bar z)\in W$, because $\bar\eta=\eta^{-1}$.
\end{proof}

For instance, if each $f_m(tz)$ has real coefficients up to a nonzero constant factor, then $V$ is real up to translation; for the lines $L_c$ of Section \ref{sec:lines} this happens exactly when $c\in\R$.

\begin{lemma}\label{lem:fixed}
A holomorphic element $g\ne\mathrm{id}$ of $\G$ has no fixed point in $(\C^*)^n$. If $g\in\G_V$ is antiholomorphic, its fixed point set in $V$ is empty or a closed real submanifold of $V$ of dimension $k$.
\end{lemma}

\begin{proof}
If $\theta z=z$ for some $z\in(\C^*)^n$ then $\theta=1$. If $g$ is antiholomorphic and preserves $V$, it restricts to an isometric involution of $V$, so by Kobayashi's theorem \cite{Kob58} its fixed point set is a closed totally geodesic submanifold whose tangent space at a fixed point $p$ is the fixed space of $dg_p$, an antilinear involution of $T_pV\simeq\C^k$; this fixed space is totally real of dimension $k$.
\end{proof}

\begin{lemma}\label{lem:infinite}
If $\G_V$ is infinite, then $\Hm^{2k}(\A(V))=0$.
\end{lemma}

\begin{proof}
If $\G_V$ is infinite, so is the closed subgroup $\G_V\cap\T$ of $\T$, whose identity component is then a subtorus of positive dimension; thus $V$ is invariant under a one-parameter group $\tau\mapsto(e^{iu_1\tau},\dots,e^{iu_n\tau})$ with $u\in\R^n\setminus\{0\}$. The vector field $X(z)=(iu_1z_1,\dots,iu_nz_n)$ is tangent to $V$ and $d\Log(X)=0$, because its flow preserves the moduli $|z_j|$. Hence $d(\Log|_V)$ has rank at most $2k-1$ everywhere, $J(\Log|_V)\equiv0$, and the area formula gives $\int_{\R^n}\#(V\cap\Log^{-1}(x))\,d\Hm^{2k}(x)=0$; since the counting function is at least one on $\A(V)$, the volume vanishes.
\end{proof}

\subsection{\texorpdfstring{The factor $|\Gamma|^{-1}$}{The factor |Gamma|\textasciicircum (-1)}}

\begin{lemma}\label{lem:count}
Let $\Gamma\subset\G_V$ be a finite subgroup and let $F_\Gamma\subset V$ be the union of the fixed point sets in $V$ of the elements of $\Gamma\setminus\{\mathrm{id}\}$. Then $\Hm^{2k}(\Log(F_\Gamma))=0$, for every $x\in\A(V)\setminus\Log(F_\Gamma)$ the fibre $V\cap\Log^{-1}(x)$ has at least $|\Gamma|$ points, and for every $I$ the set $p_I(\Log(F_\Gamma))$ has Lebesgue measure zero.
\end{lemma}

\begin{proof}
By Lemma \ref{lem:fixed}, $F_\Gamma$ is a finite union of closed submanifolds of dimension $k$, so $F_\Gamma$ and its images under the locally Lipschitz maps $\Log$ and $\Log_I$ are $\sigma$-finite for $\Hm^k$; since $k<2k$ they have vanishing $\Hm^{2k}$-measure \cite[Section 2.10.2]{Fed69}. If $x\in\A(V)\setminus\Log(F_\Gamma)$ and $\Log z=x$ with $z\in V$, then $z\notin F_\Gamma$, so the points $gz$, $g\in\Gamma$, are pairwise distinct, lie in $V$, and lie in $\Log^{-1}(x)$.
\end{proof}

\begin{theorem}\label{thm:sym}
Let $V$ be as in Theorem \ref{thm:main} and let $\Gamma\subset\G_V$ be a finite subgroup. Then

$$\Hm^{2k}(\A(V))\le\frac1{|\Gamma|}\MVol(V)\le\frac1{|\Gamma|}\sum_{|I|=2k}\MVol_I(V)\le\frac{\pi^{2k}}{|\Gamma|}\sum_{|I|=2k}\alpha_I(\Delta_1,\dots,\Delta_{n-k}).$$

Moreover $\int_{\A(V)}J^{\A(V)}_I\,d\Hm^{2k}\le|\Gamma|^{-1}\MVol_I(V)$ for each $I$.
\end{theorem}

\begin{proof}
Put $S=\A(V)$ and $N(x)=\#(V\cap\Log^{-1}(x))$. By Lemma \ref{lem:count}, $N\ge|\Gamma|$ almost everywhere on $S$, hence $\Hm^{2k}(S)\le|\Gamma|^{-1}\int_SN\,d\Hm^{2k}=|\Gamma|^{-1}\MVol(V)$. The second inequality is \eqref{eq:mvsum} and the third is Proposition \ref{prop:fibre}. For the projected statement, for $y$ outside the null set $p_I(\Log(F_\Gamma))$ every point of $S\cap p_I^{-1}(y)$ has at least $|\Gamma|$ preimages in $V\cap\Log_I^{-1}(y)$, preimages of distinct points being distinct, so $\#(S\cap p_I^{-1}(y))\le|\Gamma|^{-1}\#(V\cap\Log_I^{-1}(y))$ almost everywhere; integrating and using the area formula as in the proof of Theorem \ref{thm:main} gives the claim.
\end{proof}

\begin{corollary}\label{cor:half}
If $V$ is real up to a multiplicative translation and $n\ge2k$, then $\Hm^{2k}(\A(V))\le\frac12\MVol(V)\le\frac{\pi^{2k}}2\sum_{|I|=2k}\alpha_I$, and if moreover $\Delta_m\subset d_m\Sigma_n$ for all $m$, then $\Hm^{2k}(\A(V))\le$ $\frac12\binom n{2k}\pi^{2k}(d_1\cdots d_{n-k})^2$.
\end{corollary}

\begin{proof}
Apply Theorem \ref{thm:sym} to $\Gamma=\{\mathrm{id},\,z\mapsto\theta\bar z\}$, given by Lemma \ref{lem:real}, and argue as in Corollary \ref{cor:degree}.
\end{proof}

The factor $\frac12$ of Corollary \ref{cor:half} and the factor $\frac12$ of Mikhalkin's theorem have different origins. Mikhalkin's factor is topological: it comes from the vanishing of the degree of the proper map $\Log|_V$ when $n=2k$, and it holds for every $V$. Ours comes from the free action of the conjugation on the fibres of $\Log|_V$; it holds for every $n\ge2k$, but only for varieties that are real up to translation. Theorem \ref{thm:line} below shows that for $n>2k$ no factor $\frac12$ can hold without such a hypothesis.

\begin{example}\label{ex:power}
Let $W\subset(\C^*)^n$ be a smooth complete intersection of dimension $k$, $n\ge2k$, and $V=\phi_m^{-1}(W)$ its preimage under the power map $\phi_m(z)=(z_1^m,\dots,z_n^m)$, an unramified covering of degree $m^n$ of the torus. Then $V$ is a smooth complete intersection invariant under the group $\mu_m^n\subset\T$ of order $m^n$. Since $\Log\circ\phi_m=m\Log$, we have $\A(V)=\frac1m\A(W)$ and $\Hm^{2k}(\A(V))=m^{-2k}\Hm^{2k}(\A(W))$, while each fibre of $\Log|_V$ has $m^n$ times as many points as the corresponding fibre of $\Log|_W$, so $\MVol(V)=m^{n-2k}\MVol(W)$. Theorem \ref{thm:sym} with $\Gamma=\mu_m^n$ reads $m^{-2k}\Hm^{2k}(\A(W))\le m^{-2k}\MVol(W)$, that is the inequality $\Hm^{2k}(\A(W))\le\MVol(W)$ without loss, whereas without the symmetry the factor $m^n$ would be lost. The last bound of Theorem \ref{thm:sym} remains far from sharp here, since $\alpha_I(V)=m^{2n-2k}\alpha_I(W)$.
\end{example}

\subsection{The half-dimensional case: divisibility of the fibres}\label{sec:half}

In the half-dimensional case $n=2k$ the amoeba map $\Log|_V:V\to\R^{2k}$ is proper, because $\Log$ is proper on $(\C^*)^{2k}$ and $V$ is closed, and Mikhalkin proved that it has degree zero \cite[Corollary 7]{Mik17}, so that its fibres over regular values have an even number of points. We now show that the symmetries of $V$ multiply this parity constraint. We write $N(x)=\#(V\cap\Log^{-1}(x))$, and for a finite subgroup $\Gamma\subset\G_V$ we let $F_\Gamma\subset V$ be, as in Lemma \ref{lem:count}, the union of the fixed point sets in $V$ of the elements of $\Gamma\setminus\{\mathrm{id}\}$. By Lemma \ref{lem:fixed}, $F_\Gamma$ is a finite union of closed real submanifolds of $V$ of dimension $k$, and it is empty when $\Gamma$ contains no antiholomorphic element.

We use two topological facts. The first is the mod $2$ degree of proper maps: if $f:M\to Y$ is a proper smooth map between smooth manifolds without boundary of the same dimension and $Y$ is connected, the number of points of $f^{-1}(y)$ modulo $2$ is the same for all regular values $y$, and it vanishes if $f$ is not surjective. For compact $M$ this is \cite[Chapter 4]{Mil65}; the proof for proper maps is the same, since it only uses the compactness of the preimage of an arc joining two regular values, and it is given in \cite[Chapter 5]{Hir76} and \cite[Chapter 2, Section 4]{GP74}. No orientability is required. The second is that a closed subset of $\R^m$ of topological dimension at most $m-2$ does not separate $\R^m$ \cite[Chapter IV]{HW41}, together with the fact that the topological dimension of a subset of $\R^m$ does not exceed its Hausdorff dimension \cite[Chapter VII]{HW41}.

\begin{lemma}\label{lem:K}
Let $W\subset(\C^*)^{2k}$ be a closed subset, let $E\subset W$ be a closed subset which is a countable union of compact subsets of $C^1$ submanifolds of dimension at most $d$, and assume that $\Log|_W$ is proper. Then $\Log(E)$ is closed, has topological dimension at most $d$ and Lebesgue measure zero if $d<2k$, and its complement in $\R^{2k}$ is connected if $d\le2k-2$.
\end{lemma}

\begin{proof}
A proper continuous map into $\R^{2k}$ is closed, so $\Log(E)$ is closed. Since $\Log$ is locally Lipschitz, $\Log(E)$ is $\sigma$-finite for $\Hm^d$ \cite[Section 2.10.2]{Fed69}, hence of Hausdorff dimension at most $d$, of topological dimension at most $d$ \cite[Chapter VII]{HW41}, and of Lebesgue measure zero if $d<2k$. The last assertion is the non-separation theorem quoted above.
\end{proof}

\begin{theorem}\label{thm:half}
Let $V\subset(\C^*)^{2k}$ be a smooth algebraic variety of dimension $k$ and let $\Gamma\subset\G_V$ be a finite subgroup. Then the following hold.

(a) For almost every $x\in\A(V)$ one has $N(x)\ge2|\Gamma^+|$.

(b) If $k\ge2$, or if $\Gamma$ contains no antiholomorphic element, then $N(x)$ is divisible by $2|\Gamma|$ for almost every $x\in\R^{2k}$; in particular $N(x)\ge2|\Gamma|$ for almost every $x\in\A(V)$.

Consequently $\Vol(\A(V))\le\frac1{2m}\MVol(V)\le\frac{\pi^{2k}}{2m}\alpha(V)$, where $m=|\Gamma|$ in case (b) and $m=|\Gamma^+|$ otherwise.
\end{theorem}

\begin{proof}
For (a), let $x\in\A(V)$ be a regular value outside the null set $\Log(F_\Gamma)$ of Lemma \ref{lem:count}; almost every point of $\A(V)$ is such by Sard's theorem. Each point $z$ of the fibre has a sign $\epsilon(z)=\pm1$ according to whether $d(\Log|_V)_z$ preserves orientation, and since the degree vanishes both signs occur equally often. For $g\in\Gamma$ the identity $\Log\circ g=\Log$ gives $d(\Log|_V)_{gz}\circ dg_z=d(\Log|_V)_z$, so $\epsilon(gz)=\epsilon(z)$ when $g\in\Gamma^+$ and $\epsilon(gz)=-\epsilon(z)$ otherwise. Choosing $z$ in the fibre, the $|\Gamma^+|$ distinct points $gz$, $g\in\Gamma^+$, have the sign of $z$, and the fibre contains at least as many points of the opposite sign.

For (b), put $K_\Gamma=\Log(F_\Gamma)$ and $\Omega_\Gamma=\R^{2k}\setminus K_\Gamma$. By Lemma \ref{lem:K} with $W=V$, $E=F_\Gamma$ and $d=k$, the set $K_\Gamma$ is closed and of measure zero, and $\Omega_\Gamma$ is connected, because $k\le2k-2$ when $k\ge2$, and $K_\Gamma=\emptyset$ when $\Gamma$ has no antiholomorphic element. The open set $V^\circ=(\Log|_V)^{-1}(\Omega_\Gamma)$ is $\Gamma$-invariant, since $\Log\circ g=\Log$, and $\Gamma$ acts on it freely by diffeomorphisms, since $V^\circ\cap F_\Gamma=\emptyset$. A free action of a finite group on a manifold is properly discontinuous, so $M=V^\circ/\Gamma$ is a smooth manifold of dimension $2k$ and $V^\circ\to M$ is a covering of degree $|\Gamma|$ \cite[Chapter 21]{Lee13}. The map $\Log$ induces a smooth map $\lambda:M\to\Omega_\Gamma$, which is proper: $\Log|_{V^\circ}:V^\circ\to\Omega_\Gamma$ is proper, being the restriction of a proper map to a full preimage, and the preimage under $\lambda$ of a compact set is the image of a compact set under the projection. The regular values of $\lambda$ are the regular values of $\Log|_V$ in $\Omega_\Gamma$. The amoeba is closed with nonempty complement, since its volume is finite \cite{MN13}, and $K_\Gamma\subset\A(V)$; hence $\Omega_\Gamma$ contains regular values with empty preimage, and the mod $2$ degree of $\lambda$ vanishes. Thus $\#\lambda^{-1}(x)$ is even for every regular value $x\in\Omega_\Gamma$, and $N(x)=|\Gamma|\cdot\#\lambda^{-1}(x)$, since $V\cap\Log^{-1}(x)\subset V^\circ$ is a union of free $\Gamma$-orbits. Almost every $x$ is a regular value lying in $\Omega_\Gamma$, which proves (b).

The volume bound follows from the identity $\Vol(\A(V))=\int_VJ(\Log|_V)\,N(\Log z)^{-1}\,d\mathrm{vol}(z)$, a consequence of the area formula \cite[Theorem 3.2.3]{Fed69}, together with $\MVol(V)=\int_VJ(\Log|_V)\,d\mathrm{vol}$ and Mikhalkin's bound $\MVol(V)\le\pi^{2k}\alpha(V)$ \cite[Lemma 6]{Mik17}.
\end{proof}

In case (b) the bound $2|\Gamma|$ is at least the bound $2|\Gamma^+|$ of case (a), with strict inequality exactly when $k$ is odd and $\Gamma$ contains an antiholomorphic element; in general $2|\Gamma^+|\ge|\Gamma|$, so the theorem contains the case $n=2k$ of Theorem \ref{thm:sym}. Notice that the degree of $\Log|_V$ itself is not used in (b): the mod $2$ degree of the quotient map contains this information and more.

\begin{corollary}\label{cor:real4}
Let $V\subset(\C^*)^{2k}$ be a smooth algebraic variety of dimension $k\ge2$ which is real up to a multiplicative translation. Then $N(x)\equiv0\pmod 4$ for almost every $x$, and $\Vol(\A(V))\le\frac14\MVol(V)\le\frac{\pi^{2k}}4\alpha(V)$.
\end{corollary}

\begin{proof}
Apply Theorem \ref{thm:half}(b) to $\Gamma=\{\mathrm{id},\sigma\}$, where $\sigma$ is the antiholomorphic involution of Lemma \ref{lem:real}.
\end{proof}

\begin{remark}\label{rem:k1}
The hypothesis $k\ge2$ in Corollary \ref{cor:real4} cannot be removed: a simple Harnack curve is real and its amoeba map is two-to-one \cite{Mik00,MR01}. In the proof, the image $\Log(V(\R))$ of the real locus of a real curve is a curve, which may separate the plane; for a simple Harnack curve it is the boundary of the amoeba \cite{Mik00}. Geometrically, the two points of a conjugate pair can only collide at a real point; for $k=1$ such collisions form walls of the discriminant of the amoeba map, across which the fibre changes by $2$, while for $k\ge2$ they happen over a set of codimension at least $2$, which a generic path avoids, and the walls it crosses carry critical points in conjugate pairs, across which the fibre changes by a multiple of $4$. Only the first sentence of this remark and the divisibility statement of Corollary \ref{cor:real4} are proved; the description of the walls is heuristic, and its rigorous substitute is the description of the fibres near the tentacles of the amoeba in Lemma \ref{lem:tentacle} and Corollary \ref{cor:fibre4}.
\end{remark}

Mikhalkin's argument also bounds the fibres from above \cite[Section 2.1]{Mik17}: the fibre over $x$ is the fixed point set in $V$ of the antiholomorphic involution $\conj'_x(z)=e^{2x}/\bar z$, so it is contained in $V\cap\conj'_x(V)$, which for a regular value has at most $\beta(V)$ points, the conj$'$-degree; for a generic $k$-plane $\beta(V)=\binom{2k}k$ \cite[Corollary 3]{Mik17}.

\begin{corollary}\label{cor:sizes}
Let $V\subset(\C^*)^{2k}$ be a real generic $k$-plane. For $k=1$ almost every fibre over $\A(V)$ has exactly $2$ points and $\Vol(\A(V))=\frac{\pi^2}2$; for $k=2$ almost every fibre over $\A(V)$ has exactly $4$ points and $\Vol(\A(V))=\frac{\pi^4}4$; for $k=3$ almost every fibre over $\A(V)$ has $4$, $8$, $12$, $16$ or $20$ points; for $k=4$ the number of points is a multiple of $4$ not exceeding $68$.
\end{corollary}

\begin{proof}
For $k=1$, $N$ is even, positive and at most $\binom21=2$. For $k\ge2$, $N$ is a positive multiple of $4$ by Corollary \ref{cor:real4} and at most $\binom{2k}k$, which equals $6$, $20$, $70$ for $k=2,3,4$. The volumes follow from $\MVol(V)=\pi^{2k}$ \cite[Proposition 8]{Mik17}.
\end{proof}

\begin{example}\label{ex:planes}
For a generic real plane $V\subset(\C^*)^4$ we have $\alpha(V)=1$ and $\MVol(V)=\pi^4$ by Lemma \ref{lem:G} \cite[Proposition 8]{Mik17}, and Corollary \ref{cor:sizes} gives $\Vol(\A(V))=\pi^4/4$, the value found by Nisse and Passare \cite[Theorem 1.2]{NP17}, with a proof which only uses the conjugation, the topology of the complement of $\Log(V(\R))$ and Mikhalkin's upper bound. Thus Corollary \ref{cor:real4} is an equality for all such planes, the amoeba map being four-to-one almost everywhere, and Mikhalkin's bound $\frac{\pi^4}2$ is attained by no real plane. For $k\ge3$ the statement of \cite{NP17} that the amoeba map of a real generic $k$-plane is $2^k$-to-one does not hold: Section \ref{sec:why} explains why, and Corollary \ref{cor:fibre4} proves it. For $k=1$, Theorem \ref{thm:half} reduces to the Passare--Rullg{\aa}rd bound $\Vol(\A(V))\le\frac{\pi^2}2\alpha(V)=\pi^2\Vol(\Delta)$ \cite{PR04}, attained by the simple Harnack curves \cite{MR01,Mik00}.
\end{example}

\subsection{\texorpdfstring{Partial conjugations and the count $2^k$}{Partial conjugations and the count 2\textasciicircum k}}\label{sec:why}

Corollary \ref{cor:sizes} shows that for real generic $k$-planes the value $2^k$ is forced for $k\le2$ by arithmetic alone: $2$ is the only even positive integer not exceeding $\binom21$, and $4$ is the only positive multiple of $4$ not exceeding $\binom42=6$. For $k=3$ the multiples of $4$ not exceeding $\binom63=20$ are $4,8,12,16,20$, and Corollary \ref{cor:fibre4} shows that fibres of $4$ points occur on a set of positive measure for every real $3$-plane in general position which is not a symmetric cube of a line. We now explain where the value $2^k$ comes from and why it does not persist.

For $S\subset\{1,\dots,2k\}$ let $\tau_S$ be the partial conjugation $\tau_S(z)_j=\bar z_j$ for $j\in S$ and $\tau_S(z)_j=z_j$ otherwise. Every $\tau_S$ satisfies $\Log\circ\tau_S=\Log$, but for $\emptyset\ne S\ne\{1,\dots,2k\}$ it is neither holomorphic nor antiholomorphic and does not belong to $\G$.

\begin{lemma}\label{lem:split}
Let $W\subset\C^S\oplus\C^{S^c}$ be a complex linear subspace and let $\tau(w_S,w_{S^c})=(\bar w_S,w_{S^c})$. Then $\tau(W)$ is a complex subspace if and only if $W=(W\cap\C^S)\oplus(W\cap\C^{S^c})$.
\end{lemma}

\begin{proof}
If $W$ splits, $\tau(W)=\overline{W\cap\C^S}\oplus(W\cap\C^{S^c})$ is complex. Conversely, if $\tau(W)$ is complex and $w\in W$, then $\tau(iw)=(-i\bar w_S,iw_{S^c})$ and $i\tau(w)=(i\bar w_S,iw_{S^c})$ lie in $\tau(W)$, hence so does $(\bar w_S,0)$, and therefore $(w_S,0)\in W$; thus $W$ contains the projections of its elements to $\C^S$ and splits.
\end{proof}

If a complex submanifold $V$ is invariant under $\tau_S$, then $d\tau_S$ maps $T_pV$ onto the complex subspace $T_{\tau_S(p)}V$, so all tangent spaces of $V$ split along $\C^S\oplus\C^{S^c}$; for a linear space this means that $V$ is a product of linear spaces in $(\C^*)^S$ and $(\C^*)^{S^c}$. A generic real $k$-plane is a product for no nontrivial splitting, so its only fibre-preserving symmetries among the partial conjugations are the identity and the total conjugation, and Corollary \ref{cor:real4} is all that the real structure can give.

A product $\mathscr L=L_1\times\dots\times L_k$ of real generic lines, on the other hand, is invariant under the group $P\simeq(\Z/2)^k$ of the partial conjugations $\tau_{\{2i-1,2i\}}$ of the factors, and its fibres are products of the fibres of the factors, which have exactly $2$ points; thus $N=2^k$ almost everywhere. Here the fixed point sets of the elements of $P$ have real codimension one, their images being the boundary walls $\A(L_1)\times\dots\times\partial\A(L_i)\times\dots\times\A(L_k)$ of the product amoeba, so the argument of Theorem \ref{thm:half}(b) does not apply to $P$, which is why the answer is $2^k$ rather than $2^{k+1}$. The symmetric powers of real lines of Section \ref{sec:squares} form another family with fibres of exactly $2^k$ points, produced by the same partial conjugations acting on the factors of a symmetric product.

The proof of \cite[Proposition 5.6]{NP17} compares a real generic $k$-plane $\mathscr P$, through a real linear change of the parameters, with a product $\mathscr L$ of $k$ lines, and transfers the count $2^k$ from $\mathscr L$ to $\mathscr P$; the lower bound is attributed to a result of Madani and Nisse \cite[Proposition 5.2]{NP17}, according to which the fibres have at least $2^l$ points with $l=\dim_\C(\mathscr P\cap\overline{\mathscr P})$, and the upper bound to the injectivity of the argument map on regular points \cite[Proposition 5.4]{NP17}. A linear change of the parameters maps $\mathscr P$ to $\mathscr L$ but does not intertwine the amoeba maps, since $\Log$ is equivariant only under monomial maps, so the fibres of $\Log|_{\mathscr P}$ and of $\Log|_{\mathscr L}$ need not have the same cardinality; the count $2^k$ is a property of $\mathscr L$ coming from its symmetry group $P$, which $\mathscr P$ does not have by Lemma \ref{lem:split}. The dimension of $\mathscr P\cap\overline{\mathscr P}$ equals $k$ for every real $\mathscr P$, products included, so it cannot distinguish the product, where fibres have $2^k$ points, from the generic case, where fibres of $4$ points occur for $k=3,4$.

The geometry of a deformation from $\mathscr L$ to a generic real $k$-plane suggests how fibres of $4$ points appear. The fibre count is locally constant off the discriminant, the image of the critical locus of the amoeba map. For the product, each point of a boundary wall carries $2^{k-1}$ critical points, a real point of $L_i$ times the $2^{k-1}$ points of the fibres of the other factors, forming $2^{k-2}$ conjugate pairs when $k\ge2$, and crossing the wall changes $N$ by $2^k$; this is a highly degenerate situation in which many folds of the amoeba map lie over one hypersurface. Under a generic real deformation one expects these folds to separate into distinct walls, except that conjugate folds stay over the same wall since $\Log\circ\conj=\Log$; crossing a generic wall then changes $N$ by exactly $4$, in accordance with Corollary \ref{cor:real4}, and regions with $N=4,8,\dots$ appear between the separated walls. For $k\le2$ at most one conjugate pair of folds lies over a point of a wall, so nothing separates and the count $2^k$ persists. This description is heuristic, but a rigorous version of it, in which the walls are replaced by the tentacles of the amoeba, is given in Section \ref{sec:products}: by Corollary \ref{cor:fibre4}, every real $3$-plane in general position which is not a symmetric cube has fibres of exactly $4$ points on a set of positive measure, and for every $k\ge3$ every real $k$-plane in general position which is not a symmetric power of a line has fibres of fewer than $2^k$ points on a set of positive measure.

\subsection{\texorpdfstring{Divisibility for the projections when $n>2k$}{Divisibility for the projections when n>2k}}\label{sec:proj}

Let now $n>2k$, let $|I|=2k$, and let $V\subset(\C^*)^n$ be a smooth algebraic variety of dimension $k$. The projection $\Log_I|_V$ is not proper in general, so the argument of Theorem \ref{thm:half} does not apply to it directly. We pass through the image of $V$. Let $V_I$ be the Zariski closure of $\pi_I(V)$ in $(\C^*)^{2k}$ and assume that $\pi_I|_V:V\to V_I$ is generically finite, that is $\dim V_I=k$; this holds for generic $V$, and in particular for generic linear spaces. Since we work in characteristic zero, there is a dense Zariski open subset $U\subset V_I$, contained in the smooth locus of $V_I$, over which $\pi_I$ restricts to a finite \'etale map $V\cap\pi_I^{-1}(U)\to U$ of some degree $d_I$; this follows from generic finiteness \cite[Exercise II.3.7]{Har77} and generic smoothness \cite[Corollary III.10.7]{Har77}. Put $B_I=V_I\setminus U$, a closed algebraic subset of complex dimension at most $k-1$, and for $y\in\R^{2k}$ write $N_I(y)=\#(V\cap\Log_I^{-1}(y)),\quad S_I(y)=\#(\A(V)\cap p_I^{-1}(y)),$ the numbers of points of $V$ and of the amoeba of $V$ lying over $y$. By the area formula, $\MVol_I(V)=\int_{\R^{2k}}N_I(y)\,dy$, and, as in the proof of Theorem \ref{thm:main}, $\int_{\A(V)}J^{\A(V)}_I\,d\Hm^{2k}=\int_{\R^{2k}}S_I(y)\,dy$, while $\Vol_{2k}(p_I(\A(V)))=\Vol_{2k}(\A(V_I))$ is the measure of the set where $S_I>0$.

\begin{theorem}\label{thm:proj}
Let $V$, $I$, $V_I$ and $d_I$ be as above and let $\Gamma\subset\G_V$ be a finite subgroup, where $\G$ is now the group of maps $z\mapsto\theta z$, $z\mapsto\theta\bar z$ of $(\C^*)^n$. Assume that $k\ge2$, or that $\Gamma$ contains no antiholomorphic element. Then for almost every $y\in\R^{2k}$ the number $N_I(y)$ is divisible by $2|\Gamma|$, and it equals $d_I$ times the number of points of $V_I\cap\Log^{-1}(y)$. Consequently:

(a) the projected amoeba satisfies $\Vol_{2k}(\A(V_I))\le\frac1{2|\Gamma|}\MVol_I(V)$;

(b) the projected part of the volume satisfies $\Vol_{2k}(\A(V_I))\le\int_{\A(V)}J^{\A(V)}_I\,d\Hm^{2k}\le\frac1{|\Gamma|}\MVol_I(V)$;

(c) equality holds on the right of (b) for every $I$ if and only if $N(x)=|\Gamma|$ for $\Hm^{2k}$-almost every $x\in\A(V)$, that is, if and only if almost every fibre of the amoeba map is a single $\Gamma$-orbit; in that case $S_I(y)=N_I(y)/|\Gamma|$ is even for almost every $y$, so that $\int_{\A(V)}J^{\A(V)}_I\,d\Hm^{2k}\ge2\Vol_{2k}(\A(V_I))$.
\end{theorem}

\begin{proof}
Every $g\in\G$ acts coordinatewise, so $\pi_I\circ g=g_I\circ\pi_I$, where $g_I$ is the corresponding map of $(\C^*)^{2k}$, and $\Log_I\circ g=\Log_I$. Let $F_\Gamma\subset V$ be the union of the fixed point sets of the elements of $\Gamma\setminus\{\mathrm{id}\}$ and let $\overline{\pi_I(F_\Gamma)}$ be the closure of its image in $V_I$. The set $F_\Gamma$ is a real algebraic subset of $V$ of dimension at most $k$, so $\pi_I(F_\Gamma)$ and its closure are semialgebraic of dimension at most $k$, hence countable unions of compact subsets of $C^1$ submanifolds of dimension at most $k$. The map $\Log|_{V_I}$ is proper, since $V_I$ is closed in $(\C^*)^{2k}$. Put $K=\Log(B_I)\cup\Log\big(\overline{\pi_I(F_\Gamma)}\big),\quad\Omega=\R^{2k}\setminus K$. By Lemma \ref{lem:K} applied to $W=V_I$, the set $\Log(B_I)$ is closed of topological dimension at most $2k-2$, since $B_I$ has real dimension at most $2k-2$, and $\Log(\overline{\pi_I(F_\Gamma)})$ is closed of topological dimension at most $k$. The union of two closed sets of dimension at most $2k-2$ has dimension at most $2k-2$ \cite[Chapter III]{HW41}, which holds for the second set when $k\ge2$, and which is automatic when $F_\Gamma=\emptyset$. Hence $K$ is closed, of measure zero, and $\Omega$ is connected.

Let $V^\circ=V\cap\Log_I^{-1}(\Omega)$ and $U_\Omega=V_I\cap\Log^{-1}(\Omega)$. If $z\in V^\circ$, then $\pi_I(z)\in V_I$ and $\Log(\pi_I(z))\in\Omega$, hence $\pi_I(z)\notin B_I$, that is $\pi_I(z)\in U$; therefore $V^\circ=V\cap\pi_I^{-1}(U_\Omega)$ and $\pi_I:V^\circ\to U_\Omega$ is a finite \'etale covering of degree $d_I$, in particular a proper map. Since $U_\Omega$ is the full preimage of $\Omega$ under the proper map $\Log|_{V_I}$, the map $\Log:U_\Omega\to\Omega$ is proper, and so is the composition $\Log_I|_{V^\circ}:V^\circ\to\Omega$. It follows that, for every regular value $y\in\Omega$, $N_I(y)=d_I\cdot\#(V_I\cap\Log^{-1}(y))$.

The group $\Gamma$ preserves $V^\circ$, and it acts freely on it, because $V^\circ\cap F_\Gamma=\emptyset$: indeed, if $z\in F_\Gamma$ then $\Log_I(z)\in\Log(\pi_I(F_\Gamma))\subset K$. As in the proof of Theorem \ref{thm:half}, $M=V^\circ/\Gamma$ is a smooth manifold of dimension $2k$ and $\Log_I$ induces a proper smooth map $\lambda:M\to\Omega$, whose regular values are the regular values of $\Log_I|_V$ in $\Omega$. The amoeba $\A(V_I)=p_I(\A(V))$ is closed and has finite volume \cite{MN13}, so its complement is a nonempty open set, which meets $\Omega$ because $K$ has measure zero; hence the mod $2$ degree of $\lambda$ vanishes, $\#\lambda^{-1}(y)$ is even for every regular value $y\in\Omega$, and $N_I(y)=|\Gamma|\cdot\#\lambda^{-1}(y)$ is divisible by $2|\Gamma|$. By Sard's theorem and since $K$ has measure zero, this holds for almost every $y$.

For (a), $N_I\ge2|\Gamma|$ almost everywhere on $\A(V_I)$, so $\MVol_I(V)=\int N_I\,dy\ge2|\Gamma|\Vol_{2k}(\A(V_I))$. For (b), the left inequality holds because $S_I\ge1$ on $\A(V_I)$, and the right one is the projected statement of Theorem \ref{thm:sym}: the points of $V$ over $y$ are grouped according to the point $x\in\A(V)$ over $y$ which is their image under $\Log$, each group is the fibre $V\cap\Log^{-1}(x)$, and for almost every $y$ every such fibre is a union of free $\Gamma$-orbits by Lemma \ref{lem:count}; hence

\begin{equation}\label{eq:NS}
N_I(y)=\sum_{x\in\A(V)\cap p_I^{-1}(y)}N(x)\ \ge\ |\Gamma|\,S_I(y).
\end{equation}

For (c), equality on the right of (b) for a given $I$ means that equality holds in \eqref{eq:NS} for almost every $y$, that is $N(x)=|\Gamma|$ for every $x$ over almost every $y$. If $N=|\Gamma|$ almost everywhere on $\A(V)$ for $\Hm^{2k}$, then by the area formula for rectifiable sets \cite[Theorem 3.2.22]{Fed69} the set of $y$ over which some $x$ has $N(x)>|\Gamma|$ has measure zero, and equality holds for every $I$. Conversely, if equality holds for every $I$, let $E=\{x\in\A(V):N(x)>|\Gamma|\}$; then $\int_EJ^{\A(V)}_I\,d\Hm^{2k}=0$ for every $I$, and since $\sum_IJ_I^{\A(V)}\ge1$ by \eqref{eq:cb}, $\Hm^{2k}(E)=0$. In that case $S_I=N_I/|\Gamma|$ almost everywhere, which is even by the first part, so $S_I\ge2$ almost everywhere on $\A(V_I)$, and integrating gives the last inequality.
\end{proof}

For $\Gamma$ trivial the theorem says that $N_I$ is even almost everywhere for every smooth $V$ with $\pi_I|_V$ generically finite, which is the projected analogue of Mikhalkin's degree zero theorem; it holds also for $k=1$. For real $V$ with $k\ge2$ it says that $N_I\equiv0\pmod4$ while only $S_I\le N_I/2$ is guaranteed. Real planes in $(\C^*)^5$ show that this is the correct form: for almost every such plane $N_I=4$ and $S_I=2$ (Corollary \ref{cor:plane5} and Theorem \ref{thm:twotoone}), so the projected part of the volume is $\frac12\MVol_I(V)$ and not $\frac14\MVol_I(V)$, while for the squares of lines of Section \ref{sec:squares} one has $N_I=4$ and $S_I=1$.

\begin{corollary}\label{cor:plane5}
Let $V\subset(\C^*)^5$ be a real generic plane. Then, for every $I$ with $|I|=4$, almost every fibre of $\Log_I|_V$ over $\A(V_I)$ has exactly $4$ points, $\Vol_4(\A(V_I))=\frac{\pi^4}4$, and $\frac{\pi^4}4\le\int_{\A(V)}J^{\A(V)}_I\,d\Hm^4\le\frac{\pi^4}2$. The right inequality is an equality for every $I$ if and only if the amoeba map $\Log|_V$ is two-to-one almost everywhere on $\A(V)$, and in that case $\Hm^4(\A(V))=\frac12\MVol(V)$.
\end{corollary}

\begin{proof}
For generic $V$ the map $\pi_I|_V$ is injective, since it is the restriction of a linear map to a plane not containing the direction of its kernel, and $V_I$ is a real generic plane of $(\C^*)^4$, so $d_I=1$, $\MVol_I(V)=\pi^4$ by \cite[Proposition 8]{Mik17}, and the fibres of $\Log$ on $V_I$ have exactly $4$ points almost everywhere by Corollary \ref{cor:sizes}. Hence $N_I=4$ almost everywhere on $\A(V_I)$, $\Vol_4(\A(V_I))=\frac14\MVol_I(V)=\frac{\pi^4}4$, and the bounds follow from Theorem \ref{thm:proj}(b) with $|\Gamma|=2$. The characterisation of equality is Theorem \ref{thm:proj}(c), and the last assertion follows from the identity $\Hm^4(\A(V))=\int_VJ(\Log|_V)\,N(\Log z)^{-1}\,d\mathrm{vol}(z)$.
\end{proof}

Of the genericity of $V$ the proof only uses the injectivity of $\pi_I|_V$ and the general position of the five lines of each $V_I$. Both bounds are attained: the upper one for almost every real plane (Corollary \ref{cor:kplanes2}) and for every real plane satisfying the condition of Theorem \ref{thm:conic}, the lower one for the squares of lines (Corollary \ref{cor:squareproj}).

\subsection{Two-to-one amoeba maps of real planes}\label{sec:twotoone}

We now show that the condition of Corollary \ref{cor:plane5} is satisfied: the amoeba map of a generic real plane of $(\C^*)^5$ is two-to-one, so that its projected parts are exactly $\pi^4/2$. The argument is the analogue of the injectivity argument for lines of Section \ref{sec:injlines}. We give it for linear spaces of any dimension, first for almost every plane and then, for planes in $(\C^*)^5$, under an explicit condition on the coefficients analogous to the condition $c\notin\R$ for lines.

We measure sets of affine $k$-planes of $\C^n$ through the Lebesgue measure class of the real affine Grassmannian, for real planes, that is complexifications of affine $k$-planes of $\R^n$, and of the complex affine Grassmannian for complex planes. For $A\subset\{1,\dots,n\}$ with $|A|=k$, the planes whose projection to $\C^A$ is bijective form an open dense subset $\Pi_A$ with null complement, and each of them is the graph $z_m=\omega_{m,0}+\sum_{i\in A}\omega_{m,i}z_i$, $m\notin A$, of a unique affine map with real, respectively complex, coefficients. The map $V\mapsto\omega=(\omega_{m,i})$ is a smooth chart onto an open dense subset of $\R^{(n-k)(k+1)}$, respectively $\C^{(n-k)(k+1)}$, the transition maps are rational diffeomorphisms on the overlaps, and null sets correspond to null sets, so ``almost every plane'' can be checked in any chart $\Pi_A$.

Fix $I$ with $|I|=2k$, an index $j\notin I$, put $I'=I\cup\{j\}$, and use a chart $\Pi_A$ with $A\subset I$. The coefficients split as $\omega=(b,a)$, where $a\in\R^{k+1}$, respectively $\C^{k+1}$, is the row of $z_j$, and $b$ collects all other rows. We parametrise $V$ by $t=z_A\in\C^k$ and write $T(t)=(1,t)\in\C^{k+1}$ and $h_a(t)=a\cdot T(t)$, so that $z_j=h_a(t)$ on $V$. The projection $V_I=\pi_I(V)$ is the graph of the rows $m\in I\setminus A$ and depends only on $b$; it is the image of the injective map $t\mapsto w_b(t)$. Outside a null set of $b$, $V_I$ is in general position with respect to the coordinate hyperplanes of $\C^{2k}$ in the sense of \cite[Proposition 8 and Corollary 3]{Mik17}, and we assume this from now on.

\begin{lemma}\label{lem:hermgram}
Let $T,S\in\C^{k+1}$ with $T_0=S_0=1$. If $|a\cdot T|=|a\cdot S|$ for all $a\in\C^{k+1}$, then $T=S$.
\end{lemma}

\begin{proof}
We have $|a\cdot T|^2=a^tH_T\bar a$ with $H_T=T\,\overline T^{\,t}$, so the hypothesis says that the sesquilinear forms defined by the rank one matrices $T\overline T^{\,t}$ and $S\overline S^{\,t}$ agree on the diagonal; by polarisation the matrices are equal, hence $S=\lambda T$ with $|\lambda|=1$, by comparing images and traces, and $T_0=S_0=1$ gives $\lambda=1$.
\end{proof}

\begin{lemma}\label{lem:realgram}
Let $T,S\in\C^{k+1}$ with $T_0=S_0=1$. If $|a\cdot T|=|a\cdot S|$ for all $a\in\R^{k+1}$, then $S\in\{T,\overline T\}$.
\end{lemma}

\begin{proof}
Write $T=X+iY$, $S=U+iW$ with real vectors, so that $X_0=U_0=1$ and $Y_0=W_0=0$. For real $a$, $|a\cdot T|^2=a^t(XX^t+YY^t)a$, so the hypothesis is $Q_TQ_T^t=Q_SQ_S^t$ with $Q_T=[X\,|\,Y]$ and $Q_S=[U\,|\,W]$. Then $|Q_T^tv|=|Q_S^tv|$ for every $v\in\R^{k+1}$, so $Q_T^tv\mapsto Q_S^tv$ is a well defined isometry from the image of $Q_T^t$ onto the image of $Q_S^t$ in $\R^2$; extending it by an isometry of the orthogonal complements gives an orthogonal $2\times2$ matrix $O$ with $Q_S^t=OQ_T^t$ (compare \cite[Theorem 7.3.11]{HJ13}). If $O$ is a rotation, $S=e^{i\varphi}T$; if it is a reflection, $S=e^{i\psi}\overline T$. The coordinate $0$ gives $e^{i\varphi}=1$, respectively $e^{i\psi}=1$.
\end{proof}

Let $Z_b\subset\R^{2k}$ be the set of critical values of the amoeba map of $V_I\cap(\C^*)^{2k}$, a null set by Sard's theorem. For $y\notin Z_b$ the fibre of $V_I$ over $y$ is $\{w_b(t):t\in\Phi_b(y)\}$ for a finite set $\Phi_b(y)\subset\C^k$ independent of $a$. For fixed $(b,a)$, the image $H_{b,a}$ under $\Log$ of the points of $V_I$ whose lift to $V$ leaves the torus lies in the image of finitely many complex hyperplane sections of $V_I$, of real dimension $2k-2$, and is null; for $y\notin Z_b\cup H_{b,a}$, the points of $V\cap(\C^*)^n$ over $y$ are the lifts $z(t)$, $t\in\Phi_b(y)$, so $N_I(y)=\#\Phi_b(y)$, and two of them have different images under $\Log$ as soon as $|h_a(t)|\ne|h_a(s)|$.

\begin{proposition}\label{prop:sepfub}
In the real case, for almost every $\omega=(b,a)$ and almost every $y\in\R^{2k}\setminus Z_b$, any $t,s\in\Phi_b(y)$ with $s\notin\{t,\bar t\}$ satisfy $|h_a(t)|\ne|h_a(s)|$. In the complex case, for almost every $\omega=(b,a)$ and almost every $y$, any two distinct $t,s\in\Phi_b(y)$ satisfy $|h_a(t)|\ne|h_a(s)|$.
\end{proposition}

\begin{proof}
Consider the real case. Let $\mathcal B$ be the set of $(b,a,y)$ with $b$ generic, $y\notin Z_b$, and $t,s\in\Phi_b(y)$, $s\notin\{t,\bar t\}$, such that $|h_a(t)|^2=|h_a(s)|^2$. After the substitution $r=e^y$ it is the projection of a semialgebraic set, hence semialgebraic by the Tarski--Seidenberg theorem \cite[Theorem 2.2.1]{BCR98}, and in particular Borel. For fixed $(b,y)$ the set $\Phi_b(y)$ is finite and independent of $a$, and for each pair $(t,s)$ the function $a\mapsto|a\cdot T(t)|^2-|a\cdot T(s)|^2$ is a quadratic form on $\R^{k+1}$, not identically zero by Lemma \ref{lem:realgram}, since $T(s)\notin\{T(t),\overline{T(t)}\}$; its zero set is null. Hence every slice $\mathcal B_{b,y}$ is null, $\mathcal B$ is null by Fubini's theorem, and, by Fubini's theorem again, the slice $\mathcal B_{b,a}$ is null for almost every $(b,a)$. The complex case is identical, with Lemma \ref{lem:hermgram} and $a\in\C^{k+1}$.
\end{proof}

\begin{theorem}\label{thm:twotoone}
Let $k\ge1$ and $n>2k$. For almost every real affine $k$-plane $V\subset\C^n$ one has $S_I=N_I/2$ almost everywhere for every $I$, and the amoeba map of $V\cap(\C^*)^n$ is two-to-one at $\Hm^{2k}$-almost every point of its amoeba, the two points of the fibre being conjugate. For almost every complex affine $k$-plane one has $S_I=N_I$ almost everywhere for every $I$, and the amoeba map is injective at $\Hm^{2k}$-almost every point of its amoeba. Moreover, the exceptional set of planes is contained in a proper real algebraic subset of the space of planes.
\end{theorem}

\begin{proof}
Let $V$ be real and in the full-measure set of Proposition \ref{prop:sepfub} for every $I$ and every chart $\Pi_A$, $A\subset I$. For almost every $y$ with $N_I(y)>0$ the points of $V$ over $y$ are the lifts of the points of $\Phi_b(y)$; this set is invariant under $t\mapsto\bar t$, since $V_I$ is real, and the conjugation acts on it freely for almost every $y$, since the real locus has image of dimension $k<2k$ (Lemmas \ref{lem:fixed} and \ref{lem:count}). Points in the same conjugate pair have the same image under $\Log$, points in different pairs different images by Proposition \ref{prop:sepfub}, so $S_I(y)=N_I(y)/2$. Let $E=\{x\in\A(V):N(x)\ne2\}$. Off the null set $\Log(V(\R))$, $N$ is even and at least $2$, so up to a null set $E=\{x:N(x)\ge4\}$, and by \eqref{eq:NS}, $N_I(y)\ge2S_I(y)+2\#(E\cap p_I^{-1}(y))$; hence $E\cap p_I^{-1}(y)=\emptyset$ for almost every $y$, and the area formula \cite[Theorem 3.2.22]{Fed69} gives $\int_EJ^{\A(V)}_I\,d\Hm^{2k}=0$ for every $I$. Since $\sum_IJ_I^{\A(V)}\ge1$ by \eqref{eq:cb}, $\Hm^{2k}(E)=0$; this is also the equality case of Theorem \ref{thm:proj}(c) with $\Gamma=\{\mathrm{id},\conj\}$. The complex case is identical with the trivial group. For the last assertion, in a chart $\Pi_A$ and after the substitution $r=e^y$, the set of parameters for which the slice of the semialgebraic set $\mathcal B$ has positive measure, that is dimension $2k$, is semialgebraic, since the dimension of the fibres of a semialgebraic family is a semialgebraic function of the parameter \cite[Theorem 9.3.2]{BCR98}; being null, it has dimension smaller than that of the parameter space, and its Zariski closure is a proper real algebraic subset \cite[Proposition 2.8.2]{BCR98}.
\end{proof}

\begin{corollary}\label{cor:kplanes2}
For almost every real affine $k$-plane $V\subset\C^{2k+1}$, each of the $2k+1$ projected parts of the volume equals $\pi^{2k}/2$, $\Hm^{2k}(\A(V))=\frac12\MVol(V)$, and $\frac{\sqrt{2k+1}}2\pi^{2k}\le\Hm^{2k}(\A(V))\le\frac{2k+1}2\pi^{2k}$. For almost every complex affine $k$-plane $V\subset\C^{2k+1}$, each projected part equals $\pi^{2k}$, $\Hm^{2k}(\A(V))=\MVol(V)$, and $\sqrt{2k+1}\,\pi^{2k}\le\Hm^{2k}(\A(V))\le(2k+1)\pi^{2k}$. In particular, for almost every real plane $V\subset\C^5$ the five projected parts equal $\pi^4/2$, as asked in Corollary \ref{cor:plane5}.
\end{corollary}

\begin{proof}
For almost every $V$ and every $I$, $\pi_I|_V$ is injective and $\MVol_I(V)=\MVol(V_I)=\pi^{2k}$ \cite[Proposition 8]{Mik17}. Theorem \ref{thm:twotoone} gives $\int S_I\,dy=\frac12\int N_I\,dy=\frac{\pi^{2k}}2$ in the real case and $\pi^{2k}$ in the complex case, and $\Hm^{2k}(\A(V))=\int_VJ(\Log|_V)\,N(\Log z)^{-1}\,d\mathrm{vol}$ equals $\frac12\MVol(V)$, respectively $\MVol(V)$. The bounds follow, as in the proof of Theorem \ref{thm:main}, from $\|(\int_SJ^S_I)_I\|\le\Hm^{2k}(S)\le\sum_I\int_SJ^S_I$, where $S=\A(V)$.
\end{proof}

For $k=1$ this recovers Theorem \ref{thm:line}, there proved for every nonplanar line.

\subsection{\texorpdfstring{An explicit condition for real planes in $(\C^*)^5$}{An explicit condition for real planes in (C*)\textasciicircum 5}}\label{sec:conic}

For lines in $(\C^*)^3$ the injectivity of Section \ref{sec:injlines} holds for every line with $c\notin\R$, and the two-to-one property for every real line. For real planes in $(\C^*)^5$ we now give an explicit condition, in terms of the coefficients, under which the amoeba map is two-to-one everywhere outside the image of the real points, and we identify the planes which violate it.

Let $V\subset\C^5$ be a real affine plane, parametrised by $t\in\C^2$ through real affine functions $z_m=z_m(t)=\gamma_m+\delta_m\cdot t$, $m=1,\dots,5$, with $\gamma_m\in\R$ and $\delta_m\in\R^2$, and put $B_m=(\gamma_m,\delta_m)\in\R^3$ and $e_0=(1,0,0)$. In the projective closure $\overline V\simeq\mathbb{CP}^2$ of $V$, the points $[B_m]$ of the dual projective plane are the lines $\overline V\cap\{z_m=0\}$, and $[e_0]$ is the line at infinity $\overline V\cap\{z_0=0\}$.

\begin{theorem}\label{thm:conic}
If the six points $[e_0],[B_1],\dots,[B_5]$ of the dual projective plane lie on no common conic, that is, if the six lines in which the coordinate hyperplanes of $\mathbb{CP}^5$ cut $\overline V$ are not tangent to a common conic, then for every $z,z'\in V\cap(\C^*)^5$ with $\Log z=\Log z'$ one has $z'\in\{z,\bar z\}$. Consequently the amoeba map is two-to-one over $\A(V)\setminus\Log(V(\R))$, the projected parts equal $\frac12\MVol_I(V)$, and $\Hm^4(\A(V))=\frac12\MVol(V)$.
\end{theorem}

\begin{proof}
Let $z=z(t)$ and $z'=z(s)$, and put $u=t-s=u_1+iu_2$ and $v=t+s=v_1+iv_2$ with $u_i,v_i\in\R^2$. Since $z_m(t)-z_m(s)=\delta_m\cdot u$ and $z_m(t)+z_m(s)=2\gamma_m+\delta_m\cdot v$, and $|A|^2-|C|^2=\operatorname{Re}\big((A-C)\overline{(A+C)}\big)$, we get, using that $\gamma_m$ and $\delta_m$ are real,

\begin{equation}\label{eq:xi}
|z_m(t)|^2-|z_m(s)|^2=2\gamma_m\,\delta_m\cdot u_1+(\delta_m\cdot u_1)(\delta_m\cdot v_1)+(\delta_m\cdot u_2)(\delta_m\cdot v_2)=B_m^t\Xi B_m,
\end{equation}

where $\Xi$ is the symmetric $3\times3$ matrix with $\Xi_{00}=0$, first row and column $(0,u_1)$, and lower right block the symmetric part of $u_1v_1^t+u_2v_2^t$. The hypothesis $\Log z=\Log z'$ says that the quadratic form $\Xi$ vanishes at $B_1,\dots,B_5$, and it vanishes at $e_0$ because $\Xi_{00}=0$. Since no conic passes through the six points, $\Xi=0$. Thus $u_1=0$ and the symmetric part of $u_2v_2^t$ vanishes; for vectors $p,q\in\R^2$, $pq^t+qp^t=0$ forces $p=0$ or $q=0$, as one sees by applying it to $p$. If $u_2=0$, then $u=0$ and $s=t$; if $v_2=0$, then $t-s$ is purely imaginary and $t+s$ is real, that is $s=\bar t$. Hence $z'\in\{z,\bar z\}$. The remaining assertions follow from Theorem \ref{thm:proj}(c) and from the identity $\Hm^4(\A(V))=\int_VJ(\Log|_V)N(\Log z)^{-1}\,d\mathrm{vol}$.
\end{proof}

The condition of Theorem \ref{thm:conic} is the nonvanishing of the $6\times6$ determinant of the Veronese vectors $(x_0^2,x_0x_1,x_0x_2,x_1^2,x_1x_2,x_2^2)$ of the six points. If $V$ is written as a graph over $(z_1,z_2)$, so that $B_1=(0,1,0)$, $B_2=(0,0,1)$ and $z_m=\gamma_m+\delta_{m1}z_1+\delta_{m2}z_2$ for $m=3,4,5$, expanding along the rows of $e_0$, $B_1$ and $B_2$ reduces this determinant, up to sign, to

$$\Delta(V)=\det\begin{pmatrix}\gamma_3\delta_{31}&\gamma_3\delta_{32}&\delta_{31}\delta_{32}\\ \gamma_4\delta_{41}&\gamma_4\delta_{42}&\delta_{41}\delta_{42}\\ \gamma_5\delta_{51}&\gamma_5\delta_{52}&\delta_{51}\delta_{52}\end{pmatrix},$$

and when all the coefficients are nonzero, factoring $\gamma_m\delta_{m1}\delta_{m2}$ out of the row $m$ and exchanging the first and the last column gives $\Delta(V)=-\prod_{m=3}^5\gamma_m\delta_{m1}\delta_{m2}\cdot D(V)$, where $D(V)$ is the determinant of the $3\times3$ matrix with rows $(1/\gamma_m,\,1/\delta_{m1},\,1/\delta_{m2})$, $m=3,4,5$. Thus the amoeba map of a real plane in graph form is two-to-one outside the image of its real points as soon as the three vectors $(1/\gamma_m,1/\delta_{m1},1/\delta_{m2})$, $m=3,4,5$, are linearly independent. This is the analogue for planes of the condition on the cross-ratio of a line. The same argument shows that the amoeba map of a real $k$-plane in $(\C^*)^n$ is two-to-one outside the image of its real points whenever the $n+1$ points $[e_0],[B_1],\dots,[B_n]$ of $\mathbb{RP}^k$ lie on no common quadric, which requires $n+1\ge\binom{k+2}2$; for $k=1$ and $n=3$ it recovers the fact that every real nonplanar line of $(\C^*)^3$ has a two-to-one amoeba map (Section \ref{sec:injlines}).

The planes violating the condition have a transparent structure. Suppose that the six lines are tangent to a smooth conic $Q\subset\overline V$. Identify $\overline V$ with the symmetric square $\mathrm{Sym}^2\mathbb{CP}^1$, a point being an unordered pair $\{q_1,q_2\}$ of points of the conic $Q\simeq\mathbb{CP}^1$, namely the two points of tangency of the tangent lines through it; under this identification the tangent line at $p\in Q$ is the set of pairs containing $p$, and a linear form vanishing on it is, up to a factor, $\{q_1,q_2\}\mapsto(q_1-p)(q_2-p)$ in an affine coordinate on $Q$. If the tangency point of the line at infinity is sent to $\infty$, the coordinates of $V$ become $z_m=\nu_m\,(q_1-\varpi_m)(q_2-\varpi_m),\quad m=1,\dots,5,$ with $\nu_m\in\R^*$ and $\varpi_m\in\R$ when the configuration is real, the tangency points being real points of the real conic. Thus $V$ is, up to the multiplicative translation by $\nu=(\nu_m)$, the coordinatewise product $\ell\odot\ell=\{\ell(q_1)\ell(q_2)\}$ of the real line $\ell=\{(q-\varpi_m)_m:q\in\C\}\subset\C^5$ with itself. Its amoeba is $\Log|\nu|+\A(\ell)+\A(\ell)$, the Minkowski sum of two copies of the amoeba of the line, and, for $q_1,q_2\notin\R$ and $q_2\notin\{q_1,\bar q_1\}$, the four distinct points $\{a,b\}$ with $a\in\{q_1,\bar q_1\}$ and $b\in\{q_2,\bar q_2\}$ have the same image under $\Log$. By Theorem \ref{thm:square} below these four points form the whole fibre, so that the amoeba map of such a plane is exactly four-to-one, with $\Hm^4(\A(V))=\frac14\MVol(V)$ and projected parts $\frac14\MVol_I(V)=\frac{\pi^4}4$ (Corollary \ref{cor:squareproj}); the four points are the orbit of the partial conjugations of the two factors, the mechanism of Section \ref{sec:why}. When the conic is singular, the six lines pass through at most two points and $V$ has a product or cylinder structure.

\begin{corollary}\label{cor:hadamard}
Let $V\subset\C^5$ be a real affine plane whose six lines are in general position, no three concurrent. Then exactly one of the following holds. Either the six lines are not tangent to a common conic, and the amoeba map is two-to-one outside $\Log(V(\R))$, with $\Hm^4(\A(V))=\frac12\MVol(V)$ and projected parts $\pi^4/2$; or they are tangent to a smooth conic, $V$ is a multiplicative translate of $\ell\odot\ell$ for a real line $\ell$, its amoeba is a translate of $\A(\ell)+\A(\ell)$, the amoeba map has exactly $4$ points in almost every fibre, $\Hm^4(\A(V))=\frac14\MVol(V)$, and the projected parts equal $\pi^4/4$.
\end{corollary}

\begin{proof}
The first alternative is Theorem \ref{thm:conic}. If the six lines are tangent to a common conic, this conic is smooth, since a degenerate dual conic consists of the lines through one or two points, and six lines through at most two points include three concurrent ones. The construction above writes $V$ as a translate of $\ell\odot\ell$, and the remaining assertions are Theorem \ref{thm:square} and Corollary \ref{cor:squareproj} below.
\end{proof}

\section{Varieties of maximal volume in the half-dimensional case}\label{sec:maximal}

In this section $n=2k$ and $\T^{2k}=(S^1)^{2k}$.

\begin{lemma}\label{lem:alpha}
For almost every $\varepsilon\in\T^{2k}$ the set $V\cap\sigma_\varepsilon(V)$ consists of exactly $\alpha(V)$ points, all transversal intersections. For almost every $x\in\R^{2k}$ the number $N(x)$ is finite and even, and $N\ge2$ almost everywhere on $\A(V)$.
\end{lemma}

\begin{proof}
Since $\sigma_\varepsilon(V)=\varepsilon\conj(V)$, the first set is $V\cap\varepsilon\conj(V)$. The set of $\varepsilon\in(\C^*)^{2k}$ for which $V\cap\varepsilon\conj(V)$ is finite, transversal and of cardinality $\alpha(V)$ contains a nonempty Zariski open subset. Indeed, the points of $V\cap\varepsilon\conj(V)$ correspond bijectively to the points $(z,w)\in V\times\conj(V)$ with $z/w=\varepsilon$, so they form the fibre over $\varepsilon$ of the morphism $q:V\times\conj(V)\to(\C^*)^{2k}$, $q(z,w)=z/w$, between smooth varieties of the same dimension $2k$; here $\conj(V)$ is the variety defined by the polynomials $\bar f_j$ with conjugate coefficients. If $q$ is not dominant, its image lies in a proper algebraic subset and the intersection is empty for generic $\varepsilon$. Otherwise, by generic finiteness and generic smoothness in characteristic zero \cite[Exercise II.3.7 and Corollary III.10.7]{Har77}, applied to the components of $V\times\conj(V)$, there is a nonempty Zariski open set over which $q$ is finite \'etale, so the fibres there are finite, consist of transversal intersection points, and have constant cardinality, which is by definition the toric intersection number $\alpha(V)$ \cite[Definition 1.3]{Mik17}. The compact torus $\T^{2k}$ is Zariski dense in $(\C^*)^{2k}$, a polynomial vanishing on it vanishing identically, so its intersection with a proper algebraic subset is a proper real analytic subset of $\T^{2k}$, which is null. The second assertion is Lemma \ref{lem:parity} together with Sard's theorem, the amoeba having finite volume \cite{MN13}.
\end{proof}

\begin{proposition}\label{prop:A}
If $V$ is nondegenerate in the sense of Khovanskii, that is, for every $w\in\R^{2k}$ the initial system $\{\mathrm{in}_wf_j=0\}$ defines a smooth complete intersection of codimension $k$ in $(\C^*)^{2k}$ or the empty set, then $\alpha(V)=M$.
\end{proposition}

\begin{proof}
The polynomials $\bar f_j(z/\varepsilon)$ have the Newton polytopes $\Delta_j$, so by Bernstein's theorem \cite{Ber75} the number of solutions of the square system $f_j(z)=0$, $\bar f_j(z/\varepsilon)=0$, counted with multiplicity, equals $M$ as soon as, for every $w\ne0$, the initial system $\mathrm{in}_wf_j(z)=0$, $\mathrm{in}_w\bar f_j(z/\varepsilon)=0$ has no zero in the torus. The first half defines $W_w=\mathrm{in}_wV$, the second $\varepsilon\conj(W_w)$. For integral $w$ both are invariant under the one-parameter subgroup $\lambda_w(t)=t^w$, $t\in\C^*$; for an arbitrary $w$ the initial forms are those of an integral vector in the same cone of the normal fan of $\Delta_1+\dots+\Delta_k$, so it suffices to treat integral $w$. If $W_w\ne\emptyset$, it has dimension $k$, and the set $D_w$ of $\varepsilon\in(\C^*)^{2k}$ with $W_w\cap\varepsilon\conj(W_w)\ne\emptyset$ is the image of $W_w\times\conj(W_w)$ under $(z,u)\mapsto z/u$, which is constant along the orbits of the diagonal action of $\lambda_w$. Hence $D_w$ is contained in a constructible set of dimension at most $2k-1$, and its closure is a proper algebraic subset of $(\C^*)^{2k}$. There are finitely many cones, so for $\varepsilon$ outside a proper algebraic subset no initial system has a zero, and $\alpha(V)=M$.
\end{proof}

Proposition \ref{prop:A} complements \cite[Proposition 2]{Mik17}, which treats toric complete intersections.

\begin{proof}[Proof of Theorem \ref{thm:char}]
By Lemma \ref{lem:alpha}, for almost every $\varphi\in\T_\pi^{2k}$ the fibre $V\cap\Arg_\pi^{-1}(\varphi)=V\cap\Fix(\sigma_\varepsilon)$, $\varepsilon=e^{2i\varphi}$, is contained in the set $V\cap\sigma_\varepsilon(V)$ of $\alpha(V)$ points, because a fixed point $z=\sigma_\varepsilon(z)$ of $V$ lies in $\sigma_\varepsilon(V)$. Hence, by \eqref{eq:rolled}, $\MVol(V)\le\pi^{2k}\alpha(V)$, with equality if and only if for almost every $\varepsilon$ all the $\alpha(V)$ points are fixed, which is condition $(\mathrm B)$; the map $\varphi\mapsto e^{2i\varphi}$ is a diffeomorphism $\T_\pi^{2k}\to\T^{2k}$, so null sets correspond. By Lemma \ref{lem:alpha} again, $1\le N/2$ almost everywhere on $\A(V)$, and $N$ vanishes off $\A(V)$, so $\Vol(\A(V))\le\frac12\int_{\A(V)}N\,dx=\frac12\MVol(V)$, with equality if and only if $N=2$ almost everywhere on $\A(V)$. Since the fibre of $\Log|_V$ over $x$ is $V\cap\Fix(\sigma'_x)\subset V\cap\sigma'_x(V)$, this is condition $(\mathrm C)$. Finally $\alpha(V)\le M$ by Bernstein's theorem. Altogether $\Vol(\A(V))\le\frac12\MVol(V)\le\frac{\pi^{2k}}2\alpha(V)\le\frac{\pi^{2k}}2M$, and the outer terms are equal if and only if the three inequalities are equalities, that is, if and only if $(\mathrm C)$, $(\mathrm B)$ and $(\mathrm A)$ hold. The last assertion of $(\mathrm A)$ is Proposition \ref{prop:A}.

For the reformulation of $(\mathrm C)$ as a covering, let $Z\subset\R^{2k}$ be the set of critical values of $\Log|_V$, closed because $\Log|_V$ is proper and its critical set is closed, and null by Sard's theorem. Over $\R^{2k}\setminus Z$ the proper local diffeomorphism $\Log|_V$ is a finite covering, so $N$ is constant on each connected component of $\R^{2k}\setminus Z$, and these components are open. Hence $N=2$ almost everywhere on $\A(V)$ if and only if $N=2$ on $\A(V)\setminus Z$, that is, if and only if $\Log$ is a two-sheeted covering of $\A(V)\setminus Z$ by its preimage.
\end{proof}

\subsection{Curves}\label{sec:harnack}
For $k=1$ the polytope $\Delta$ is a lattice polygon, $M=2\Vol(\Delta)$, and the bound reads $\Vol(\A(V))\le\pi^2\Vol(\Delta)$, the inequality of Passare and Rullg{\aa}rd \cite{PR04}. Mikhalkin and Rullg{\aa}rd proved \cite{MR01} that equality holds if and only if, after a multiplicative translation, $V$ is a simple Harnack curve \cite{Mik00}. By Theorem \ref{thm:char}, the simple Harnack curves are therefore exactly the curves satisfying $(\mathrm A)$, $(\mathrm B)$ and $(\mathrm C)$, the deck involution of the covering in $(\mathrm C)$ being complex conjugation. In particular the hypothesis $k\ge2$ in Corollary \ref{cor:real4} cannot be removed (compare Remark \ref{rem:k1}).

\subsection{Symmetries and maximal volume}\label{sec:maxsym}

The divisibility results of Section \ref{sec:half} give a strong restriction on the varieties of maximal volume when $k\ge2$.

\begin{corollary}\label{cor:nomaxreal}
Let $n=2k$ with $k\ge2$, let $V\in\Vc(\Delta)$ with $M>0$, and assume that $V$ has maximal volume. Then the stabiliser $\G_V$ is trivial. In particular $V$ is not real up to a multiplicative translation.
\end{corollary}

\begin{proof}
By Theorem \ref{thm:char}, condition $(\mathrm C)$ holds, that is $N=2$ almost everywhere on $\A(V)$, and $\Vol(\A(V))=\frac{\pi^{2k}}2M>0$. If $\G_V$ were infinite, the volume would vanish by Lemma \ref{lem:infinite}. If $\G_V$ were finite and nontrivial, Theorem \ref{thm:half}(b), which applies because $k\ge2$, would give $N\equiv0\pmod{2|\G_V|}$, hence $N\ge4$, almost everywhere on $\A(V)$, a set of positive measure, contradicting $N=2$. The last assertion follows from Lemma \ref{lem:real}, a variety real up to translation having an antiholomorphic element in its stabiliser.
\end{proof}

For $k=1$ the situation is opposite: the varieties of maximal volume are the simple Harnack curves, which are real (Section \ref{sec:harnack}). Corollary \ref{cor:nomaxreal} also shows that conditions $(\mathrm A)$ and $(\mathrm B)$ do not imply $(\mathrm C)$, as the following example makes explicit.

\begin{example}\label{ex:products}
Let $C_1,C_2\subset(\C^*)^2$ be simple Harnack curves with Newton polygons $P_1,P_2$ and $V=C_1\times C_2\subset(\C^*)^4$, with $\Delta_1=P_1\times\{0\}$, $\Delta_2=\{0\}\times P_2$. Then $M=4\Vol(P_1)\Vol(P_2)$, as the coefficient of $\lambda^2\mu^2$ in $\Vol(\lambda\Delta_1+\mu\Delta_2)$ shows, $\MVol(V)=\MVol(C_1)\MVol(C_2)=\pi^4M$ and $N=N_1N_2=4$ almost everywhere on $\A(V)$ \cite{MR01,Mik00}. So $V$ satisfies $(\mathrm A)$ and $(\mathrm B)$, that is it is multiHarnack, in every degree \cite[Problem 3.3]{Mik17}, but $\Vol(\A(V))=\frac{\pi^4}4M$ is half of the bound, as Corollary \ref{cor:nomaxreal} requires.
\end{example}

\subsection{\texorpdfstring{Planes in $(\C^*)^4$ never have maximal volume}{Planes in (C*)\textasciicircum 4 never have maximal volume}}\label{sec:planesmax}

Let $V=\{z_3=a+bz_1+cz_2,\ z_4=d+ez_1+fz_2\}\subset(\C^*)^4$ be a plane in general position: $a,\dots,f$, $bf-ce$, $af-cd$, $ae-bd$ are nonzero, that is, the five boundary lines $H_1=\{z_1=0\}$, $H_2=\{z_2=0\}$, $H_3=\{z_3=0\}$, $H_4=\{z_4=0\}$ and the line at infinity $H_5$ of $V\subset\mathbb P^2$ are distinct and no three are concurrent. Then $M=\alpha(V)=1$ and $\MVol(V)=\pi^4$ \cite[Proposition 8]{Mik17}, so $V$ has maximal volume if and only if $N=2$ almost everywhere on $\A(V)$. For $\rho>0$ put
$$L_1(w)=\big(\re\tfrac{bw}a,\re\tfrac{ew}d\big),\quad L_2(w)=\big(\re\tfrac{cw}a,\re\tfrac{fw}d\big),$$
and $\Phi_\rho(\theta)=L_1(e^{i\theta_1})+\rho L_2(e^{i\theta_2})$,
the first-order model of the normal sections of $\A(V)$ near the cone of $\Trop(V)$ where $z_1,z_2\to0$.

\begin{lemma}\label{lem:firstorder}
If some $y_0$ is a regular value of some $\Phi_{\rho_0}$ with at least four preimages, then $N\ge4$ on a nonempty open set.
\end{lemma}

\begin{proof}
Put $x(s,\rho,y)=(\log s,\log s\rho,\log|a|+sy_1,\log|d|+sy_2)$. The points of $V$ over $x(s,\rho,y)$ are the $(se^{i\theta_1},s\rho e^{i\theta_2})$ with $\Phi_{s,\rho}(\theta)=y$, where $\Phi_{s,\rho}(\theta)=s^{-1}\big(\log|1+(bse^{i\theta_1}+cs\rho e^{i\theta_2})/a|,\log|1+(ese^{i\theta_1}+fs\rho e^{i\theta_2})/d|\big)$ extends analytically to $s=0$ with $\Phi_{0,\rho}=\Phi_\rho$, because $\log|1+sw|=s\re w+O(s^2)$. The implicit function theorem at the four regular preimages of $y_0$ gives four distinct preimages of $y$ under $\Phi_{s,\rho}$ for $|s|,|\rho-\rho_0|,|y-y_0|<\delta$, hence $N\ge4$ on the image of $(0,\delta)\times(\rho_0-\delta,\rho_0+\delta)\times B(y_0,\delta)$ under the injective map $x$, whose Jacobian determinant is $s/\rho\ne0$.
\end{proof}

\begin{lemma}\label{lem:ellipses}
Let $L_1,L_2:\C\to\R^2$ be nonzero real linear maps with distinct images if both have rank one. Then some $\Phi_\rho$ has a regular value with at least four preimages, unless $L_1,L_2$ are invertible and $L_2=\lambda L_1\circ\varsigma$ with $\lambda>0$ and $\varsigma$ orthogonal.
\end{lemma}

\begin{proof}
The preimages of $0$ correspond to the intersection points of the curves $E_1=L_1(S^1)$ and $\rho E_2$, counted with their numbers of preimages in $S^1$, and $0$ is regular if the curves cross transversally with nonzero velocities. An invertible $L_j$ gives an ellipse centred at $0$, each point having one preimage; a map of rank one gives a segment centred at $0$, each interior point having two preimages and nonzero velocity. If both are ellipses and are not homothetic, write $E_j=\{q_j=1\}$ with non-proportional positive definite forms, let $\mu_1<\mu_2$ be the eigenvalues of $q_2$ relative to $q_1$, and take $\mu_1<\rho^2<\mu_2$: then $q_2-\rho^2q_1$ is indefinite and nondegenerate, and $E_1\cap\rho E_2=\{q_1=1,q_2=\rho^2\}$ consists of the four points where its two zero lines meet $E_1$; by B\'ezout's theorem two conics meeting in four points meet transversally. If $E_1$ is a segment and $E_2$ an ellipse, for $\rho$ small the segment crosses $\rho E_2$ transversally at two interior points, each with two preimages; and two segments with distinct directions cross at $0$, which has $2\times2$ preimages with independent velocities. Finally $E_1,E_2$ are homothetic ellipses exactly when $\lambda^{-1}L_1^{-1}L_2$ preserves the unit circle.
\end{proof}

With $\beta_1=b/a$, $\beta'_1=e/d$, $\beta_2=c/a$, $\beta'_2=f/d$, the map $L_j$ is invertible if and only if $\beta_j/\beta'_j\notin\R$, two maps of rank one have the same image if and only if $\beta_1/\beta'_1=\beta_2/\beta'_2$, that is $bf=ce$, which is excluded, and $L_2=\lambda L_1\circ\varsigma$ means $\beta_2/\beta'_2=\beta_1/\beta'_1$ for a rotation (excluded likewise) and $\beta_2/\beta'_2=\overline{\beta_1/\beta'_1}$ for a reflection. So if $N=2$ almost everywhere, Lemmas \ref{lem:firstorder} and \ref{lem:ellipses} give
\begin{equation}\label{eq:mirror}
\kappa:=\frac{bd}{ae}\notin\R,\qquad\frac{cd}{af}=\bar\kappa .
\end{equation}
With the cross-ratio $[P,Q;R,S]=\frac{(R-P)(S-Q)}{(S-P)(R-Q)}$, so that $[0,\infty;R,S]=R/S$, we have $\kappa=[H_2\cap H_1,H_2\cap H_5;H_2\cap H_4,H_2\cap H_3]$ on $H_2$ (coordinate $z_1$, points $0,\infty,-d/e,-a/b$) and $\frac{cd}{af}=[H_1\cap H_2,H_1\cap H_5;H_1\cap H_4,H_1\cap H_3]$ on $H_1$. The other nine cones correspond to the other points $H_i\cap H_j$; putting a third line $H_k$ at infinity and using the coordinates $\ell_i/\ell_k$, $\ell_j/\ell_k$, where $\ell_r$ defines $H_r$, changes the amoeba by an invertible integral linear map, which preserves the counting function and the null sets. Hence, if $N=2$ almost everywhere, then for all distinct $i,j,k,m,m'$
\begin{equation}\label{eq:mirrorij}
\begin{gathered}[H_i\cap H_j,H_i\cap H_k;H_i\cap H_{m'},H_i\cap H_m]\\=\overline{[H_j\cap H_i,H_j\cap H_k;H_j\cap H_{m'},H_j\cap H_m]}\notin\R .\end{gathered}
\end{equation}

\begin{lemma}\label{lem:conic}
If no three of the lines $H_1,\dots,H_5$ are concurrent, there are distinct $t_1,\dots,t_5\in\mathbb P^1$ with $[H_j\cap H_i,H_j\cap H_k;H_j\cap H_{m'},H_j\cap H_m]=[t_i,t_k;t_{m'},t_m]$ for all distinct $i,j,k,m,m'$.
\end{lemma}

\begin{proof}
The dual points $h_r$ of the lines $H_r$ lie on a unique smooth conic $Q\simeq\mathbb P^1$; let $t_r=h_r$. The points $H_j\cap H_r$ of $H_j$ correspond to the lines $h_jh_r$, with the same cross-ratios, and by Chasles' theorem the cross-ratio of the lines joining $h_j\in Q$ to four other points of $Q$ is the cross-ratio of these points on $Q$ \cite[Chapter 16]{Ber87}.
\end{proof}

By Lemma \ref{lem:conic}, \eqref{eq:mirrorij} reads $[t_j,t_k;t_{m'},t_m]=\overline{[t_i,t_k;t_{m'},t_m]}\notin\R$. If $g$ is the M\"obius map sending $t_k,t_{m'},t_m$ to $\infty,1,0$, then $[x,\infty;1,0]=(x-1)/x$ has real coefficients, so the condition says $g(t_j)=\overline{g(t_i)}\notin\R$: $t_i$ and $t_j$ are symmetric with respect to the circle through the three other points, and do not lie on it.

\begin{lemma}\label{lem:fivepoints}
No five distinct points of $\mathbb P^1$ are such that any two are symmetric with respect to the circle through the three others.
\end{lemma}

\begin{proof}
Normalise $t_3=0$, $t_4=1$, $t_5=\infty$. The pair $\{1,2\}$ gives $t_2=\bar t_1\notin\R$; the pair $\{3,4\}$, symmetric in the line $\re z=\re t_1$, gives $\re t_1=\frac12$; the pair $\{3,5\}$ forces the circle through $t_1,\bar t_1,1$ to be centred at $0$, so $|t_1|=1$, and after conjugation, which preserves circles and symmetry, $t_1=e^{i\pi/3}$. The pair $\{1,3\}$ then requires $t_1$ to be the reflection $1+u^2\,\overline{(0-1)}=1-u^2=\frac32-\frac{\sqrt3}2i$ of $0$ in the line through $1$ with direction $u=e^{-2i\pi/3}$, which contains $\bar t_1$; a contradiction.
\end{proof}

\begin{theorem}\label{thm:planesgp}
For every plane $V$ in general position in $(\C^*)^4$, real or not, $N\ge4$ on a set of positive measure and $\Vol(\A(V))<\frac12\MVol(V)=\frac{\pi^4}2$; no such plane has maximal volume.
\end{theorem}

\begin{proof}
If $N=2$ almost everywhere, \eqref{eq:mirrorij} holds for the ten pairs $\{i,j\}$, contradicting Lemmas \ref{lem:conic} and \ref{lem:fivepoints}. So $N\ge4$, $N$ being even and positive (Lemma \ref{lem:alpha}), on a set of positive measure, and $\Vol(\A(V))\le\int_{\A(V)}\frac N2-\Vol(\{N\ge4\})<\frac12\MVol(V)$.
\end{proof}

This answers in part \cite[Problem 3.1]{Mik17}. We expect the same argument, with Minkowski sums of $k$ ellipses in $\R^k$, to exclude $k$-planes for every $k\ge2$.

\section{Symmetric powers of lines}\label{sec:squares}

The planes which violate the condition of Theorem \ref{thm:conic} are coordinatewise squares of lines. In this section we study, more generally, coordinatewise products of $k$ points of one line, for any $k\ge1$ and any number $n\ge2k$ of variables. Let $\varpi_1,\dots,\varpi_n$ be distinct complex numbers and let $\ell=\{\ell(q):q\in\C\}$, $\ell(q)=(q-\varpi_1,\dots,q-\varpi_n)$, be the corresponding line of $\C^n$; we call $\ell$ real if all $\varpi_m$ are real. Every line of $\C^n$ which meets the coordinate hyperplanes in $n$ distinct finite points has this form after an affine change of the parameter and a multiplicative translation, since each coordinate is an affine function of the parameter which can be rescaled to be monic; for a real line both can be chosen real. A complex affine change $q\mapsto\alpha q+\beta$ replaces $\varpi$ by $\alpha\varpi+\beta$ and $\ell$ by its multiplicative translate by $\alpha^{-1}(1,\dots,1)$, so everything below depends on $\varpi$ only up to complex affine maps, and up to real affine maps for real lines.

For $q=(q_1,\dots,q_k)\in\C^k$ put $\Phi(q)=\ell(q_1)\odot\dots\odot\ell(q_k)$, the coordinatewise product, so that

$$\Phi(q)_m=\prod_{i=1}^k(q_i-\varpi_m)=(-1)^kf_q(\varpi_m),\qquad f_q(x)=\prod_{i=1}^k(x-q_i)=x^k+c_1x^{k-1}+\dots+c_k .$$

The coordinates of $\Phi(q)$ are affine functions of the coefficients $c=(c_1,\dots,c_k)$ of $f_q$, and the map $c\mapsto\big((-1)^kf(\varpi_m)\big)_m$ is injective as soon as $n\ge k$, by the Vandermonde determinant. Since $\{q_1,\dots,q_k\}\mapsto c$ identifies the $k$-th symmetric power of $\C$ with $\C^k$, the image $\ell^{\odot k}=\{\Phi(q):q\in\C^k\}\subset\C^n$ is an affine $k$-plane, real when $\ell$ is real, which we call the $k$-th symmetric power of $\ell$; the map $\Phi$ is $k!$-to-one off the big diagonal. The points of $\ell^{\odot k}$ in the torus are the $\Phi(q)$ with all $q_i\notin\{\varpi_1,\dots,\varpi_n\}$, and, writing $\lambda(q)=\Log\ell(q)$, we have $\Log\Phi(q)=\lambda(q_1)+\dots+\lambda(q_k)$. Hence the amoeba of $\ell^{\odot k}$ is the Minkowski sum $\A(\ell)^{[k]}=\A(\ell)+\dots+\A(\ell)\quad(k\ \text{summands}).$ A multiplicative translate $\nu\cdot\ell^{\odot k}$, $\nu\in(\C^*)^n$, has the amoeba $\Log|\nu|+\A(\ell)^{[k]}$, and its fibres are the images under multiplication by $\nu$ of those of $\ell^{\odot k}$, so all statements below hold for translates. For $k=1$ we recover the line itself, and for $k=2$ the squares $\ell\odot\ell$ of Section \ref{sec:conic}.

The symmetric powers are recognised through their hyperplane sections. The projective closure of $\ell^{\odot k}$ is the projective space of polynomials of degree at most $k$; the coordinate hyperplane $\{z_m=0\}$ cuts it in the hyperplane $\{f:f(\varpi_m)=0\}$, which corresponds in the dual projective space to the point $(1,\varpi_m,\dots,\varpi_m^k)$, and the hyperplane at infinity corresponds to $(0,\dots,0,1)$. These $n+1$ points lie on the rational normal curve $t\mapsto(1,t,\dots,t^k)$, the hyperplane at infinity corresponding to $t=\infty$. Conversely, if the $n+1$ points of the dual projective space defined by the coordinate hyperplanes and the hyperplane at infinity lie on a rational normal curve of degree $k$, a projective change of its parameter sends the point at infinity to $t=\infty$, and the $k$-plane is a multiplicative translate of a symmetric power; for a real $k$-plane and a real curve, the parameters of its real points are real. Through $k+3$ points in general linear position of $\mathbb{P}^k$ passes a unique rational normal curve of degree $k$ \cite[Theorem 1.18]{Har92}, which is real if the points are real, since its conjugate is another such curve. Hence every $k$-plane of $\C^{k+2}$ whose $k+3$ hyperplane sections are in general position is a multiplicative translate of a symmetric power of a line. For $k=2$ this is the case $n=4$ treated in Corollary \ref{cor:square4} below; for $k\ge3$ one has $2k>k+2$, and in the half-dimensional case $n=2k$ the symmetric powers form a proper family. For instance, real affine $3$-planes of $\C^6$ form a family of dimension $12$, while the symmetric cubes of real lines and their multiplicative translates by real vectors depend on $4+6=10$ parameters, namely $\varpi$ up to real affine maps and the translation, so they have codimension $2$.

\subsection{Real symmetric powers}

\begin{theorem}\label{thm:square}
Let $\ell$ be a real line, $k\ge1$, $n\ge2k$, and let $q=(q_1,\dots,q_k)$ and $r=(r_1,\dots,r_k)$ have entries in $\C\setminus\{\varpi_1,\dots,\varpi_n\}$. Then $\sum_i\lambda(q_i)=\sum_i\lambda(r_i)$ if and only if, after renumbering the $r_i$, one has $r_i\in\{q_i,\bar q_i\}$ for every $i$. Consequently the fibre of the amoeba map of $V=\ell^{\odot k}$ through $\Phi(q)$ is the set of points $\Phi(a)$ with $a_i\in\{q_i,\bar q_i\}$ for every $i$; it has exactly $2^k$ points when no $q_i$ is real and $q_j\notin\{q_i,\bar q_i\}$ for $i\ne j$, and in general as many points as there are distinct multisets $\{a_1,\dots,a_k\}$ of this form. The amoeba map of $V$ is exactly $2^k$-to-one at $\Hm^{2k}$-almost every point of $\A(V)$, and $\Hm^{2k}(\A(\ell)^{[k]})=2^{-k}\MVol(V)$.
\end{theorem}

\begin{proof}
Put $f=f_q$ and $g=f_r$, and let $\bar f,\bar g$ be the polynomials with conjugate coefficients, so that $f(x)\bar f(x)=|f(x)|^2$ and $g(x)\bar g(x)=|g(x)|^2$ for real $x$. The hypothesis says that $|f(\varpi_m)|=|g(\varpi_m)|$ for $m=1,\dots,n$. The polynomial $P=f\bar f-g\bar g$ has real coefficients, and since $f\bar f$ and $g\bar g$ are both monic of degree $2k$, it has degree at most $2k-1$. It vanishes at the $n\ge2k$ distinct real numbers $\varpi_m$, hence $P=0$, and by unique factorisation in $\C[x]$ the multisets of roots of $f\bar f$ and $g\bar g$ agree: $\biguplus_{i=1}^k\{q_i,\bar q_i\}=\biguplus_{i=1}^k\{r_i,\bar r_i\}.$ We show by induction on $k$ that this forces, after renumbering, $r_i\in\{q_i,\bar q_i\}$. For $k=0$ there is nothing to prove. The number $r_k$ belongs to the left-hand side, so after renumbering the $q_i$ we may assume $r_k\in\{q_k,\bar q_k\}$; then $\{r_k,\bar r_k\}=\{q_k,\bar q_k\}$ as multisets, and removing this multiset from both sides gives the same identity for $k-1$ pairs, to which the induction hypothesis applies. Conversely, $|a-\varpi_m|=|\bar a-\varpi_m|$ for every real $\varpi_m$, so $\lambda(\bar a)=\lambda(a)$, and all the points $\Phi(a)$ with $a_i\in\{q_i,\bar q_i\}$ have the same image under $\Log$.

Since $\Phi$ induces a bijection from multisets onto $\ell^{\odot k}$, the fibre through $\Phi(q)$ consists of the points $\Phi(a)$, counted as multisets. Under the stated genericity, two choices $a,a'$ giving the same multiset satisfy $a_i=a'_{\pi(i)}$ for a permutation $\pi$, and $a'_{\pi(i)}\in\{q_{\pi(i)},\bar q_{\pi(i)}\}$ meets $\{q_i,\bar q_i\}$ only if $\pi(i)=i$; hence $\pi$ is the identity and $a=a'$, and the $2^k$ choices give $2^k$ distinct points, the choices within each pair being distinct since $q_i\notin\R$. Let $E\subset(\C\setminus\{\varpi_m\})^k$ be the set of $q$ such that some $q_i$ is real or $q_j\in\{q_i,\bar q_i\}$ for some $i\ne j$; it is a finite union of real submanifolds of dimension at most $2k-1$, invariant under the conjugation of any single entry. If a point $x$ of $\A(V)$ has a preimage $\Phi(q)$ with $q\in E$, then all its preimages $\Phi(a)$ have $a\in E$, so the set of points of $\A(V)$ over which the fibre does not have $2^k$ points is $\Log(\Phi(E))$, which is $\sigma$-finite for $\Hm^{2k-1}$ and therefore $\Hm^{2k}$-null \cite[Section 2.10.2]{Fed69}. The identity $\Hm^{2k}(\A(V))=\int_VJ(\Log|_V)N(\Log z)^{-1}\,d\mathrm{vol}(z)$, a consequence of the area formula \cite[Theorem 3.2.3]{Fed69}, gives the last assertion.
\end{proof}

The statement holds at every point, not only almost everywhere, like the injectivity criterion $c\notin\R$ for lines of Section \ref{sec:injlines}. The $2^k$ points of a generic fibre form the orbit of $q$ under the group $(\Z/2)^k$ generated by the conjugations of the single entries $q_i$; only the conjugation of all entries, which is the conjugation $z\mapsto\bar z$ of $\C^n$, belongs to $\G$, the others being neither holomorphic nor antiholomorphic on $V$. These are the partial conjugations of the $k$ factors, the mechanism which produces the count $2^k$ for products of real lines in Section \ref{sec:why}, here acting on a symmetric product instead of a product; as there, their fixed point sets $\{q_i\in\R\}$ have real codimension one in $V$. In the half-dimensional case $n=2k$ the theorem shows that the amoeba map of a real symmetric power is exactly $2^k$-to-one, so that the value $\pi^{2k}/2^k$ stated in \cite[Theorem 1.2]{NP17} for real generic $k$-planes is correct for symmetric powers (Corollary \ref{cor:squareproj}), although the $2^k$-to-one property from which it was derived fails for all other real $k$-planes in general position when $k\ge3$ (Theorem \ref{thm:classif} and Corollary \ref{cor:fibre4}); for $k=3$ the fibres of $8$ points are among those allowed by Corollary \ref{cor:sizes}. The hypothesis $n\ge2k$ is used in the proof, since for $n=2k-1$ the polynomial $P$ may vanish at $n$ points without vanishing identically; for $k=2$ and $n=5$ the failure of the condition of Theorem \ref{thm:conic} appears in its proof as follows: for two points $z(t),z(s)$ of a fibre with $s\notin\{t,\bar t\}$, the quadratic form $\Xi$ of \eqref{eq:xi} is nonzero and vanishes at the six points $[e_0],[B_1],\dots,[B_5]$, that is, it is an equation of the conic.

To compute the volume, write $v_m(q)=\frac{d}{dq}\log(q-\varpi_m)=1/(q-\varpi_m)$ and $v(q)=a(q)+ib(q)$ with $a(q),b(q)\in\R^n$, and put $a_i=a(q_i)$, $b_i=b(q_i)$. As in Section \ref{sec:lines}, the differential of $\lambda$ at $q=u_1+iu_2$ sends $\partial_{u_1}$ to $a(q)$ and $\partial_{u_2}$ to $-b(q)$, so the Jacobian of $q\mapsto\sum_i\lambda(q_i)$ with respect to the Lebesgue measure $dA(q)=dA(q_1)\cdots dA(q_k)$ of $\C^k$ is the norm of the $2k$-vector $a_1\wedge b_1\wedge\dots\wedge a_k\wedge b_k\in\Lambda^{2k}\R^n$. For $|I|=2k$ we denote by $J_I(q)$ the absolute value of its Pl\"ucker coordinate of index $I$, that is of the minor of rows $I$ of the matrix $[a_1\,|\,b_1\,|\,\dots\,|\,a_k\,|\,b_k]$. For a regular point $x$ of $\A(\ell)$ we denote by $\xi(x)$ a unit $2$-vector spanning the tangent plane of $\A(\ell)$ at $x$, which is $\pm a(q)\wedge b(q)/|a(q)\wedge b(q)|$ for $x=\lambda(q)$.

\begin{proposition}\label{prop:squarevol}
Let $\ell$ be a real line, $k\ge1$, $n\ge2k$ and $V=\ell^{\odot k}$. Then

$$\Hm^{2k}\big(\A(\ell)^{[k]}\big)=\frac{\MVol(V)}{2^k}=\frac1{2^kk!}\int_{\C^k}|a_1\wedge b_1\wedge\dots\wedge a_k\wedge b_k|\,dA(q),$$

$$\Hm^{2k}\big(\A(\ell)^{[k]}\big)=\frac1{k!}\int_{\A(\ell)^k}|\xi(x_1)\wedge\dots\wedge\xi(x_k)|\,d\Hm^2(x_1)\cdots d\Hm^2(x_k),$$

where $|a_1\wedge\dots\wedge b_k|=(\sum_{|I|=2k}J_I^2)^{1/2}$ and, for $q_1,\dots,q_k\notin\R$,

$$J_I(q)=\frac{\prod_i|\operatorname{Im}q_i|\,\prod_{i<j}|q_i-q_j|^2|q_i-\bar q_j|^2\,\prod_{a<b\in I}|\varpi_a-\varpi_b|}{\prod_{m\in I}\prod_{i}|q_i-\varpi_m|^2}.$$
\end{proposition}

\begin{proof}
By the chain rule for Jacobians, $J(\Log|_V)(\Phi(q))\,J\Phi(q)=|a_1\wedge\dots\wedge b_k|$, where $J\Phi$ is the Jacobian of $\Phi$ for the Lebesgue measure of $\C^k$ and the induced metric of $V$. Since $\Phi$ is $k!$-to-one off the big diagonal, which is null, the area formula \cite[Theorem 3.2.3]{Fed69} gives $\int_{\C^k}|a_1\wedge\dots\wedge b_k|\,dA=k!\int_VJ(\Log|_V)\,d\mathrm{vol}=k!\MVol(V)$, and Theorem \ref{thm:square} gives the first two equalities. For the third, the amoeba map of the real line $\ell$ is two-to-one off the real axis: if $|q-\varpi_m|=|q'-\varpi_m|$ for all $m$, the real polynomial $(x-q)(x-\bar q)-(x-q')(x-\bar q')$, of degree at most $1$, vanishes at $n\ge2$ points, so $q'\in\{q,\bar q\}$. Hence the area formula gives $\int_{\A(\ell)}h\,d\Hm^2=\frac12\int_\C h(\lambda(q))\,|a(q)\wedge b(q)|\,dA(q)$ for every nonnegative Borel function $h$, and applying it in each factor gives $\int_{\A(\ell)^k}|\xi_1\wedge\dots\wedge\xi_k|=2^{-k}\int_{\C^k}|a_1\wedge\dots\wedge b_k|\,dA$. The expression of the norm through the minors is the Cauchy--Binet identity of Section \ref{sec:notation}.

For the minors, replace the columns $a(q_i),b(q_i)$ by $v(q_i)=a+ib$ and $\bar v(q_i)=a-ib$; the change of basis $[v\,|\,\bar v]=[a\,|\,b]\begin{pmatrix}1&1\\ i&-i\end{pmatrix}$ has determinant $-2i$, so the minor of rows $I$ of the matrix $[v(q_1)\,|\,\bar v(q_1)\,|\,\dots\,|\,v(q_k)\,|\,\bar v(q_k)]$ has absolute value $2^kJ_I$. Since the $\varpi_m$ are real, $\bar v_m(q)=1/(\bar q-\varpi_m)$, so this matrix is the Cauchy matrix $\big(1/(y_c-x_r)\big)$ with $x=(\varpi_m)_{m\in I}$ and $y=(q_1,\bar q_1,\dots,q_k,\bar q_k)$, whose determinant is $\pm\prod_{r<r'}(x_r-x_{r'})\prod_{c<c'}(y_c-y_{c'})/\prod_{r,c}(y_c-x_r)$ by Cauchy's classical evaluation (see \cite[Chapter 0]{HJ13}). Here $\prod_{c<c'}|y_c-y_{c'}|=2^k\prod_i|\operatorname{Im}q_i|\prod_{i<j}|q_i-q_j|^2|q_i-\bar q_j|^2$, because $|q_i-\bar q_i|=2|\operatorname{Im}q_i|$, $|\bar q_i-\bar q_j|=|q_i-q_j|$ and $|\bar q_i-q_j|=|q_i-\bar q_j|$, and $\prod_{r,c}|y_c-x_r|=\prod_{m\in I}\prod_i|q_i-\varpi_m|^2$, which gives the formula.
\end{proof}

The formula shows that $J_I>0$ whenever no $q_i$ is real and the $q_i$ are pairwise distinct and not conjugate, so the amoeba of every projection of a real symmetric power has full dimension $2k$. The projected multivolumes are computed directly from the rolled coamoeba, without appeal to the genericity in \cite[Proposition 8]{Mik17}.

\begin{proposition}\label{prop:coamsym}
Let $\ell$ be a real line, $k\ge1$, $n\ge2k$, $V=\ell^{\odot k}$ and $|I|=2k$. Then almost every fibre of $\Arg_{\pi,I}|_V$ has exactly one point, so that $\MVol_I(V)=\pi^{2k}$. Equivalently, for any $2k$ distinct real numbers $\varpi_m$, $m\in I$,

$$\int_{\C^k}\frac{\prod_i|\operatorname{Im}q_i|\,\prod_{i<j}|q_i-q_j|^2|q_i-\bar q_j|^2\,\prod_{a<b\in I}|\varpi_a-\varpi_b|}{\prod_{m\in I}\prod_{i}|q_i-\varpi_m|^2}\,dA(q_1)\cdots dA(q_k)=k!\,\pi^{2k}.$$
\end{proposition}

\begin{proof}
Write the points of $V$ as $z_m=(-1)^kf(\varpi_m)$ with $f(x)=x^k+c_1x^{k-1}+\dots+c_k$, $c\in\C^k$. A point lies in $\Arg_{\pi,I}^{-1}(\varphi)$, $\varphi\in\T^I_\pi$, if and only if $z_m\ne0$ and $\operatorname{Im}(e^{-i\varphi_m}z_m)=0$ for $m\in I$. These are $2k$ real affine equations in the $2k$ real unknowns $(\operatorname{Re}c,\operatorname{Im}c)$, and the determinant $D(\varphi)$ of their linear part is a trigonometric polynomial in $\varphi$. Split $I=A\sqcup B$ with $|A|=|B|=k$ and take $\varphi_m=0$ for $m\in A$ and $\varphi_m=\pi/2$ for $m\in B$. Since the $\varpi_m$ are real, the equations become $\sum_j\operatorname{Im}(c_j)\varpi_m^{k-j}=0$ for $m\in A$ and $\sum_j\operatorname{Re}(c_j)\varpi_m^{k-j}=-\varpi_m^k$ for $m\in B$, a block diagonal system whose determinant is, up to sign, the product of the two Vandermonde determinants of $(\varpi_m)_{m\in A}$ and $(\varpi_m)_{m\in B}$, which is nonzero. Hence $D$ is not identically zero, and, being real analytic, it vanishes only on a null set. For $\varphi$ with $D(\varphi)\ne0$ the system has exactly one solution $c(\varphi)$, and it gives a point of $V\cap(\C^*)^n$ unless $z_m(c(\varphi))=0$ for some $m$. For $m\notin I$, the set of $\varphi$ for which this happens is contained in the image under the rolled argument map $\Arg_{\pi,I}$ of the hyperplane section $\{z_m=0\}$ of $V$, of real dimension $2k-2$, hence it is null. For $m\in I$, if $z_m(c(\varphi))=0$ then $c(\varphi)$ lies in $\{z_m=0\}$ and $\varphi_{m'}$ is the rolled argument of $z_{m'}(c(\varphi))$ for $m'\in I\setminus\{m\}$, so $\varphi$ lies in the image of a manifold of dimension $(2k-2)+1$ under a smooth map, again a null set. Thus almost every fibre has exactly one point, and \eqref{eq:rolled} gives $\MVol_I(V)=\pi^{2k}$. Since $\pi_I|_V$ is injective, $\MVol_I(V)=\MVol(\pi_I(V))$, and the area formula applied to $\pi_I\circ\Phi$, which is $k!$-to-one, as in the proof of Proposition \ref{prop:squarevol}, gives the integral identity.
\end{proof}

For $k=1$ the identity is \eqref{eq:passare}, the identity used by Passare in his proof of $\zeta(2)=\pi^2/6$ \cite{Pas08}; Proposition \ref{prop:coamsym} is its extension to Cauchy determinants of any size.

\begin{corollary}\label{cor:squareproj}
Let $\ell$ be a real line, $k\ge1$, $n\ge2k$, $V=\ell^{\odot k}$ and $|I|=2k$. Then $\pi_I|_V$ is injective, $V_I=\pi_I(V)$ is the $k$-th symmetric power of the line $\pi_I(\ell)$ of $\C^{2k}$, almost every point of $\A(V_I)$ has $N_I=2^k$ preimages in $V$ and $S_I=1$ preimage in $\A(V)$, and

$$\Vol_{2k}(\A(V_I))=\int_{\A(V)}J^{\A(V)}_I\,d\Hm^{2k}=\frac{\pi^{2k}}{2^k}=\frac{\MVol_I(V)}{2^k}.$$

In particular, for every real line $\ell\subset\C^{2k}$ meeting the coordinate hyperplanes in $2k$ distinct finite points, the Minkowski sum $\A(\ell)^{[k]}$ of $k$ copies of its amoeba has volume $\pi^{2k}/2^k$. For $n\ge2k$ one has $2^{-k}\binom n{2k}^{1/2}\pi^{2k}\le\Hm^{2k}(\A(\ell)^{[k]})\le2^{-k}\binom n{2k}\pi^{2k}$; for $k=2$ and $n=5$ this reads $\frac{\sqrt5}4\pi^4\le\Hm^4(\A(\ell)+\A(\ell))\le\frac54\pi^4$.
\end{corollary}

\begin{proof}
The injectivity of $\pi_I|_V$ follows from the Vandermonde determinant, as observed above, and $\pi_I\circ\Phi$ is the map $\Phi$ of the line $\pi_I(\ell)$. By Proposition \ref{prop:coamsym}, $\MVol_I(V)=\pi^{2k}$. By Theorem \ref{thm:square} applied with $n=2k$ to $\pi_I(\ell)$, almost every fibre of the amoeba map of $V_I$ has $2^k$ points, so $N_I=2^k$ almost everywhere on $\A(V_I)$ and $\Vol_{2k}(\A(V_I))=2^{-k}\MVol_I(V)$. These points are the $\Phi(a)$ with $a_i\in\{q_i,\bar q_i\}$, which by Theorem \ref{thm:square} for $V$ have the same image under $\Log$, so $S_I=1$, and $\int_{\A(V)}J^{\A(V)}_I\,d\Hm^{2k}=\int S_I\,dy=\Vol_{2k}(\A(V_I))$. The bounds follow from Theorem \ref{thm:main} as in the proof of Corollary \ref{cor:kplanes2}.
\end{proof}

Thus for $k=2$ the projected parts of the volume of a square of a line attain the lower bound $\frac{\pi^4}4$ of Corollary \ref{cor:plane5}, whereas those of almost every real plane attain the upper bound $\frac{\pi^4}2$ (Corollary \ref{cor:kplanes2}). Theorem \ref{thm:proj} allows both, and the two examples show that the inequality $S_I\le N_I/2$ of Theorem \ref{thm:proj} and the trivial inequality $S_I\ge1$ are both sharp.

\begin{corollary}\label{cor:square4}
Let $V\subset\C^4$ be a real affine plane such that no three of the five lines cut on its projective closure by the coordinate hyperplanes and the hyperplane at infinity are concurrent. Then $V$ is the multiplicative translate of $\ell\odot\ell$ by a real vector, for a real line $\ell\subset\C^4$, and the fibres of its amoeba map are described by Theorem \ref{thm:square}; in particular they have exactly $4$ points outside the image of a real $3$-dimensional subset of $V$, and $\Vol(\A(V))=\frac{\pi^4}4$. Equivalently, for every real line $\ell\subset(\C^*)^4$ meeting the coordinate hyperplanes in four distinct finite points, the Minkowski sum $\A(\ell)+\A(\ell)$ has volume $\frac{\pi^4}4$.
\end{corollary}

\begin{proof}
In the notation of Section \ref{sec:conic}, the five points $[e_0],[B_1],\dots,[B_4]$ of the dual real projective plane are in general position, so they lie on a unique conic, a rational normal curve of degree $2$; it is real, its equation being the solution of a real linear system with one-dimensional solution space, and smooth, since a pair of lines through five points would contain three of them. By the discussion at the beginning of this section, $V$ is the multiplicative translate of $\ell\odot\ell$ by a real vector $\nu$, with distinct real $\varpi_m$, distinct points of the conic having distinct parameters, and Corollary \ref{cor:squareproj} with $k=2$ and $n=4$ gives the rest.
\end{proof}

Corollary \ref{cor:square4} gives an elementary proof of the case $k=2$ of Corollary \ref{cor:sizes}, and hence of the value $\pi^4/4$ of Nisse and Passare \cite[Theorem 1.2]{NP17}, which uses neither the topology of the amoeba map nor Mikhalkin's upper bound for the fibres, and which identifies every fibre explicitly as an orbit of partial conjugations. It also explains Example \ref{ex:planes}: the generic real plane of $(\C^*)^4$ is a symmetric square in disguise. In higher dimension Corollary \ref{cor:squareproj} gives, for instance, real $3$-planes of $(\C^*)^6$ whose amoebas have volume exactly $\pi^6/8$.

\subsection{\texorpdfstring{The volume of $\A(\ell)+\A(\ell)$ in $(\C^*)^5$}{The volume of A(l)+A(l) in (C*)\textasciicircum 5}}

For $k=2$ and $n=5$ the volume $\Hm^4(\A(\ell)+\A(\ell))$ is a function of the configuration $\varpi$ of five distinct real numbers up to real affine maps and permutations, that is of three real moduli; normalising $\varpi=(0,a,b,c,1)$ with $0<a<b<c<1$, the moduli space is an open simplex, on which the reflection $q\mapsto1-q$ acts by $(a,b,c)\mapsto(1-c,1-b,1-a)$. It is the analogue for squares of the function $\Mv$ of Section \ref{sec:cross}, which governs lines in $(\C^*)^3$, and its boundary behaviour is described through $\Mv$. We first compute its limits when two points collide.

\begin{theorem}\label{thm:collide}
Let $\varpi_1,\varpi_3,\varpi_4,\varpi_5$ be distinct real numbers, let $\ell_\varepsilon$ be the real line with points $\varpi^\varepsilon=(\varpi_1,\varpi_1+\varepsilon,\varpi_3,\varpi_4,\varpi_5)$, and put $c=(\varpi_4-\varpi_3)/(\varpi_5-\varpi_3)$. Then

$$\lim_{\varepsilon\to0}\Hm^4\big(\A(\ell_\varepsilon)+\A(\ell_\varepsilon)\big)=\frac{\sqrt2}4\,\pi^4+\frac{\pi^2}4\,\Mv(c),$$

where $\Mv(c)=\MVol(L_c)$ is the multivolume of the line $L_c$ of Section \ref{sec:lines}.
\end{theorem}

\begin{proof}
By Proposition \ref{prop:squarevol}, $\Hm^4(\A(\ell_\varepsilon)+\A(\ell_\varepsilon))=\frac18\int_{\C^2}F_\varepsilon\,dA$, with $F_\varepsilon=(\sum_IJ_I^\varepsilon{}^2)^{1/2}$ and $J_I^\varepsilon$ given by the Cauchy formula for $\varpi^\varepsilon$; we also put $G_\varepsilon=\sum_IJ^\varepsilon_I$, so that $F_\varepsilon\le G_\varepsilon$ and, by Proposition \ref{prop:coamsym}, $\int_{\C^2}J_I^\varepsilon\,dA=2\pi^4$ for every $I$ and every $\varepsilon\ne0$, hence $\int G_\varepsilon=10\pi^4$. Fix $\delta>0$, let $0<\varepsilon<\delta/2$, and split $\C^2$ into $D_0=\{|q_1-\varpi_1|>\delta,\ |q_2-\varpi_1|>\delta\}$, $D_1=\{|q_1-\varpi_1|\le\delta<|q_2-\varpi_1|\}$, its mirror image $D_1'$ under $(q_1,q_2)\mapsto(q_2,q_1)$, on which the integrands take the same values, and $D_{11}=\{|q_1-\varpi_1|\le\delta,\ |q_2-\varpi_1|\le\delta\}$.

On $D_0$ the functions $J_I^\varepsilon$ converge pointwise, as $\varepsilon\to0$, to the functions $J^0_I$ given by the Cauchy formula for the collided configuration $(\varpi_1,\varpi_1,\varpi_3,\varpi_4,\varpi_5)$. For the three sets $I\supset\{1,2\}$ the factor $|\varpi_1-\varpi_2|$ gives $J^0_I=0$; for $I=\{1,3,4,5\}$ and $I=\{2,3,4,5\}$ both limits equal the Jacobian $J'$ of the square of the line $\ell'$ of $\C^4$ with points $(\varpi_1,\varpi_3,\varpi_4,\varpi_5)$. Hence $F_\varepsilon\to F_0=\sqrt2\,J'$ and $G_\varepsilon\to G_0=2J'$ pointwise on $D_0$, and $\int_{\C^2}J'\,dA=2\pi^4$ by Proposition \ref{prop:coamsym}.

On $D_1$ we substitute $q_1=\varpi_1+\varepsilon u$, so that $dA(q_1)=\varepsilon^2dA(u)$ and $D_1$ becomes $\{|u|\le\delta/\varepsilon\}\times\{|q_2-\varpi_1|>\delta\}$, which increases to $\C\times\{|q_2-\varpi_1|>\delta\}$. For $I=\{1,2,m,m'\}$ with $\{m,m'\}\subset\{3,4,5\}$ the Cauchy formula gives, as $\varepsilon\to0$ with $u$ and $q_2$ fixed,

$$\varepsilon^2J_I^\varepsilon(\varpi_1+\varepsilon u,q_2)\longrightarrow\frac{|\operatorname{Im}u|}{|u|^2|u-1|^2}\cdot\frac{|\operatorname{Im}q_2|\,|\varpi_m-\varpi_{m'}|}{|q_2-\varpi_m|^2|q_2-\varpi_{m'}|^2},$$

because the factors $|\operatorname{Im}q_1|=\varepsilon|\operatorname{Im}u|$ and $|\varpi_1-\varpi_2|=\varepsilon$ of the numerator and $|q_1-\varpi_1|^2|q_1-\varpi_2|^2=\varepsilon^4|u|^2|u-1|^2$ of the denominator carry all the powers of $\varepsilon$, while $|q_1-q_2|^2|q_1-\bar q_2|^2\to|q_2-\varpi_1|^4$ cancels against $|q_2-\varpi_1|^2|q_2-\varpi_2|^2$, and $|\varpi_1-\varpi_m||\varpi_2-\varpi_m|\to|\varpi_1-\varpi_m|^2$ cancels against $|q_1-\varpi_m|^2$, and similarly for $m'$. The first factor is the function $|w_{12}|$ of Lemma \ref{lem:formula} for the line $\{(u,u-1)\}$, whose integral is $\pi^2$ by \eqref{eq:passare}, and the second is $|w_{mm'}(q_2)|$ for the line $\ell''=\{(q-\varpi_3,q-\varpi_4,q-\varpi_5)\}$, as one sees from $\frac12|\operatorname{Im}q||\varpi_m-\varpi_{m'}|=\frac12\rho_m\rho_{m'}|\sin\theta_{mm'}|$, twice the area of a triangle computed in two ways. For $I=\{1,3,4,5\}$ or $\{2,3,4,5\}$, the factor $|\operatorname{Im}q_1|=\varepsilon|\operatorname{Im}u|$ and the single factor $\varepsilon^2|u|^2$, respectively $\varepsilon^2|u-1|^2$, of the denominator give $\varepsilon^2J_I^\varepsilon=O(\varepsilon)\to0$. Hence on $D_1$, in the coordinates $(u,q_2)$, the functions $\varepsilon^2F_\varepsilon$ converge pointwise to $A(u)\,W(q_2)$ and $\varepsilon^2G_\varepsilon$ to $A(u)\sum_{m<m'}|w_{mm'}(q_2)|$, where $A(u)=|\operatorname{Im}u|/(|u|^2|u-1|^2)$ and $W=(\sum_{m<m'}w_{mm'}^2)^{1/2}$ is the integrand of $\Mv$ for $\ell''$, which is affinely equivalent to $L_c$, so that $\int_\C W\,dA=\Mv(c)$.

Let $\Lambda(\delta)$ be the sum of the integrals of the pointwise limits of $G_\varepsilon$ over $D_0$, $D_1$ and $D_1'$, computed in the coordinates above. By the previous paragraphs and monotone convergence, $\Lambda(\delta)\to2\cdot2\pi^4+2\cdot\pi^2\cdot3\pi^2=10\pi^4$ as $\delta\to0$, using $\int|w_{mm'}|\,dA=\pi^2$ for each of the three pairs, again by \eqref{eq:passare}. Put $\eta(\delta)=10\pi^4-\Lambda(\delta)\ge0$. Fatou's lemma on $D_0$, $D_1$ and $D_1'$, together with $\int G_\varepsilon=10\pi^4$, gives $\limsup_{\varepsilon\to0}\int_RG_\varepsilon\le\int_R\lim G_\varepsilon+\eta(\delta)$ for $R=D_0,D_1,D_1'$, and $\limsup_{\varepsilon\to0}\int_{D_{11}}G_\varepsilon\le\eta(\delta)$. We now use the following form of Pratt's lemma \cite{Pra60}: if $0\le f_\varepsilon\le g_\varepsilon$, $f_\varepsilon\to f$ and $g_\varepsilon\to g$ pointwise and $\limsup\int g_\varepsilon\le\int g+\eta<\infty$, then $\int f\le\liminf\int f_\varepsilon\le\limsup\int f_\varepsilon\le\int f+\eta$; the left inequality is Fatou's lemma, and the right one follows from Fatou's lemma applied to $g_\varepsilon-f_\varepsilon\ge0$. Applied to $F_\varepsilon\le G_\varepsilon$ on $D_0$, $D_1$ and $D_1'$, and combined with $\int_{D_{11}}F_\varepsilon\le\int_{D_{11}}G_\varepsilon$, it shows that every limit point of $\int_{\C^2}F_\varepsilon\,dA$ lies within $4\eta(\delta)$ of

$$\int_{D_0}\sqrt2\,J'\,dA+2\pi^2\int_{|q_2-\varpi_1|>\delta}W(q_2)\,dA(q_2).$$

As $\delta\to0$ this tends to $2\sqrt2\,\pi^4+2\pi^2\Mv(c)$ and $\eta(\delta)\to0$, so $\int_{\C^2}F_\varepsilon\,dA\to2\sqrt2\,\pi^4+2\pi^2\Mv(c)$, and dividing by $8$ gives the claim.
\end{proof}

The two terms have a geometric meaning. The first is the volume $\sqrt2\cdot\frac{\pi^4}4$ of the image of the amoeba of $\ell'\odot\ell'$, of volume $\frac{\pi^4}4$ by Corollary \ref{cor:square4}, under the linear embedding $(x_1,x_3,x_4,x_5)\mapsto(x_1,x_1,x_3,x_4,x_5)$, whose Jacobian is $\sqrt2$. The second is the volume of the piece of the amoeba coming from pairs $\{q_1,q_2\}$ with $q_1$ at distance of order $\varepsilon$ from the colliding points: there the first two coordinates describe a translate of the amoeba of a real line of $(\C^*)^2$, of area $\frac{\pi^2}2$, and the last three the amoeba of the real line $\ell''$, of area $\frac12\Mv(c)$, and the piece is a sheared product of the two.

\begin{corollary}\label{cor:supsq}
For real lines $\ell\subset(\C^*)^5$, the supremum of $\Hm^4(\A(\ell)+\A(\ell))$ is at least $\frac{1+2\sqrt2}4\pi^4$, and the limits along single collisions of Theorem \ref{thm:collide} lie in the interval $\big(\frac{\sqrt2+\sqrt3}4\pi^4,\frac{\sqrt2+3}4\pi^4\big)$.
\end{corollary}

\begin{proof}
By Theorem \ref{thm:line} and Theorem \ref{thm:linesummary}, $\sqrt3\,\pi^2<\Mv(c)<3\pi^2$, which gives the interval. By Theorem \ref{thm:degen}, $\Mv(c)\to(1+\sqrt2)\pi^2$ as $c\to0$, so the limits of Theorem \ref{thm:collide} come arbitrarily close to $\frac{\sqrt2}4\pi^4+\frac{1+\sqrt2}4\pi^4=\frac{1+2\sqrt2}4\pi^4$, and so do the volumes themselves.
\end{proof}

The value $\frac{1+2\sqrt2}4\pi^4$ is also the limit expected along double collisions, when $\varpi_2-\varpi_1$ and $\varpi_4-\varpi_3$ both tend to $0$: the same decomposition gives no contribution from $D_0$, since the limiting line has only three distinct points and the corresponding square has a $3$-dimensional amoeba, the contribution $2\cdot\pi^2\cdot\sqrt2\,\pi^2$ for each of the two clusters, the line $\ell''$ now having two equal coordinates, and the contribution $2\pi^4$ from pairs with one point near each cluster, for a total of $\frac18(4\sqrt2+2)\pi^4=\frac{1+2\sqrt2}4\pi^4$; we have not written the details of this second computation.

We do not know the extremal values of $\Hm^4(\A(\ell)+\A(\ell))$ on the moduli space of configurations. What we can prove is the lower bound $\frac{\sqrt5}4\pi^4$ of Corollary \ref{cor:squareproj}, which is not sharp, and the lower bound $\frac{1+2\sqrt2}4\pi^4$ for the supremum of Corollary \ref{cor:supsq}. The comparison with lines in $(\C^*)^3$ is instructive: there the multivolume tends to $(1+\sqrt2)\pi^2$ at collisions and stays below this value near them (Theorem \ref{thm:degen} and Proposition \ref{prop:rate}), and Theorem \ref{thm:collide} shows that the collisions of two points of $\ell$ play the same role for $\A(\ell)+\A(\ell)$.

\subsection{Symmetric powers of nonreal lines}

When the points $\varpi_m$ are not all real, the conjugation of a single entry $q_i$ no longer preserves $|q_i-\varpi_m|$, and the polynomial argument of Theorem \ref{thm:square} has to be replaced by its analogue in two real variables. For $k\ge1$ let $\mathscr P_k$ be the real vector space of functions on $\C$ spanned by $\operatorname{Re}(w^\alpha\bar w^\beta)$ and $\operatorname{Im}(w^\alpha\bar w^\beta)$ for $0\le\beta\le\alpha\le k$ and $(\alpha,\beta)\ne(k,k)$; the functions $\operatorname{Im}(w^\alpha\bar w^\alpha)$ vanish, and since the monomials $w^\alpha\bar w^\beta$ are linearly independent, $\dim\mathscr P_k=k+2\binom{k+1}2=k^2+2k$. For monic $f,g$ of degree $k$ with coefficients $c_\alpha,d_\alpha$ ($c_k=d_k=1$), the expansion $|f(w)|^2=\sum_{\alpha,\beta}c_\alpha\bar c_\beta w^\alpha\bar w^\beta=\sum_\alpha|c_\alpha|^2|w|^{2\alpha}+2\sum_{\alpha>\beta}\operatorname{Re}(c_\alpha\bar c_\beta w^\alpha\bar w^\beta)$ shows that $P_{f,g}(w)=|f(w)|^2-|g(w)|^2$ belongs to $\mathscr P_k$, the terms $|w|^{2k}$ cancelling. Moreover $P_{f,g}=0$ only if $f=g$: if $|f|=|g|$ on $\C$, then $f/g$ is a meromorphic function of modulus one, hence a constant of modulus one, which equals $1$ since $f$ and $g$ are monic. For $k=1$, $\mathscr P_1$ is the space of real affine functions of $w$, and $P_{f,g}=0$ is the perpendicular bisector of the roots of $f$ and $g$; for $k=2$, $\mathscr P_2$ has dimension $8$ and the curves $P=0$ are the real cubics of the form $|w|^2\operatorname{Re}(\beta w)+Q(w)=0$ with $Q$ of degree at most $2$, that is the circular cubics, which pass through the two circular points at infinity.

\begin{theorem}\label{thm:nonrealsq}
Let $k\ge1$, let $\varpi\in\C^n$ have distinct entries and let $V=\ell^{\odot k}$.

(a) If $n\ge k^2+2k$ and no nonzero element of $\mathscr P_k$ vanishes at all the points $\varpi_1,\dots,\varpi_n$, then $\Log|_{V\cap(\C^*)^n}$ is injective. For $k=1$ this condition says that the $\varpi_m$ are not collinear; for $n=3$ this is the condition $c\notin\R$ of Section \ref{sec:injlines}.

(b) If $n\ge2k+1$, then for almost every $\varpi\in\C^n$, the exceptional set being contained in a proper real algebraic subset of $\C^n$, the amoeba map of $V$ is injective at $\Hm^{2k}$-almost every point of $\A(V)$, and $\Hm^{2k}(\A(\ell)^{[k]})=\MVol(V)=\frac1{k!}\int_{\C^k}|a_1\wedge b_1\wedge\dots\wedge a_k\wedge b_k|\,dA(q)$.
\end{theorem}

\begin{proof}
Two points $\Phi(q)$ and $\Phi(r)$ of $V$ in the torus have the same image under $\Log$ if and only if $P_{f_q,f_r}(\varpi_m)=0$ for every $m$. Under the hypothesis of (a) this forces $P_{f_q,f_r}=0$, hence $f_q=f_r$ and $\Phi(q)=\Phi(r)$. For $k=1$, $\mathscr P_1$ is spanned by $1,\operatorname{Re}w,\operatorname{Im}w$, and a nonzero element vanishes on a line.

For (b), identify the coefficient vectors $c$ of monic polynomials of degree $k$ with $\C^k$ and consider the set

$$Z=\{(\varpi,c,c')\in\C^n\times\C^k\times\C^k:\ c\ne c',\ P_{f_c,f_{c'}}(\varpi_m)=0\ \text{for}\ m=1,\dots,n\},$$

which is semialgebraic. For fixed $(c,c')$ with $c\ne c'$, the function $P_{f_c,f_{c'}}$ is a nonzero real polynomial in $(\operatorname{Re}w,\operatorname{Im}w)$, whose zero set in $\C$ has dimension at most $1$; hence the fibre of $Z$ over $(c,c')$, a product of $n$ such zero sets, has dimension at most $n$, and $\dim Z\le n+4k$ \cite[Theorem 9.3.2]{BCR98}. The projection $Z'$ of $Z$ to $\C^n\times\C^k$ is semialgebraic of dimension at most $n+4k$ \cite[Chapter 2]{BCR98}, which is smaller than $2n+2k$ because $n>2k$; hence $Z'$ is null, and by Fubini's theorem, for almost every $\varpi$ the slice $Z'_\varpi=\{c:(\varpi,c)\in Z'\}$ is a null subset of $\C^k$. For such $\varpi$, the set $B\subset V$ of points whose fibre under $\Log$ contains another point of $V$ corresponds to a subset of $Z'_\varpi$, so it is null in $V$, and the set of points of $\A(V)$ with more than one preimage is $\Log(B)$, because every preimage of such a point lies in $B$; since $\Log$ is locally Lipschitz, $\Log(B)$ is $\Hm^{2k}$-null. The volume formula follows from the area formula as in the proof of Proposition \ref{prop:squarevol}, now with $N=1$. Finally, the set of $\varpi$ for which $Z'_\varpi$ has dimension $2k$ is semialgebraic, since the dimension of the fibres of a semialgebraic family is a semialgebraic function of the parameter \cite[Theorem 9.3.2]{BCR98}; being null, it has dimension smaller than $2n$, and its Zariski closure is a proper real algebraic subset \cite[Proposition 2.8.2]{BCR98}.
\end{proof}

Thus the expectation that squares of nonreal lines should have two-to-one amoeba maps is not correct: the partial conjugations which produce the fibres of $2^k$ points of real symmetric powers disappear together with the total conjugation, and for generic nonreal lines the amoeba map of $\ell^{\odot k}$ is injective almost everywhere when $n\ge2k+1$, exactly as for generic complex $k$-planes in Theorem \ref{thm:twotoone}, which does not apply directly since symmetric powers form a null family of $k$-planes. For $k=2$ part (a) requires eight points not on a common circular cubic, while part (b) covers $n=5,6,7$ only for almost every $\varpi$. When $n=2k$ the amoeba map of a symmetric power of a nonreal line is proper, of degree zero, so it is never injective.

\begin{corollary}\label{cor:nonrealproj}
Let $k\ge1$, $n\ge2k+1$, and let $\varpi\in\C^n$ be outside the exceptional set of Theorem \ref{thm:nonrealsq}(b) and outside the set $Z_G$ of those $\varpi$ for which the amoeba of some projection $\pi_I(V)$, $|I|=2k$, has volume zero; the set $Z_G$ is contained in a proper real algebraic subset of $\C^n$. Then for every $I$ with $|I|=2k$ one has $\int_{\A(V)}J^{\A(V)}_I\,d\Hm^{2k}=\MVol_I(V)=\pi^{2k}$, and $\binom n{2k}^{1/2}\pi^{2k}\le\Hm^{2k}(\A(\ell)^{[k]})\le\binom n{2k}\pi^{2k}$.
\end{corollary}

\begin{proof}
The projection $\pi_I(V)$ is the $k$-th symmetric power of $\pi_I(\ell)$, a $k$-plane of $(\C^*)^{2k}$, which is the image of the affine map $c\mapsto(f_c(\varpi_m))_{m\in I}$ from the space $\C^k$ of coefficients of monic polynomials of degree $k$, with coefficients polynomial in $\varpi$. By Lemma \ref{lem:G}, $\MVol_I(V)\in\{0,\pi^{2k}\}$, with $\MVol_I(V)=\pi^{2k}$ exactly when $\pi_I(V)$ satisfies (G), that is when $\varpi\notin Z_G$. Since the linear part of this affine map has rank $k$ by the Vandermonde determinant, $\pi_I(V)$ satisfies (G) if and only if the square real linear system $(f_c(\varpi_m))_{m\in I}=D_\varphi r$ of size $4k$ in the unknowns $(r,\operatorname{Re}c,\operatorname{Im}c)$ has a determinant which is not identically zero in $\varphi$ and a solution none of whose coordinates $r_j$ vanishes identically. As in the proof of Lemma \ref{lem:G}, the determinant and the Cramer numerators of the $r_j$ are trigonometric polynomials in $\varphi$ whose coefficients are real polynomials in $(\operatorname{Re}\varpi,\operatorname{Im}\varpi)$, so $Z_G$ is contained in a real algebraic subset of $\C^n$; this subset is proper because it contains no real $\varpi$ with distinct entries, for which $\MVol_I(V)=\pi^{2k}$ by Proposition \ref{prop:coamsym}. Hence $\MVol_I(V)=\pi^{2k}$. Since the amoeba map of $V$ is injective almost everywhere and $\Log_I(B)$ is null for the set $B$ of the proof of Theorem \ref{thm:nonrealsq}, one has $S_I=N_I$ almost everywhere, and $\int S_I\,dy=\int N_I\,dy=\MVol_I(V)$. The bounds follow from Theorem \ref{thm:main}.
\end{proof}

\section{\texorpdfstring{Products of symmetric powers and $2^k$-to-one amoeba maps}{Products of symmetric powers and 2\textasciicircum k-to-one amoeba maps}}\label{sec:products}

The real $k$-planes of $(\C^*)^{2k}$ known so far to have $2^k$-to-one amoeba maps are the products of real lines in disjoint pairs of coordinates (Section \ref{sec:why}) and the symmetric powers of real lines (Theorem \ref{thm:square}). Both are special cases of products of symmetric powers in disjoint sets of coordinates, and it is natural to ask whether these are the only real $k$-planes with $2^k$-to-one amoeba maps. In this section we show that the answer depends on the position of the coordinate hyperplanes. Without any condition the answer is negative: products of symmetric powers admit fibred extensions, which are not products and still have $2^k$-to-one amoeba maps (Theorem \ref{thm:fibred}). Under the natural condition of general position, which excludes all products with more than one factor, the answer is positive in the strongest form: a real $k$-plane in general position in $\C^{2k}$ whose amoeba map has at least $2^k$ points in almost every fibre is a symmetric power of a real line (Theorem \ref{thm:classif}). As a consequence, every real $3$-plane in general position in $(\C^*)^6$ which is not a symmetric cube has fibres of exactly $4$ points on a set of positive measure.

Throughout, a real affine $k$-plane $V\subset\C^n$, $n\ge k$, is parametrised as $z_m=\gamma_m+\delta_m\cdot t$, $t\in\C^k$, with $\gamma_m\in\R$ and $\delta_m\in\R^k$, and as in Section \ref{sec:conic} its dual points are the points $p_0=[e_0]$ and $p_m=[B_m]$, $B_m=(\gamma_m,\delta_m)$, of the dual projective space $(\mathbb P^k)^*$, corresponding to the hyperplane at infinity and to the hyperplanes $\{z_m=0\}$ of the projective closure $\overline V\simeq\mathbb P^k$. They determine $V$ up to a multiplicative translation by a real vector and an affine change of the parameter. We say that $V$ is in general position if any $k+1$ of the $n+1$ points $p_0,\dots,p_n$ are linearly independent, that is, if any $k+1$ of the $n+1$ hyperplanes of $\overline V$ have empty intersection. By Section \ref{sec:squares}, $V$ is a multiplicative translate of a symmetric power of a real line if and only if its dual points lie on a real rational normal curve of degree $k$, the parameters of the points being distinct and real; since any $k+1$ distinct points of a rational normal curve are linearly independent, by the Vandermonde determinant, symmetric powers are in general position. In $\C^{2k}$ the real $k$-planes form, up to multiplicative translations by real vectors, a family of dimension $(2k+1)k-((k+1)^2-1)=k^2-k$, namely configurations of $2k+1$ points of $(\mathbb P^k)^*$ up to projective transformations, while the symmetric powers form a family of dimension $2k+1-3=2k-2$, namely $2k+1$ points of the projective line up to projective transformations; the codimension of the symmetric powers is therefore $(k-1)(k-2)$, which vanishes for $k=2$ and equals $2$ for $k=3$.

\subsection{Products of symmetric powers}

Let $\ell_1,\dots,\ell_r$ be real lines, $\ell_a\subset\C^{n_a}$ with distinct real points $\varpi^{(a)}\in\R^{n_a}$, let $k_1,\dots,k_r\ge1$ with $n_a\ge2k_a$, and put $k=k_1+\dots+k_r$ and $n=n_1+\dots+n_r$. The product $V=\ell_1^{\odot k_1}\times\dots\times\ell_r^{\odot k_r}\subset\C^{n_1}\times\dots\times\C^{n_r}=\C^n$ is a real affine $k$-plane, and we call it, as well as its images under permutations of the coordinates and multiplicative translations, a product of symmetric powers. For $k_a=1$ and $n_a=2$ the factors are real lines of $(\C^*)^2$, and for $r=1$ we recover the symmetric powers.

\begin{proposition}\label{prop:prodfib}
Let $V$ be a product of symmetric powers as above. Then the fibre of the amoeba map of $V$ through a point $(\Phi_1(q^{(1)}),\dots,\Phi_r(q^{(r)}))$ is the product of the fibres of the amoeba maps of the factors, described by Theorem \ref{thm:square}; it is an orbit of the group $(\Z/2)^k$ of the conjugations of the single entries $q^{(a)}_i$, and it has exactly $2^k$ points at $\Hm^{2k}$-almost every point of $\A(V)$. Moreover $\A(V)=\prod_a\A(\ell_a)^{[k_a]}$, $\MVol(V)=\prod_a\MVol(\ell_a^{\odot k_a})$ and

$$\Hm^{2k}(\A(V))=2^{-k}\MVol(V)=\prod_{a=1}^r\Hm^{2k_a}\big(\A(\ell_a)^{[k_a]}\big);$$

in particular $\Hm^{2k}(\A(V))=\pi^{2k}/2^k$ when $n_a=2k_a$ for every $a$. If $r\ge2$, then $V$ is not in general position.
\end{proposition}

\begin{proof}
Since $\Log$ acts coordinatewise, the amoeba map of $V$ is the product of the amoeba maps of the factors, so its fibres are the products of their fibres, and $\A(V)$ is the product of the amoebas $\A(\ell_a^{\odot k_a})=\A(\ell_a)^{[k_a]}$. By Theorem \ref{thm:square} each factor has fibres which are orbits of the conjugations of the single entries, with exactly $2^{k_a}$ points outside the image of a set of dimension $2k_a-1$; the product of these exceptional sets with the other factors has dimension at most $2k-1$, so the fibres of $V$ have exactly $2^k$ points outside an $\Hm^{2k}$-null set. The Jacobian of a product of maps is the product of the Jacobians, and the volume element of a product of Riemannian manifolds is the product of the volume elements, so $\MVol(V)=\prod_a\MVol(\ell_a^{\odot k_a})$ by Fubini's theorem. The identity $\Hm^{2k}(\A(V))=\int_VJ(\Log|_V)N(\Log z)^{-1}\,d\mathrm{vol}(z)$ with $N=2^k$ almost everywhere gives $\Hm^{2k}(\A(V))=2^{-k}\MVol(V)$, and the product formula follows from Theorem \ref{thm:square} for each factor; the value $\pi^{2k}/2^k$ follows from Corollary \ref{cor:squareproj}. Finally, write $V=V_A\times V_B$ with $A$ the first block of coordinates, $V_A=\ell_1^{\odot k_1}$, and $B$ the others, and parametrise $V$ by $t=(t_A,t_B)\in\C^{k_1}\times\C^{k-k_1}$. The coordinates $z_j$, $j\in A$, do not depend on $t_B$, so the $n_1+1\ge2k_1+1$ points $p_0$ and $p_j$, $j\in A$, lie in the projective subspace $\{$coefficients of $t_B$ vanish$\}$, of dimension $k_1$; hence $k_1+2$ of them are linearly dependent, and $k_1+2\le k+1$ when $r\ge2$.
\end{proof}

The last argument shows more generally that a $k$-plane which is a product $V_A\times V_B$ of planes of positive dimensions $k_A,k_B$ in disjoint sets $A,B$ of coordinates has its dual points in two projective subspaces $\Lambda_A\ni p_0$ and $\Lambda_B\ni p_0$ of dimensions $k_A$ and $k_B$ with $\Lambda_A\cap\Lambda_B=\{p_0\}$, namely $\Lambda_A=\{$coefficients of $t_B$ vanish$\}$ and $\Lambda_B=\{$coefficients of $t_A$ vanish$\}$, with $p_j\in\Lambda_A$ for $j\in A$ and $p_j\in\Lambda_B$ for $j\in B$.

\subsection{Fibre extensions}

Products of symmetric powers are not the only real $k$-planes with $2^k$-to-one amoeba maps. The reason is that the fibres of the amoeba map of a line of $(\C^*)^2$ depend only on the distance between the two points where it meets the coordinate axes, and that such distances can be read off from the base in a way which is invariant under the partial conjugations.

Let $W$ be a product of symmetric powers of real lines as above, of dimension $k_0$, in $\C^{n_0}$, with parameters $q^{(a)}$. We call a function $\kappa$ on $W$ admissible if it is either a nonzero real constant or of the form $\kappa(w)=\mu\prod_{i=1}^{k_a}\big(q^{(a)}_i-\varpi'\big)=\mu\,(-1)^{k_a}f_{q^{(a)}}(\varpi')$ for some factor $a$, some $\mu\in\R^*$ and some real $\varpi'\notin\{\varpi^{(a)}_1,\dots,\varpi^{(a)}_{n_a}\}$. An admissible function is a real affine function of the point $w\in W$, since the coefficients of $f_{q^{(a)}}$ are affine coordinates on the factor, and its modulus is invariant under the conjugation of any single entry $q^{(a)}_i$, since $|\bar q-\varpi'|=|q-\varpi'|$. For admissible functions $\kappa_1,\dots,\kappa_b$ we define the fibred extension

$$V=\big\{(w,\,s_1,\,s_1-\kappa_1(w),\,\dots,\,s_b,\,s_b-\kappa_b(w)):\ w\in W,\ s\in\C^b\big\}\subset\C^{n_0}\times\C^{2b},$$

a real affine plane of dimension $k=k_0+b$ in $\C^{n_0+2b}$, half-dimensional when $n_0=2k_0$. When all $\kappa_j$ are constant, $V$ is the product of $W$ with $b$ real lines of $(\C^*)^2$.

\begin{theorem}\label{thm:fibred}
Let $V$ be a fibred extension of a product of symmetric powers $W$ as above. Then the amoeba map of $V$ is exactly $2^k$-to-one at $\Hm^{2k}$-almost every point of $\A(V)$. More precisely, if $(w,s)\in V$ is such that the fibre of the amoeba map of $W$ through $w$ has $2^{k_0}$ points and $s_j\notin\R\kappa_j(w)$ for every $j$, then the fibre through $(w,s)$ consists of the points $(w',s')$ with $w'$ in the fibre of $w$ and $s'_j$ one of the two intersection points of the circles $\{|\sigma|=|s_j|\}$ and $\{|\sigma-\kappa_j(w')|=|s_j-\kappa_j(w)|\}$.
\end{theorem}

\begin{proof}
A point $(w',s')$ of $V$ has the same image under $\Log$ as $(w,s)$ if and only if $\Log w'=\Log w$, $|s'_j|=|s_j|$ and $|s'_j-\kappa_j(w')|=|s_j-\kappa_j(w)|$ for every $j$. The first condition says that $w'$ lies in the fibre of $w$, an orbit of the conjugations of the single entries of the parameters of $W$ (Proposition \ref{prop:prodfib}), and since the $\kappa_j$ are admissible, $|\kappa_j(w')|=|\kappa_j(w)|$. For fixed $w'$ and $j$, the remaining conditions say that $s'_j$ lies on the circle of radius $\rho_1=|s_j|$ centred at $0$ and on the circle of radius $\rho_2=|s_j-\kappa_j(w)|$ centred at $\kappa_j(w')$. Two circles whose centres are at distance $d>0$ meet in exactly two points if $|\rho_1-\rho_2|<d<\rho_1+\rho_2$, in one point if equality holds on one side, and in no point otherwise. For $w'=w$ the point $s_j$ lies on both circles, and the triangle with vertices $0$, $\kappa_j(w)$, $s_j$ is nondegenerate because $s_j\notin\R\kappa_j(w)$, so the strict inequalities hold with $d=|\kappa_j(w)|$; since they involve only $\rho_1,\rho_2$ and $d=|\kappa_j(w')|=|\kappa_j(w)|$, they hold for every $w'$ in the fibre. Hence each of the $2^{k_0}$ points $w'$ and each $j$ contribute exactly two values of $s'_j$, and the fibre has $2^{k_0}\cdot2^b=2^k$ points. The excluded points form the union of the preimage of the exceptional set of $W$ and of the sets $\{s_j\in\R\kappa_j(w)\}$, which have real codimension at least one in $V$; their images under the locally Lipschitz map $\Log$ are $\Hm^{2k}$-null \cite[Section 2.10.2]{Fed69}, and every point of $\A(V)$ outside these images has a fibre of exactly $2^k$ points, because a fibre containing one excluded point consists of excluded points, the conditions being invariant under the substitutions above.
\end{proof}

The fibred extensions with a nonconstant $\kappa_j$ are not products of symmetric powers. The simplest example is the real plane

$$V_\mu=\{(q,\,q-1,\,s,\,s-\mu(q-\varpi')):\ q,s\in\C\}\subset\C^4,\qquad\mu\in\R^*,\ \varpi'\in\R\setminus\{0,1\},$$

a fibred extension of the real line $W=\{(q,q-1)\}$ with $b=1$. Its dual points in $(\mathbb P^2)^*$, with coordinates dual to $(1,q,s)$, are $p_0=[1,0,0]$, $p_1=[0,1,0]$, $p_2=[-1,1,0]$, $p_3=[0,0,1]$ and $p_4=[\mu\varpi',-\mu,1]$. The points $p_0,p_1,p_2$ are collinear, so $V_\mu$ is not in general position and is not a square of a line. It is not a product either: a decomposition into two lines in disjoint pairs of coordinates would give two lines $\Lambda_A,\Lambda_B$ through $p_0$, meeting only at $p_0$, containing two of the points $p_1,\dots,p_4$ each; the line through $p_0$ and $p_1$ contains $p_2$ and neither $p_3$ nor $p_4$, while the line through $p_0$ and $p_3$ does not contain $p_4$ since $\mu\ne0$, so no such decomposition exists. By Theorem \ref{thm:fibred} the amoeba map of $V_\mu$ is four-to-one. In higher dimension, the same argument shows that the real $3$-plane

$$\{(q,\,q-1,\,s_1,\,s_1-\mu_1(q-\varpi'_1),\,s_2,\,s_2-\mu_2(q-\varpi'_2))\}\subset\C^6,\qquad\mu_1\mu_2\ne0,$$

has an eight-to-one amoeba map and is not a product of symmetric powers. Its dual points in $(\mathbb P^3)^*$, with coordinates dual to $(1,q,s_1,s_2)$, are $p_0=[1,0,0,0]$, $p_1=[0,1,0,0]$, $p_2=[-1,1,0,0]$, $p_3=[0,0,1,0]$, $p_4=[\mu_1\varpi'_1,-\mu_1,1,0]$, $p_5=[0,0,0,1]$, $p_6=[\mu_2\varpi'_2,-\mu_2,0,1]$; if $V$ were a product $V_A\times V_B$, the subspaces $\Lambda_A,\Lambda_B$ of the remark after Proposition \ref{prop:prodfib} would separate these points, but $p_1$ and $p_2$ must lie in the same one, say $\Lambda_A$, since $p_2$ lies on the line through $p_0$ and $p_1$. If $p_3\in\Lambda_B$, then either $p_4\in\Lambda_A$, and $p_3=p_4-\mu_1\varpi'_1p_0+\mu_1p_1$ (as vectors) lies in $\Lambda_A\cap\Lambda_B$, or $p_4\in\Lambda_B$, and $\mu_1p_1=-(p_4-p_3-\mu_1\varpi'_1p_0)$ lies in $\Lambda_A\cap\Lambda_B$; both are impossible, so $p_3\in\Lambda_A$, and then $p_4\in\Lambda_A$ by the same relation. The same argument with $\mu_2\ne0$ puts $p_5,p_6$ in $\Lambda_A$, so the decomposition is trivial. Fibre extensions can be iterated only in a restricted way: if a further line were attached through a function of the variable $s_j$, the two points $s'_j$ of the fibre above, which are exchanged by the reflection in the line $\R\kappa_j(w)$ rather than by a conjugation, would give in general different distances, and the count $2^k$ would be lost.

In particular the answer to the question at the beginning of this section is negative as soon as general position is not required, already for $k=2$. The classification of the real $k$-planes with $2^k$-to-one amoeba maps without any condition on the position of the hyperplanes requires a study of the degenerate configurations, which we carry out for $k=2$ and, on the stratum of configurations with a collinear triple, for $k=3$ in Section \ref{sec:degenerate}.

\subsection{The classification in general position}

We now prove that in general position the symmetric powers are the only real $k$-planes with $2^k$-to-one amoeba maps. In fact we prove a stronger statement, in which only the inequality $N\ge2^k$ is assumed, for all $n\ge2k$; this form is needed for the induction on $k$. For a real $p$-plane $P\subset\C^N$ we write $N_P(x)=\#(P\cap\Log^{-1}(x))$.

\begin{theorem}\label{thm:classif}
Let $p\ge2$, $N\ge2p$, and let $P\subset\C^N$ be a real affine $p$-plane in general position. Then $N_P(x)\ge2^p$ for $\Hm^{2p}$-almost every $x\in\A(P)$ if and only if $P$ is a multiplicative translate of the $p$-th symmetric power of a real line, that is, if and only if its $N+1$ dual points lie on a rational normal curve of degree $p$; in that case the amoeba map of $P$ is exactly $2^p$-to-one almost everywhere.
\end{theorem}

The proof rests on four lemmas. The first passes from $N>2p$ to the half-dimensional case through coordinate projections, the second is a statement on rational normal curves, the third controls the amoeba map near a coordinate hyperplane, and the fourth deduces from it the behaviour of the hyperplane sections.

\begin{lemma}\label{lem:projfib}
Let $P\subset\C^N$ be a real $p$-plane, $N>2p$, and $|I|=2p$ with $\pi_I|_P$ injective. If $N_P\ge c$ almost everywhere on $\A(P)$, then $N_{\pi_I(P)}\ge c$ almost everywhere on $\A(\pi_I(P))$.
\end{lemma}

\begin{proof}
Since $\pi_I|_P$ is injective, $N_{\pi_I(P)}(y)$ is the number $N_I(y)$ of points of $P$ in $\Log_I^{-1}(y)$. Let $E=\{x\in\A(P):N_P(x)<c\}$, an $\Hm^{2p}$-null set; then $p_I(E)$ is Lebesgue null, since $p_I$ is Lipschitz. For $y\in\A(\pi_I(P))\setminus p_I(E)$ there is at least one point of $\A(P)$ over $y$, and by \eqref{eq:NS}, $N_I(y)=\sum_{x}N_P(x)\ge c$, the sum running over the points $x\in\A(P)\cap p_I^{-1}(y)$.
\end{proof}

\begin{lemma}\label{lem:cones}
Let $k\ge3$ and let $p_0,\dots,p_{2k}$ be points of $\mathbb P^k$ in general linear position. Suppose that for two indices $a\ne b$ the projections from $p_a$ and from $p_b$ map the other $2k$ points onto rational normal curves of degree $k-1$ of $\mathbb P^{k-1}$. Then all the points lie on a rational normal curve of degree $k$.
\end{lemma}

\begin{proof}
The projected points are in general linear position in $\mathbb P^{k-1}$, since $k$ of them are dependent only if, together with the centre of projection, $k+1$ of the original points are dependent. Let $C$ be the rational normal curve of degree $k$ through $p_a$, $p_b$ and $k+1$ further points of the configuration, which exists and is unique by \cite[Theorem 1.18]{Har92}, since these are $k+3$ points in general position. The projection $\pi_a$ from $p_a$ maps $C$ onto a rational normal curve $C_a$ of degree $k-1$, which passes through the $k+2$ projections of the other chosen points; by hypothesis there is a rational normal curve $D_a$ of degree $k-1$ through the projections of all $2k$ points, and since $k+2=(k-1)+3$ points in general position of $\mathbb P^{k-1}$ lie on a unique rational normal curve of degree $k-1$, $D_a=C_a$. Hence every point $p$ of the configuration lies on the cone $\Sigma_a$ with vertex $p_a$ over $C$, the union of the lines joining $p_a$ to the points of $C$, together with the tangent line of $C$ at $p_a$; similarly $p\in\Sigma_b$. We claim that $\Sigma_a\cap\Sigma_b=C\cup\overline{p_ap_b}$. Let $x\in\Sigma_a\cap\Sigma_b$ not on the line $\overline{p_ap_b}$. Then the line $\overline{p_ax}$ meets $C$ in a point $c_a$ other than $p_a$, or is tangent to $C$ at $p_a$, and similarly the line $\overline{p_bx}$ meets $C$ in $c_b$ or is tangent at $p_b$. The two lines lie in the plane spanned by $p_a,p_b,x$, so this plane contains the divisor $p_a+c_a+p_b+c_b$ of $C$, where $c_a=p_a$ means the tangent line, of degree $4$. Any divisor of degree at most $k+1$ on a rational normal curve of degree $k$ spans a projective space of dimension one less than its degree, because the corresponding vectors $(1,t,\dots,t^k)$ and their derivatives form a confluent Vandermonde matrix of maximal rank; since $4\le k+1$, the divisor cannot lie in a plane unless two of its points coincide in a way which lowers the degree, that is, unless $c_a=c_b$. If $c_a=c_b=c$, then $c$ lies on the two distinct lines $\overline{p_ax}$ and $\overline{p_bx}$, which meet only at $x$, so $x=c\in C$; the cases $c_a=p_b$ or $c_b=p_a$ would put $x$ on the line $\overline{p_ap_b}$. This proves the claim. No point of the configuration other than $p_a,p_b$ lies on $\overline{p_ap_b}$, by general position, so all of them lie on $C$.
\end{proof}

For the third lemma let $V\subset\C^{2k}$ be a real $k$-plane in general position, $k\ge2$, fix an index $m$, and let $H=V\cap\{z_m=0\}$, regarded as a real $(k-1)$-plane of $\C^{2k-1}$ in the coordinates $z'=(z_j)_{j\ne m}$. Its dual points are the projections from $p_m$ of the points $p_j$, $j\ne m$, since the dual of the hyperplane $\{z_m=0\}$ of $\overline V$ is the quotient of $(\mathbb P^k)^*$ by $p_m$; hence $H$ is in general position. Choose a real vector $u$ in the direction space of $V$ with $u_m=1$, so that every point of $V$ is $h+tu$ with $h\in H$ and $t=z_m\in\C$. The map $\Log|_H:H^\circ\to\R^{2k-1}$ is proper, since $\Log$ is proper on $(\C^*)^{2k-1}$ and $H^\circ$ is closed there. At a point $h\in H^\circ$ where it is an immersion, let $\eta(h)\in\R^{2k-1}$ be a unit normal vector of the hypersurface $\Log(U)$, $U$ a small neighbourhood of $h$, at $\Log h$, and put $\chi(h)=\sum_{j\ne m}\eta_j(h)\,\frac{u_j}{h_j}\in\C$. Replacing $u$ by $u+v$ with $v$ in the direction space of $H$ does not change $\chi$, because $\eta(h)$ is orthogonal to $\operatorname{Re}(e^{i\phi}v'/h')$ for every $\phi$, where $v'/h'=(v_j/h_j)_{j\ne m}$, that is $\sum_j\eta_jv_j/h_j=0$; and $|\chi(h)|$ does not depend on the sign of $\eta$.

\begin{lemma}\label{lem:chi}
The map $\Log|_H$ is an immersion on a dense open subset of $H^\circ$, and $\chi\ne0$ on a dense open subset of it whose complement is null.
\end{lemma}

\begin{proof}
Let $V'\subset\C^{2k-1}$ be the image of $V$ under the projection forgetting $z_m$, which is injective, and let $D$ be its direction space, the complexification of a real $k$-dimensional space $D_\R$; then $D$ is spanned by the direction space of $H$ and by $u'$. Call a point $z\in V'^\circ$ critical if there is a real vector $\eta\ne0$ with $\sum_j\eta_jw_j/z_j=0$ for all $w\in D$, that is, if the differential of $\Log|_{V'}$ at $z$ has rank less than $2k-1$. If $h\in H^\circ$ is a point where $\Log|_H$ is an immersion and $\chi(h)=0$, then $\eta(h)$ annihilates both the direction space of $H$ divided by $h$ and $u'/h'$, so $h$ is critical. Criticality is the vanishing of the $(2k-1)\times(2k-1)$ minors of the real matrix of $w\mapsto\operatorname{Re}(w/z)$ on $D_\R\otimes\C$, a real analytic condition on the connected set $H^\circ$, the complement of a hyperplane arrangement in $H$; hence either every point of $H^\circ$ is critical, or the critical points of $H^\circ$ form a proper real analytic subset, which is closed, nowhere dense and null. The same applies to the points where $\Log|_H$ fails to be an immersion. It therefore suffices to exhibit points of $H^\circ$ which are not critical for $V'$ and at which $\Log|_H$ is an immersion.

Choose a set $K$ of $k-1$ indices different from $m$ and let $L$ be the set of the $k$ remaining indices different from $m$. By general position, the point $p=H\cap\{z_j=0,\ j\in K\}$ is a single real point of $\C^{2k-1}$, its coordinates $p_j$, $j\in L$, are nonzero, the coordinates $(z_j)_{j\in K}$ are affine coordinates on $H$, and the projection $D\to\C^L$, $w\mapsto(w_j)_{j\in L}$, is an isomorphism, as are $D_\R\to\R^L$ and the projection of the direction space of $H$ to $\C^K$. For small $\varepsilon>0$ let $h(\varepsilon)\in H$ be the point with $h_j(\varepsilon)=i\varepsilon$ for $j\in K$; then $h(\varepsilon)\in H^\circ$ and $h_j(\varepsilon)=p_j+O(\varepsilon)$ for $j\in L$. For $j_0\in L$ let $w^{(j_0)}\in D_\R$ be the vector with $w^{(j_0)}_j=\delta_{jj_0}$ for $j\in L$, and put $a_{ij_0}=w^{(j_0)}_i$ for $i\in K$, a real $K\times L$ matrix of rank $k-1$, because $D_\R\to\R^K$ is onto. If $h(\varepsilon)$ is critical with a unit vector $\eta=\eta(\varepsilon)$, the conditions for the vectors $w^{(j_0)}$ read

$$\frac{\eta_{j_0}}{h_{j_0}(\varepsilon)}=-\sum_{i\in K}\frac{\eta_ia_{ij_0}}{i\varepsilon}=\frac{i}{\varepsilon}\sum_{i\in K}\eta_ia_{ij_0},\qquad j_0\in L .$$

The right-hand side is purely imaginary, while the real part of the left-hand side is $\eta_{j_0}/p_{j_0}+O(\varepsilon)$; hence $\eta_{j_0}=O(\varepsilon)$ for $j_0\in L$. Then the imaginary part of the left-hand side is $\eta_{j_0}\operatorname{Im}(1/h_{j_0}(\varepsilon))=O(\varepsilon^2)$, because $\operatorname{Im}(1/h_{j_0}(\varepsilon))=O(\varepsilon)$, so $\sum_{i\in K}\eta_ia_{ij_0}=O(\varepsilon^3)$ for all $j_0\in L$, and since the matrix $(a_{ij_0})$ has rank $k-1=|K|$, $\eta_i=O(\varepsilon^3)$ for $i\in K$. This contradicts $|\eta|=1$ for small $\varepsilon$, so $h(\varepsilon)$ is not critical. At the same points $\Log|_H$ is an immersion: if $v$ is in the direction space of $H$ and $\operatorname{Re}(v_j/h_j(\varepsilon))=0$ for all $j$, then $v_j/(i\varepsilon)$ is purely imaginary for $j\in K$, so $(v_j)_{j\in K}$ is real, hence $v$ is real, being determined by its coordinates in $K$; then $\operatorname{Re}(v_j/h_j(\varepsilon))=v_j\operatorname{Re}(1/h_j(\varepsilon))$ with $\operatorname{Re}(1/h_j(\varepsilon))=1/p_j+O(\varepsilon)\ne0$ for $j\in L$, so $v_j=0$ for $j\in L$, and $v=0$ because the $k$ coordinates in $L$ cannot all vanish on a nonzero vector of the $(k-1)$-dimensional direction space of $H$, by general position.
\end{proof}

\begin{lemma}\label{lem:tentacle}
Let $V$, $H$ and $m$ be as above, and let $y\in\A(H)$ be such that all points $h_1,\dots,h_M$ of $H^\circ$ over $y$ are points where $\Log|_H$ is an immersion with $\chi(h_\alpha)\ne0$, the images of their neighbourhoods having the same germ $S$ at $y$, a smooth hypersurface of $\R^{2k-1}$ with unit normal $\eta$ at $y$. Then there are a neighbourhood $Y$ of $y$ in $S$ and $r_0>0$ such that, for $y'\in Y$, $0<r<r_0$ and real $d$ with $|d|<r\min_\alpha|\chi(h_\alpha)|/2$, the fibre of the amoeba map of $V$ over the point $x$ with $x_m=\log r$ and $(x_j)_{j\ne m}=y'+d\,\eta(y')$ has exactly $2M$ points. In particular $N_V=2M$ on a set of positive Lebesgue measure.
\end{lemma}

\begin{proof}
Let $U$ be a ball around $y$ in $\R^{2k-1}$. A point of $V$ over $x$ is $h+tu$ with $h\in H$, $|t|=r$ and $\Log(h'+tu')=x'$, where $x'=(x_j)_{j\ne m}$. If $x'\in U$, then $h'$ lies at distance $O(r)$ from the compact set $\Log^{-1}(\overline U)$, so $\Log h'\in U$ for $r$ small, uniformly; since $\Log|_H$ is proper with finite fibre over $y$, after shrinking $U$ every such $h$ lies in a small neighbourhood $U_\alpha$ of one of the $h_\alpha$. Fix $\alpha$, and let $\varphi:\R^{2k-2}\supset\Omega\to U_\alpha$ be a local parametrisation of $H$ with $\varphi(0)=h_\alpha$. In tubular coordinates $x'=\sigma(\xi)+d\,\eta(\sigma(\xi))$ around $S$, with $\sigma$ a parametrisation of $S$ by $\xi\in\R^{2k-2}$, the map $(\omega,\theta,r)\mapsto\Log(\varphi(\omega)'+re^{i\theta}u')$ is smooth up to $r=0$ and equals $\Log\varphi(\omega)'+r\operatorname{Re}\big(e^{i\theta}u'/\varphi(\omega)'\big)+O(r^2)$, the error being smooth with its derivatives. Its tubular coordinates are therefore $\xi=\xi_0(\omega)+O(r)$, where $\xi_0$ is the local diffeomorphism $\sigma^{-1}\circ\Log\circ\varphi$, and $d=r\big(\operatorname{Re}(e^{i\theta}\chi(\varphi(\omega)))+O(r)\big)$, with $\eta$ oriented consistently. Solving the first equation for $\omega$ by the implicit function theorem, the points of $V$ near $h_\alpha$ over $x$ correspond to the solutions $\theta$ of $d/r=\operatorname{Re}(e^{i\theta}\chi_\alpha(\xi,\theta,r))+O(r)$, where $\chi_\alpha(\xi,\theta,r)=\chi(h_\alpha)+o(1)$ uniformly as $y'\to y$ and $r\to0$, with the same control on the $\theta$-derivative. For $|d/r|<|\chi(h_\alpha)|/2$ and $r$ small, the function $\theta\mapsto|\chi(h_\alpha)|\cos(\theta+\arg\chi(h_\alpha))$ takes the value $d/r$ at exactly two points, at which its derivative is at least $\frac{\sqrt3}2|\chi(h_\alpha)|$ in absolute value, and a $C^1$-small perturbation keeps exactly two solutions. Summing over $\alpha$ gives $2M$ points. The set of such $x$ is the image of $\{(y',d,\log r)\}$ under a diffeomorphism onto an open subset of $\R^{2k}$, so it has positive measure.
\end{proof}

\begin{proposition}\label{prop:sections}
Let $k\ge2$, let $V\subset\C^{2k}$ be a real $k$-plane in general position with $N_V\ge2^k$ almost everywhere on $\A(V)$, and let $m$ be a coordinate index. Then $N_H\ge2^{k-1}$ almost everywhere on $\A(H)$ for $H=V\cap\{z_m=0\}$, a real $(k-1)$-plane in general position in $\C^{2k-1}$.
\end{proposition}

\begin{proof}
By Lemma \ref{lem:chi} and the area formula, the set $Z\subset H^\circ$ of points where $\Log|_H$ is not an immersion or $\chi=0$ is null, so its image is $\Hm^{2k-2}$-null, and so is the image of $H(\R)$, of dimension $k-1$. For $\Hm^{2k-2}$-almost every $y\in\A(H)$ the preimages avoid $Z$, and the images of their neighbourhoods have the same germ at $y$: two real analytic hypersurface germs through $y$ which do not coincide meet in a set of dimension at most $2k-3$, and there are countably many pairs of such germs to consider by the properness of $\Log|_H$ and the second countability of $H$. If $N_H(y)=M<2^{k-1}$ for such a $y$, Lemma \ref{lem:tentacle} gives a set of positive measure on which $N_V=2M<2^k$, a contradiction.
\end{proof}

\begin{proposition}\label{prop:planesN}
Let $N\ge5$ and let $P\subset\C^N$ be a real plane in general position. If $N_P\ge4$ almost everywhere on $\A(P)$, then $P$ is a multiplicative translate of the square of a real line.
\end{proposition}

\begin{proof}
The proof of Theorem \ref{thm:conic} applies verbatim to $N\ge5$: if the $N+1$ dual points of $P$ lie on no conic, then $\Log z=\Log z'$ implies $z'\in\{z,\bar z\}$, so $N_P=2$ outside the image of $P(\R)$, which is $\Hm^4$-null, contradicting the hypothesis since $\A(P)$ has positive measure. Hence the dual points lie on a conic, which is smooth because no three of them are collinear, and the construction of Section \ref{sec:conic} writes $P$ as a multiplicative translate of $\ell\odot\ell$ for a real line $\ell$.
\end{proof}

\begin{proof}[Proof of Theorem \ref{thm:classif}]
The symmetric powers satisfy $N_P=2^p$ almost everywhere by Theorem \ref{thm:square}, and their dual points lie on a rational normal curve, so only the direct implication requires proof. We prove by induction on $p\ge2$ the statement $\mathrm{A}(p)$, that the direct implication holds for $N=2p$, together with the statement $\mathrm{B}(p)$, that it holds for every $N>2p$.

For $p=2$, $\mathrm A(2)$ is Corollary \ref{cor:square4}, which applies to every real plane in general position of $\C^4$, and $\mathrm B(2)$ is Proposition \ref{prop:planesN}.

Let $p\ge3$ and assume $\mathrm B(p-1)$. Let $V\subset\C^{2p}$ be a real $p$-plane in general position with $N_V\ge2^p$ almost everywhere. For each coordinate index $m$, Proposition \ref{prop:sections} shows that $H_m=V\cap\{z_m=0\}$, a real $(p-1)$-plane in general position in $\C^{2p-1}$, satisfies $N_{H_m}\ge2^{p-1}$ almost everywhere, and $2p-1>2(p-1)$, so by $\mathrm B(p-1)$ it is a multiplicative translate of a symmetric power, that is, its $2p$ dual points, the projections from $p_m$ of the points $p_j$, $j\ne m$, lie on a rational normal curve of degree $p-1$. Applying this to two coordinate indices $a\ne b$, Lemma \ref{lem:cones} shows that the $2p+1$ dual points of $V$ lie on a rational normal curve of degree $p$. It is real, being unique through $p+3$ of the real points, and its real points have real parameters in a real parametrisation; moving the parameter of $p_0$ to $\infty$ by a real projective transformation, the construction of Section \ref{sec:squares} shows that $V$ is a multiplicative translate of a symmetric power of a real line. This proves $\mathrm A(p)$.

Assume $\mathrm A(p)$ with $p\ge3$, and let $P\subset\C^N$, $N>2p$, be a real $p$-plane in general position with $N_P\ge2^p$ almost everywhere. For every $I$ with $|I|=2p$, the projection $\pi_I|_P$ is injective and $\pi_I(P)$ is a real $p$-plane in general position in $\C^{2p}$, whose dual points are $p_0$ and $p_j$, $j\in I$. By Lemma \ref{lem:projfib}, $N_{\pi_I(P)}\ge2^p$ almost everywhere, and by $\mathrm A(p)$ the points $p_0$ and $p_j$, $j\in I$, lie on a rational normal curve $C_I$ of degree $p$. Fix a set $T$ of $p+2$ coordinate indices. Since $2p\ge p+3$, for every index $j\notin T$ there is $I\supset T\cup\{j\}$ with $|I|=2p$, and $C_I$ passes through the $p+3$ points $p_0$ and $p_i$, $i\in T$, so it is the unique rational normal curve $C_T$ through them \cite[Theorem 1.18]{Har92}. Hence all dual points of $P$ lie on $C_T$, and as before $P$ is a multiplicative translate of a symmetric power. This proves $\mathrm B(p)$ and completes the induction.
\end{proof}

\begin{corollary}\label{cor:fibre4}
Let $k\ge3$ and let $V\subset\C^{2k}$ be a real $k$-plane in general position which is not a multiplicative translate of a symmetric power of a real line. Then the amoeba map of $V$ has fewer than $2^k$ points, and a number of points divisible by $4$, in every fibre over a set of positive measure. For $k=3$ these fibres have exactly $4$ points; more precisely, for every coordinate index $m$ except at most one, the section $H_m=V\cap\{z_m=0\}$ is not a square of a line, its amoeba map is two-to-one, and the fibres have exactly $4$ points near the tentacle of $\A(V)$ in the direction $-e_m$ described by Lemma \ref{lem:tentacle}.
\end{corollary}

\begin{proof}
The first assertion is the contrapositive of Theorem \ref{thm:classif}, together with Corollary \ref{cor:real4}. For $k=3$, if two coordinate sections $H_a,H_b$ were squares, Lemma \ref{lem:cones} would make $V$ a symmetric cube. For any other $m$, the six dual points of $H_m$, a real plane in general position in $\C^5$, lie on no conic, since a conic through them would be smooth and $H_m$ a square by the construction of Section \ref{sec:conic}; by Theorem \ref{thm:conic} its amoeba map is two-to-one outside an $\Hm^4$-null set, and Lemma \ref{lem:tentacle} with $M=2$, applied at almost every point of $\A(H_m)$ as in the proof of Proposition \ref{prop:sections}, gives fibres of exactly $4$ points on a set of positive measure.
\end{proof}

By Corollary \ref{cor:fibre4}, for every real generic $k$-plane of $(\C^*)^{2k}$ with $k\ge3$, the amoeba map is not $2^k$-to-one, contrary to \cite[Proposition 5.6]{NP17} and to the lower bound $2^k$ of \cite[Proposition 5.2]{NP17} used in its proof, and for $k=3$ fibres of exactly $4$ points occur on a set of positive measure, which solves the case $k=3$ of the second problem of Section \ref{sec:open}. The volume formula \cite[Theorem 1.2(ii)]{NP17} is not decided by this argument, as explained below. It also gives a precise meaning to the heuristic description of Section \ref{sec:why}: the tentacles of the amoeba of $V$ are thickenings of the amoebas of the hyperplane sections, and the fibre over a point of a tentacle has twice as many points as the fibre of the section, so the count $2^k$ can only be inherited from sections which themselves have $2^{k-1}$-to-one amoeba maps. For $k=2$ the statement of Corollary \ref{cor:fibre4} is empty, in accordance with Corollary \ref{cor:square4}, and for $k\ge4$ the fibres of fewer than $2^k$ points have $4$, $8,\dots$ or $2^k-4$ points; we do not know whether fibres of exactly $4$ points always occur.

Theorem \ref{thm:classif} does not decide which volumes the amoebas of the other real $k$-planes have, since fibres of more than $2^k$ points may also occur; but it shows that among real $k$-planes in general position the value $\pi^{2k}/2^k$ of Corollary \ref{cor:squareproj}, obtained for symmetric powers with every fibre of $2^k$ points, cannot be obtained in the same way by any other plane.

\subsection{Planes not in general position}\label{sec:degenerate}

We now turn to real $k$-planes whose dual points are not in general position, and classify those with $2^k$-to-one amoeba maps completely for $k=2$ and on the main stratum of configurations for $k=3$. The classification is naturally made up to a symmetry which exchanges the hyperplane at infinity with the coordinate hyperplanes.

For $a\in\{1,\dots,n\}$ let $\mu_a:(\C^*)^n\to(\C^*)^n$ be the map with $\mu_a(z)_a=1/z_a$ and $\mu_a(z)_j=z_j/z_a$ for $j\ne a$. In the homogeneous coordinates $(Z_0:\dots:Z_n)$ of $\mathbb P^n$, with $z_j=Z_j/Z_0$, it is the change of affine chart from $Z_0\ne0$ to $Z_a\ne0$, followed by the exchange of the names of $Z_0$ and $Z_a$.

\begin{lemma}\label{lem:monomial}
The map $\mu_a$ is a group automorphism of $(\C^*)^n$ with real coefficients, and $\Log\circ\mu_a=L_a\circ\Log$, where $L_a(x)_a=-x_a$ and $L_a(x)_j=x_j-x_a$ for $j\ne a$ is a linear automorphism of $\R^n$ with determinant $-1$. If $V\subset\C^n$ is a real affine $k$-plane on which $z_a$ is not constant, then $\mu_a(V\cap(\C^*)^n)$ is the intersection with the torus of a real affine $k$-plane $V^{(a)}$, whose dual points are those of $V$ with $p_0$ and $p_a$ exchanged, and $N_{V^{(a)}}(L_ax)=N_V(x)$ for every $x\in\R^n$. In particular, $\A(V)$ has positive $\Hm^{2k}$-measure if and only if $\A(V^{(a)})$ has, and the amoeba map of $V$ is $c$-to-one almost everywhere if and only if that of $V^{(a)}$ is.
\end{lemma}

\begin{proof}
The first assertions are immediate, $\log|1/z_a|=-\log|z_a|$ and $\log|z_j/z_a|=\log|z_j|-\log|z_a|$. The projective closure $\overline V\subset\mathbb P^n$ is a projective $k$-space defined over $\R$, and $\mu_a$ is the restriction to the torus of the projective transformation which exchanges $Z_0$ and $Z_a$; the image of $\overline V$ is again a real projective $k$-space, not contained in $\{Z_0=0\}$ because $z_a$ is not constant on $V$, and its intersection with $\{Z_0\ne0\}$ is the real affine $k$-plane $V^{(a)}$. The hyperplanes of $\overline V$ cut by the coordinate hyperplanes are the same, only their labels $0$ and $a$ being exchanged, so the dual points are exchanged accordingly. Since $\mu_a$ is a bijection of the torus mapping $V\cap(\C^*)^n$ onto $V^{(a)}\cap(\C^*)^n$ and the fibre of $\Log$ over $x$ onto the fibre over $L_ax$, the fibre counts agree; since $L_a$ is a linear isomorphism, it preserves the classes of $\Hm^{2k}$-null sets on $2k$-rectifiable sets \cite[Section 2.10.11]{Fed69}.
\end{proof}

We call two real affine $k$-planes of $\C^n$ monomially equivalent if they are related by a finite sequence of maps $\mu_a$, permutations of the coordinates and multiplicative translations by vectors of $(\R^*)^n$. Monomial equivalence permutes the $n+1$ dual points arbitrarily, and preserves fibre counts, reality and the nondegeneracy of the amoeba. The configuration of dual points up to projective transformations, and in particular its pattern of collinear and coplanar subsets, is an invariant.

\begin{theorem}\label{thm:planes4}
Let $V\subset\C^4$ be a real affine plane with $\Hm^4(\A(V))>0$. Then its five dual points are distinct, no four of them are collinear, at most two triples of them are collinear, and $V$ is monomially equivalent to exactly one of the following: (a) if no three dual points are collinear, the square $\ell\odot\ell$ of a real line; (b) if exactly three are collinear, the fibred extension $V_\mu=\{(q,\,q-1,\,s,\,s-\mu(q-\varpi'))\}$ with $\mu\ne0$ and $\varpi'\notin\{0,1\}$; (c) if two triples are collinear, the product of two real lines of $(\C^*)^2$. In every case the amoeba map of $V$ is exactly four-to-one almost everywhere, $\MVol(V)=\pi^4$ and $\Vol(\A(V))=\pi^4/4$.
\end{theorem}

\begin{proof}
If $p_i=p_j$ with $i,j\ge1$, then $z_i=cz_j$ on $V$ for a real $c$, and $\A(V)$ lies in the hyperplane $\{x_i-x_j=\log|c|\}$; if $p_j=p_0$, then $z_j$ is constant on $V$; in both cases $\Hm^4(\A(V))=0$. Five points spanning the plane cannot be collinear. If four of them are collinear, Lemma \ref{lem:monomial} lets us assume that $p_0$ is one of them, say $p_0,p_1,p_2,p_3$; then $z_1,z_2,z_3$ are functions of one real affine form $q$ on $V$, and taking $s=z_4$ as the second coordinate, $V$ is the product of a real line $\ell\subset\C^3$ with $\C$, so $\A(V)=\A(\ell)\times\R$ has dimension $3$, a contradiction. Two distinct collinear triples share exactly one point, since they lie on two distinct lines and no four points are collinear; if there were three triples, the third would have to contain points of both other lines outside their common point, and then two of its points would lie on one of these lines, which would contain four points. Hence there are at most two triples.

In case (c), Lemma \ref{lem:monomial} lets us move the common point of the two triples to $p_0$; then $z_1,z_2$ depend on one real affine form $q$ and $z_3,z_4$ on another form $s$, independent of $q$, so $V$ is a product of two real lines, and Proposition \ref{prop:prodfib} applies. In case (b), we move one point of the triple to $p_0$, so that $z_1,z_2$ depend only on a real affine form $q$, which we normalise, by an affine change of $q$ and real multiplicative translations, so that $z_1=q$ and $z_2=q-1$. Since $p_3$ is not on the line of the triple, $(q,s)$ with $s=z_3$ are affine coordinates on $V$, and $z_4=\alpha s+\beta q+\gamma$ with real coefficients; $\alpha\ne0$ because $p_4$ is not on the line of the triple, and $\beta\ne0$ because $p_0,p_3,p_4$ are not collinear. A real multiplicative translation of $z_4$ gives $z_4=s-\mu(q-\varpi')$ with $\mu\ne0$, and $\varpi'\notin\{0,1\}$ because for $\varpi'=0$, respectively $1$, the points $p_1,p_3,p_4$, respectively $p_2,p_3,p_4$, would form a second collinear triple. Case (a) is Corollary \ref{cor:square4}. The three cases are distinguished by the number of collinear triples, an invariant of monomial equivalence.

The amoeba map is four-to-one in case (a) by Theorem \ref{thm:square}, in case (b) by Theorem \ref{thm:fibred} and in case (c) by Proposition \ref{prop:prodfib}, and these properties are invariant under monomial equivalence by Lemma \ref{lem:monomial}. For the volume, the linear maps $L_a$ of Lemma \ref{lem:monomial} have determinant $\pm1$, so they preserve the Lebesgue measure of $\R^4$, and $\Vol(\A(V^{(a)}))=\Vol(\A(V))$; it therefore suffices to compute the volume for the normal forms. For (a) it is $\pi^4/4$ by Corollary \ref{cor:square4}, and for (c) it is $(\pi^2/2)^2=\pi^4/4$, the area of the amoeba of a real line of $(\C^*)^2$ being $\pi^2/2$ \cite{PR04}. For (b), \eqref{eq:rolled} shows that $\MVol(V_\mu)=\pi^4$ as soon as almost every fibre of the rolled coamoeba map has exactly one point, and indeed the conditions $\operatorname{Im}(e^{-i\varphi_1}q)=\operatorname{Im}(e^{-i\varphi_2}(q-1))=0$ determine $q$ uniquely when $\varphi_1\ne\varphi_2$, as the intersection of two distinct real lines through $0$ and $1$, and then $\operatorname{Im}(e^{-i\varphi_3}s)=\operatorname{Im}(e^{-i\varphi_4}(s-\mu(q-\varpi')))=0$ determine $s$ uniquely when $\varphi_3\ne\varphi_4$ and $\mu(q-\varpi')\ne0$; hence $\Vol(\A(V_\mu))=\frac14\MVol(V_\mu)=\pi^4/4$. Finally $\MVol(V)=4\Vol(\A(V))=\pi^4$ for every $V$ in the theorem, by the identity $\Vol(\A(V))=\int_VJ(\Log|_V)N(\Log z)^{-1}\,d\mathrm{vol}$ with $N=4$.
\end{proof}

Thus, for $k=2$, every real plane of $(\C^*)^4$ with a nondegenerate amoeba has a four-to-one amoeba map, and the classes of the question at the beginning of this section are exactly the squares, the products of two lines and the fibred extensions $V_\mu$, up to monomial equivalence. The fibred extensions $V_\mu$ are not products, and are not monomially equivalent to products, since they have one collinear triple instead of two. The statement $\Vol(\A(V))=\pi^4/4$ extends Corollary \ref{cor:sizes} and the value of \cite[Theorem 1.2]{NP17} from generic real planes to all real planes of $(\C^*)^4$ with nondegenerate amoebas.

For $k=3$ we use the following reduction, valid for any $k$, when three dual points are collinear. Suppose, after a monomial equivalence, that $p_0,p_1,p_2$ are collinear, so that $z_1=q-\varpi_1$ and $z_2=\nu(q-\varpi_2)$ for a real affine form $q$ on $V$, real $\nu\ne0$ and distinct real $\varpi_1,\varpi_2$. For $q\in\C$ let $P_q\subset\C^{2k-2}$ be the set of the coordinates $(z_3,\dots,z_{2k})$ of the points of $V$ with the given value of $q$, a complex affine $(k-1)$-plane with real direction space, and with $P_{\bar q}=\overline{P_q}$. For $x\in\R^{2k}$ write $x=(x_1,x_2,x'')$.

\begin{lemma}\label{lem:fibline}
Let $V$ be as above, and let $x$ be a point such that $(x_1,x_2)$ lies in the amoeba of the real line $\{(q-\varpi_1,\nu(q-\varpi_2))\}$ of $(\C^*)^2$ outside the image of its real points, the two points of its fibre being $q$ and $\bar q$. Then $N_V(x)=2N_{P_q}(x'')$.
\end{lemma}

\begin{proof}
A point of $V$ over $x$ has its value of $q$ in the fibre $\{q,\bar q\}$ of the amoeba map of the line, and for each of the two values its last coordinates form a point of $P_q$, respectively $P_{\bar q}$, over $x''$. Conjugation maps the points of $P_q$ over $x''$ bijectively onto those of $P_{\bar q}$ over $x''$.
\end{proof}

Hence $V$ is $2^k$-to-one if and only if the complex planes $P_q$ are $2^{k-1}$-to-one for almost every $q$. For $k=3$ the planes $P_q$ are complex planes of $\C^4$, which can be written explicitly. If $p_0,p_3,p_4,p_5,p_6$ are in general position, then $B_3,\dots,B_6$ form a basis of $\R^4$, and $e_0=\sum_jc_jB_j$ with all $c_j\ne0$, since $p_0$ lies in no plane spanned by three of the $p_j$, $j\ge3$. As functions on $V$ this reads $\sum_jc_jz_j=1$, and after the real multiplicative translation $z_j\mapsto c_jz_j$ the projection of $V$ to $(z_3,\dots,z_6)$ is contained in the real hyperplane $U=\{z_3+z_4+z_5+z_6=1\}$; since $B_3,\dots,B_6$ are independent this projection is all of $U$, and the form $q$ is the restriction to $U$ of a real linear form $\lambda(z)=\sum_jl_jz_j$, with $l\in\R^4$ determined up to adding multiples of $(1,1,1,1)$, which do not change $\lambda$ on $U$. Thus

$$V=\big\{(q-\varpi_1,\ \nu(q-\varpi_2),\ z):\ z\in U,\ \lambda(z)=q\big\},\qquad P_q=\{z\in\C^4:\ z_3+z_4+z_5+z_6=1,\ l\cdot z=q\},$$

up to real multiplicative translations, where $l$ is not a multiple of $(1,1,1,1)$, because $q$ is not constant on $V$.

\begin{theorem}\label{thm:triple}
Let $V\subset\C^6$ be a real affine $3$-plane with $\Hm^6(\A(V))>0$, whose dual configuration contains three collinear points $p_a,p_b,p_c$ such that, for one of them, say $p_c$, the five points $p_c$ and $p_j$, $j\notin\{a,b,c\}$, are in general position. Then the following are equivalent:

(i) $N_V\ge8$ almost everywhere on $\A(V)$;

(ii) the amoeba map of $V$ is exactly eight-to-one almost everywhere;

(iii) the line through $p_a,p_b,p_c$ meets two opposite edges of the tetrahedron whose vertices are the four points $p_j$, $j\notin\{a,b,c\}$;

(iv) $V$ is monomially equivalent to a fibred extension $\{(q-\varpi_1,\,q-\varpi_2,\,s_1,\,s_1-\kappa_1(q),\,s_2,\,s_2-\kappa_2(q))\}$ of a real line by two lines with nonconstant $\kappa_1,\kappa_2$.

If these conditions fail, the amoeba map of $V$ has fibres of exactly $4$ points on a set of positive measure, near the tentacles of $\A(V)$ in the directions of the coordinate hyperplanes corresponding to the four points $p_j$, $j\notin\{a,b,c\}$, after the monomial equivalence sending $p_c$ to $p_0$.
\end{theorem}

\begin{proof}
By Lemma \ref{lem:monomial} we may assume $p_c=p_0$ and $\{a,b\}=\{1,2\}$, and we use the normal form above, with $l=(l_3,l_4,l_5,l_6)$.

Assume (iii). In the normalised coordinates the vectors $B_3,\dots,B_6$ satisfy $\sum_jB_j=e_0$, and the identity $\lambda=q$ on $V$ says that $\sum_jl_jB_j=Q$, the vector of coefficients of the form $q$; the line of the triple is spanned by $e_0$ and $Q$. If it meets the edges $\overline{p_ip_j}$ and $\overline{p_rp_s}$, with $\{i,j,r,s\}=\{3,4,5,6\}$, in the points $X=[x_iB_i+x_jB_j]$ and $Y=[y_rB_r+y_sB_s]$, then $X\ne Y$ since the edges are skew, the line is spanned by $X$ and $Y$, and writing $e_0=\sum B_j=\sigma X+\tau Y$ in the basis $B_3,\dots,B_6$ gives $x_i=x_j$ and $y_r=y_s$; then $Q\in\operatorname{span}(B_i+B_j,B_r+B_s)$ means that $l_i=l_j=\alpha$ and $l_r=l_s=\beta$, with $\alpha\ne\beta$ since $l$ is not a multiple of $(1,1,1,1)$. Conversely these equalities imply (iii) by the same computation. Under them, on $U$ we have $\lambda=\alpha(z_i+z_j)+\beta(1-z_i-z_j)$, so $z_i+z_j=\kappa_1(q)$ and $z_r+z_s=1-\kappa_1(q)=:\kappa_2(q)$ with $\kappa_1(q)=(q-\beta)/(\alpha-\beta)$; putting $s_1=z_i$ and $s_2=z_r$ and multiplying $z_j$ and $z_s$ by $-1$, $V$ becomes the fibred extension in (iv), with $\kappa_1,\kappa_2$ nonconstant. This proves (iii)$\Rightarrow$(iv), and (iv)$\Rightarrow$(ii) is Theorem \ref{thm:fibred} together with Lemma \ref{lem:monomial}, while (ii)$\Rightarrow$(i) is trivial.

Assume (i), and fix $m\in\{3,4,5,6\}$; let $\{a',b',c'\}$ be the other three indices. For nonreal $q$, the section $L_q=P_q\cap\{z_m=0\}$ is the complex line $\{z_{a'}+z_{b'}+z_{c'}=1,\ l_{a'}z_{a'}+l_{b'}z_{b'}+l_{c'}z_{c'}=q\}$ of $\C^3$. If $l_{a'},l_{b'},l_{c'}$ are pairwise distinct, $L_q$ meets the three coordinate planes and the plane at infinity in four distinct points, and in the affine coordinate $z_{b'}$ on $L_q$ these points are $0$, $\beta_{a'}=(q-l_{c'})/(l_{b'}-l_{c'})$, $\beta_{c'}=(q-l_{a'})/(l_{b'}-l_{a'})$ and $\infty$. Their cross-ratio is real if and only if $\beta_{a'}/\beta_{c'}\in\R$, that is, if and only if $(q-l_{c'})/(q-l_{a'})\in\R$, which fails for nonreal $q$ since $l_{a'}\ne l_{c'}$. By the injectivity criterion of Section \ref{sec:injlines}, applied to $L_q$ after an affine change of parameter and a multiplicative translation which bring it to the normal form $L_c$ with $c\notin\R$, the amoeba map of $L_q$ is injective. We now apply Lemma \ref{lem:tentacle} to the complex plane $P_q\subset\C^4$ and its section $L_q$; the proof of that lemma uses neither the reality nor the general position of the plane, only that the points of the section over $y$ are points where the amoeba map of the section is an immersion with $\chi\ne0$, where $\chi$ is defined as before with a vector $u$ of the direction space of $P_q$ with $u_m=1$, which exists since $z_m$ is not constant on $P_q$. The amoeba map of $L_q$ is an immersion outside a proper analytic subset, since its amoeba is $2$-dimensional. As in the proof of Lemma \ref{lem:chi}, $\chi(h)=0$ implies that $h$ is a critical point of the amoeba map of the projection $P'_q$ of $P_q$ forgetting $z_m$, which is the plane $\{\sum_{j\ne m}b_jz_j=e\}$ of $\C^3$ with $b_j=l_j-l_m$ and $e=q-l_m\notin\R$. A point $z$ of this plane is critical if and only if there is a real vector $\eta\ne0$ with $(\eta_j/z_j)_j$ proportional to $(b_j)_j$, that is, if and only if the numbers $b_jz_j$ with $b_j\ne0$ lie on a common real line through $0$; since they sum to $e$, this line is $\R e$. If every point of $L_q$ were critical, the holomorphic functions $b_jz_j/e$ would be real on $L_q$, hence constant, for every $j$ with $b_j\ne0$; at most one of $l_{a'},l_{b'},l_{c'}$ equals $l_m$, so at least two coordinates would be constant on $L_q$, and then the third too, because their sum is $1$, which is absurd. Hence $\chi\ne0$ outside a proper analytic subset of $L_q$, and Lemma \ref{lem:tentacle} with $M=1$ gives, for every nonreal $q$, a set of positive measure in $\R^4$ on which $N_{P_q}=2$. We show that this produces an open subset of $\R^6$ on which $N_V=4$. Fix $q_0$ with $\operatorname{Im}q_0>0$. Since the critical values of the amoeba map of $P_{q_0}\cap(\C^*)^4$ form a null set by Sard's theorem, there is a regular value $x''_0$ of this map with $N_{P_{q_0}}(x''_0)=2$; let $z^{(1)},z^{(2)}$ be its two preimages. Consider the real analytic manifold $\mathcal P=\{(q,z):\operatorname{Im}q>0,\ z\in P_q\cap(\C^*)^4\}$, which is the image of $\{\operatorname{Im}q>0\}\times\C^2$ under an affine parametrisation of the planes $P_q$ depending affinely on $q$, intersected with the open set $\{z\in(\C^*)^4\}$, and the map $\Psi:\mathcal P\to\{\operatorname{Im}q>0\}\times\R^4$, $\Psi(q,z)=(q,\Log z)$. The map $\Psi$ is proper: if $(q_n,z_n)\in\mathcal P$ and $\Psi(q_n,z_n)$ converges to $(q_*,x''_*)$ with $\operatorname{Im}q_*>0$, then the moduli of the coordinates of $z_n$ are bounded and bounded away from $0$, so a subsequence of $z_n$ converges to a point $z_*\in(\C^*)^4$, which lies in $P_{q_*}$ because the conditions $\sum_jz_j=1$, $l\cdot z=q$ are closed. At $(q_0,z^{(i)})$ the differential of $\Psi$ is invertible, because its restriction to the fibre $\{q_0\}\times P_{q_0}$ is the differential of the amoeba map of $P_{q_0}$ at a regular point; by the inverse function theorem there are disjoint neighbourhoods $U_i$ of $(q_0,z^{(i)})$ on which $\Psi$ is a diffeomorphism onto a neighbourhood of $(q_0,x''_0)$. By properness, there is a neighbourhood $O$ of $(q_0,x''_0)$ such that $\Psi^{-1}(O)\subset U_1\cup U_2$: otherwise there would be points $(q_n,z_n)\notin U_1\cup U_2$ with $\Psi(q_n,z_n)\to(q_0,x''_0)$, and a limit point of $(q_n,z_n)$ would be a third preimage of $(q_0,x''_0)$, outside $U_1\cup U_2$. Shrinking $O$, we obtain $N_{P_q}(x'')=2$ for every $(q,x'')\in O$. Now the map $q\mapsto(\log|q-\varpi_1|,\log|\nu(q-\varpi_2)|)$ is a local diffeomorphism on $\{\operatorname{Im}q>0\}$, its Jacobian being $\operatorname{Im}((q-\varpi_1)^{-1}\overline{(q-\varpi_2)^{-1}})\ne0$ there, so the image of $O$ under $(q,x'')\mapsto(\log|q-\varpi_1|,\log|\nu(q-\varpi_2)|,x'')$ contains a nonempty open set of $\R^6$. Every point $x$ of this set satisfies the hypothesis of Lemma \ref{lem:fibline}, the fibre of the amoeba map of the real line $\{(q-\varpi_1,\nu(q-\varpi_2))\}$ over $(x_1,x_2)$ being $\{q,\bar q\}$ with $\operatorname{Im}q>0$, and Lemma \ref{lem:fibline} gives $N_V(x)=2N_{P_q}(x'')=4$. This contradicts (i). Therefore, for every $m\in\{3,4,5,6\}$, two of the three numbers $l_j$, $j\ne m$, are equal.

A vector $(l_3,l_4,l_5,l_6)$ such that every three of its entries contain two equal ones has either two pairs of equal entries, or three equal entries, or all four equal. The last case is excluded. If $l_i=l_j=l_r\ne l_s$, then $\lambda=l_i(1-z_s)+l_sz_s$ on $U$, so $z_s$ is an affine function of $q$ and $p_s$ lies on the line of the triple; then $V$ consists of the points whose coordinates $z_1,z_2,z_s$ lie on a real line of $\C^3$, of $2$-dimensional amoeba, and whose other three coordinates lie on the plane $\{z_i+z_j+z_r=1-z_s\}$ of $\C^3$, whose amoeba is $3$-dimensional, so that $\A(V)$ has dimension at most $5$, a contradiction. Hence $l$ has two pairs of equal entries with different values, which is (iii). The last assertion of the theorem was proved along the way.
\end{proof}

Condition (iii) has codimension one among the configurations of the theorem: the lines through $p_0$ meeting the line $\overline{p_3p_4}$ form a plane, and those meeting also $\overline{p_5p_6}$ form a single line, while the line of the triple can be any line through $p_0$. So, as for symmetric powers in general position, the $8$-to-one property is exceptional in the stratum of configurations with a collinear triple, and it holds exactly on the fibred extensions.

The other strata of degenerate configurations of seven points in $\mathbb P^3$ contain the remaining examples known to us. The products $\ell\times W$ of a real line of $(\C^*)^2$ and a real plane $W$ of $(\C^*)^4$ with a nondegenerate amoeba are eight-to-one, by Proposition \ref{prop:prodfib} and Theorem \ref{thm:planes4}; their dual points lie on a line and a plane through $p_0$ meeting only at $p_0$. The fibred extensions of a square $\ell\odot\ell\subset\C^4$ by one line with $\kappa=\mu f_q(\varpi')$ are eight-to-one by Theorem \ref{thm:fibred}; their dual points are five points $p_0,\dots,p_4$ on a conic $C$ in a plane $\Pi$, a point $p_5$ outside $\Pi$, and a point $p_6$ on the line joining $p_5$ to the point of $C$ with parameter $\varpi'$, because $B_6=B_5-\mu(\varpi'^2,\varpi',1,0)$. By Theorem \ref{thm:fibred}, the fibred extensions of the line $\{(q,q-1)\}$ by two lines, with constant or nonconstant $\kappa_j$, the products and the symmetric cubes, these exhaust the constructions of the previous subsections for $k=3$. We conjecture that there are no others.

\begin{conjecture}\label{conj:k3}
A real affine $3$-plane $V\subset\C^6$ with $\Hm^6(\A(V))>0$ has an eight-to-one amoeba map if and only if it is monomially equivalent to a symmetric cube of a real line, to a product $\ell\times W$ of a real line of $(\C^*)^2$ and a real plane of $(\C^*)^4$ with nondegenerate amoeba, to a fibred extension of a real line by two lines, or to a fibred extension of the square of a real line by a line with $\kappa=\mu f_q(\varpi')$.
\end{conjecture}

Theorems \ref{thm:classif} and \ref{thm:triple} prove the conjecture for configurations in general position and for configurations with a collinear triple whose complement is in general position with one point of the triple. The method extends to the other strata: in each of them, Lemmas \ref{lem:monomial} and \ref{lem:fibline}, or a fibration over the projection to the coordinates of a coplanar quadruple through $p_0$, reduce the eight-to-one property to properties of the fibres, which can then be tested near the tentacles by Lemma \ref{lem:tentacle}, provided its nondegeneracy hypothesis can be checked; the number of strata is large, and we have not carried out this case analysis.

\section{\texorpdfstring{Lines in $(\C^*)^3$}{Lines in (C*)\textasciicircum 3}}\label{sec:lines}

\subsection{Normal form}
Let $L\subset(\C^*)^3$ be the intersection with the torus of an affine line of $\C^3$ which meets the three coordinate planes in three distinct points, none at infinity, and which is not contained in a translate of a two-dimensional subtorus. Parametrise it as $z_j=a_j+b_jt$ with $b_1\ne0$; after the change of parameter $t\mapsto t-a_1/b_1$ and a rescaling of $z_1$ we may assume $z_1=s$; then $z_2=a_2+b_2s$ with $a_2b_2\ne0$, and after rescaling $z_2$ and the parameter, and then $z_1$ again, we obtain $z_1=s$, $z_2=1-s$; rescaling $z_3$ gives $z_3=1-cs$. Multiplicative translations translate the amoeba, so we may assume

$$L=L_c=\{(s,\,1-s,\,1-cs):\ s\in\C\setminus\{0,1,1/c\}\},\qquad c\in\C\setminus\{0,1\}.$$

The parameter $c$ is the cross-ratio of the four points $0,1,1/c,\infty$ where $L_c$ meets the toric boundary. For $c=0$ or $c=1$ the line lies in a translate of a two-dimensional subtorus ($z_3=1$ or $z_3=z_2$); we call these lines planar. The line $L_c$ is the smooth complete intersection of $f_1=1-z_1-z_2$ and $f_2=1-cz_1-z_3$, with Newton polytopes $\operatorname{conv}(0,e_1,e_2)$ and $\operatorname{conv}(0,e_1,e_3)$, and it is real up to translation if and only if $c\in\R$.

\subsection{The three projections}
The coordinate projections of $L_c$ are, up to finitely many points, the plane lines $\{1-z_1-z_2=0\}$, $\{1-cz_1-z_3=0\}$ and $\{(1-c)+cz_2-z_3=0\}$, all with three nonzero coefficients since $c\ne0,1$, and $\pi_I|_{L_c}$ is injective, being the restriction of a linear map to a line not parallel to its kernel. For $I=\{1,2\}$ the system \eqref{eq:doubled} reads $1-e^{i\varphi_1}r_1-e^{i\varphi_2}r_2=0$, $1-e^{-i\varphi_1}r_1-e^{-i\varphi_2}r_2=0$, $1-ce^{i\varphi_1}r_1-z_3=0$, $1-\bar ce^{-i\varphi_1}r_1-w_3=0$; the first two equations are linear in $(r_1,r_2)$ and the last two determine $z_3,w_3$, so for generic coefficients there is exactly one solution and $\alpha_{12}=1$. The same elimination gives $\alpha_{13}=\alpha_{23}=1$, so Theorem \ref{thm:main} yields $\Hm^2(\A(L_c))\le 3\pi^2$. Each projected multivolume can be computed exactly: since $\pi_I|_{L_c}$ is injective, $\MVol_I(L_c)$ is the multivolume of the amoeba of a generic line in $(\C^*)^2$, which is twice its area $\pi^2/2$ \cite{PR04,Pas08} (see also \cite[Proposition 8]{Mik17}). Hence $\MVol_{12}(L_c)=\MVol_{13}(L_c)=\MVol_{23}(L_c)=\pi^2$, and Proposition \ref{prop:fibre} is an equality for each projection.

\subsection{An explicit formula}
Write $v(s)=\frac{d}{ds}\log z(s)$, where $z(s)=(s,1-s,1-cs)$. Then $v_i=1/(s-q_i)$ with $q=(0,1,1/c)$, because $-\frac1{1-s}=\frac1{s-1}$ and $-\frac{c}{1-cs}=\frac1{s-1/c}$. Writing $v=a+ib$ with $a,b\in\R^3$, the differential of $\Log|_{L_c}$ sends the coordinate vectors of $s=u_1+iu_2$ to $a$ and $-b$, so $J(\Log|_{L_c})\,d\mathrm{vol}=|a\wedge b|\,dA(s)$, where $dA$ is the Lebesgue measure of the parameter plane. Put $w_{ij}=\operatorname{Im}(\bar v_iv_j)=a_ib_j-a_jb_i$ and $\rho_i=|s-q_i|$.

\begin{lemma}\label{lem:formula}
For $s\notin\{q_1,q_2,q_3\}$ one has $|a\wedge b|=(\sum_{i<j}w_{ij}^2)^{1/2}=\frac12\big((\sum_i|v_i|^2)^2-|\sum_iv_i^2|^2\big)^{1/2}$, and $|w_{ij}(s)|=|\sin\theta_{ij}(s)|/(\rho_i\rho_j)$, where $\theta_{ij}(s)$ is the angle at $s$ of the triangle with vertices $q_i,s,q_j$. Moreover $|a\wedge b|>0$ unless $q_1,q_2,q_3$ and $s$ are collinear.
\end{lemma}

\begin{proof}
The first identity is the Lagrange identity $|a\wedge b|^2=|a|^2|b|^2-(a\cdot b)^2$. For the second, $\sum|v_i|^2=|a|^2+|b|^2$ and $\sum v_i^2=|a|^2-|b|^2+2i\,a\cdot b$, so the difference of squares is $4|a|^2|b|^2-4(a\cdot b)^2$. For the third, $\bar v_iv_j=(s-q_i)\overline{(s-q_j)}/(\rho_i^2\rho_j^2)$ and $\operatorname{Im}((s-q_i)\overline{(s-q_j)})=\rho_i\rho_j\sin(\arg(s-q_i)-\arg(s-q_j))$. Finally $|a\wedge b|=0$ forces $\sin\theta_{ij}=0$ for every pair.
\end{proof}

By Lemma \ref{lem:formula} and the identity \eqref{eq:passare} of Section \ref{sec:notation}, $\int_\C|w_{ij}|\,dA=\pi^2$ for every pair $i<j$ and every $c$, which is the statement $\MVol_I(L_c)=\pi^2$.

\subsection{Injectivity of the amoeba map}\label{sec:injlines}
Let $s\neq s'$ satisfy $\Log z(s)=\Log z(s')$. The conditions $|s|=|s'|$ and $|1-s|=|1-s'|$ give $s'=\bar s$ with $\operatorname{Im}s\neq0$, and since $|1-cs|^2-|1-c\bar s|^2=-2\operatorname{Re}(c(s-\bar s))=4\operatorname{Im}(c)\operatorname{Im}(s)$, the third condition forces $\operatorname{Im}c=0$. Thus $\Log|_{L_c}$ is injective if $c\notin\R$, and exactly two-to-one off the real points if $c\in\R$, and by the area formula

$$\Hm^2(\A(L_c))=\MVol(L_c)\ \text{ if }c\notin\R,\qquad \Hm^2(\A(L_c))=\tfrac12\MVol(L_c)\ \text{ if }c\in\R.$$

The real case is the equality case of Corollary \ref{cor:half}. The same argument shows that for $c\notin\R$ a generic point of $\A(\pi_{ij}(L_c))$ has two preimages in $\A(L_c)$, the images of $s$ and $\bar s$, separated by the third coordinate, so $\int_{\A(L_c)}J^{\A(L_c)}_{ij}\,d\Hm^2=\pi^2$, while for $c\in\R$ this integral equals $\pi^2/2$.

\begin{theorem}\label{thm:line}
If $c\notin\R\cup\{0,1\}$, then $\sqrt3\,\pi^2<\Hm^2(\A(L_c))<3\pi^2$. If $c\in\R\setminus\{0,1\}$, then $\frac{\sqrt3}{2}\pi^2<\Hm^2(\A(L_c))<\frac32\pi^2$. In particular, for nonreal $c$ the bound $\frac{\pi^2}2\sum_I\alpha_I=\frac32\pi^2$ fails.
\end{theorem}

\begin{proof}
For $c\notin\R$ put $S=\A(L_c)$. By \eqref{eq:cb}, $1\le\sum_{i<j}J^S_{ij}\le\sqrt3$, and $\int_SJ^S_{ij}\,d\Hm^2=\pi^2$, so $\sqrt3\,\pi^2\le\Hm^2(S)\le3\pi^2$. Equality on the right would force the tangent planes of $S$ to be coordinate planes almost everywhere, and equality on the left would force $J^S_{ij}=1/\sqrt3$ almost everywhere; neither happens, because the tangent plane at $\Log z(s)$ is spanned by $a(s),b(s)$ and varies real analytically and non-trivially (near $s=0$ it tends to a plane containing $e_1$, so $J^S_{23}\to0$). The real case is identical with all integrals halved.
\end{proof}

\section{The multivolume of lines as a function of the cross-ratio}\label{sec:cross}

We now study $\Mv(c)=\MVol(L_c)=\int_\C(\sum_{i<j}w_{ij}^2)^{1/2}\,dA$, which by Section \ref{sec:injlines} equals $\Hm^2(\A(L_c))$ for $c\notin\R$ and $2\Hm^2(\A(L_c))$ for $c\in\R$, and which satisfies $\sqrt3\,\pi^2<\Mv(c)<3\pi^2$. This is the case $k=1$, $n=3$ of the question of the extremal volumes of amoebas of linear spaces, which for half-dimensional linear spaces was posed by Mikhalkin \cite[Problem 3.1]{Mik17}.

\subsection{Symmetries and the space of triangles}

\begin{proposition}\label{prop:sym12}
$\Mv(c)$ depends only on the similarity class, reflections allowed and vertices unlabelled, of the triangle with vertices $0,1,c$. Equivalently, $\Mv$ is invariant under the group $\Sg$ of order $12$ generated by $c\mapsto1-c$, $c\mapsto1/c$ and $c\mapsto\bar c$.
\end{proposition}

\begin{proof}
By Lemma \ref{lem:formula}, $\Mv(c)$ is a symmetric function of $q_1,q_2,q_3$. Under a similarity $s=\lambda\sigma+\mu$ the points become $(q_i-\mu)/\lambda$, $v$ becomes $\lambda^{-1}v'$, the $w_{ij}$ are multiplied by $|\lambda|^{-2}$ and $dA$ by $|\lambda|^2$; under $s\mapsto\bar s$, $v$ becomes $\bar v$ and $|w_{ij}|$ is unchanged. Hence $\Mv$ depends only on the similarity class of $\{0,1,1/c\}$, which after multiplication by $c$ is the class of $\{0,1,c\}$. The transformations of $c$ preserving this class are generated by the relabellings $c\mapsto1-c$, $c\mapsto1/c$, which generate the anharmonic group of order $6$ \cite{Ahl79}, together with $c\mapsto\bar c$.
\end{proof}

Although $c$ is the cross-ratio of $0,1,1/c,\infty$, $\Mv$ is not invariant under all permutations of these four points: the tentacle of the amoeba corresponding to $\infty$ has direction $(1,1,1)$ while the other three have directions $-e_1,-e_2,-e_3$, and a permutation exchanging $(1,1,1)$ with some $-e_i$ is induced by a lattice automorphism of $\Z^3$ which is not an isometry. A fundamental domain for $\Sg$ is $D=\{c:\operatorname{Im}c\ge0,\ \operatorname{Re}c\le\frac12,\ |c|\ge1\}$. Its boundary consists of the half-line $(-\infty,-1]$ of collinear triangles, ending at the harmonic configuration $c=-1$ in which one vertex is the midpoint of the other two; the arc $\{e^{i\alpha}:\frac\pi3\le\alpha\le\pi\}$ of isosceles triangles with apex $0$ and apex angle $\alpha\ge60^\circ$; and the half-line $\{\frac12+iy:y\ge\frac{\sqrt3}2\}$ of isosceles triangles with apex angle smaller than $60^\circ$. The circles and lines $|c|=1$, $|c-1|=1$, $\operatorname{Re}c=\frac12$ are the fixed point sets of the reflections of $\Sg$, the orbit $\{-1,\frac12,2\}$ consists of the harmonic configurations, the orbit $\{e^{\pm i\pi/3}\}$ of the equilateral ones, and $c=0,1,\infty$ correspond to the collision of two vertices.

\subsection{Regularity}

\begin{proposition}\label{prop:reg}
$\Mv$ is continuous on $\C\setminus\{0,1\}$, smooth on $\R\setminus\{0,1\}$, and real analytic on $\C\setminus\R$.
\end{proposition}

\begin{proof}
Let $c\to c_0$. Off $\{0,1,1/c_0\}$ the integrand converges pointwise and is dominated by $g_c=|w_{12}|+|w_{13}|+|w_{23}|$, which converges pointwise to $g_{c_0}$ and has integral $3\pi^2$ for every $c$ by \eqref{eq:passare}. Pratt's lemma \cite{Pra60}, the version of dominated convergence with variable dominating functions whose integrals converge, gives continuity.

For smoothness we use a partition of unity on the Riemann sphere subordinate to small discs around $q_1,q_2,q_3,\infty$ and a set avoiding them. Away from these points the integrand is smooth in $(s,c)$, and real analytic for noncollinear points, being the square root of a positive real analytic function. Near $q_k$ we use moving coordinates $s=q_k+re^{i\varphi}$; then $\bar v_k=e^{i\varphi}/r$ and $w_{kj}=r^{-1}h_{kj}$ with $h_{kj}=\operatorname{Im}(e^{i\varphi}/(q_k-q_j+re^{i\varphi}))$ real analytic up to $r=0$, while the remaining $w_{jl}$ is bounded, so the integrand times $dA=r\,dr\,d\varphi$ equals $(h_{kj}^2+h_{kl}^2+r^2w_{jl}^2)^{1/2}dr\,d\varphi$. At $r=0$ the expression under the root vanishes only if $e^{i\varphi}$ is parallel to both $q_k-q_j$ and $q_k-q_l$, which is impossible for noncollinear points. Near $\infty$, with $s=1/\sigma$, $\sigma=re^{i\varphi}$, we get $w_{ij}=r^3(\operatorname{Im}((q_j-q_i)e^{i\varphi})+O(r))$ up to sign and $dA=r^{-3}dr\,d\varphi$, and the integrand becomes the square root of a real analytic function equal at $r=0$ to $\sum_{i<j}\operatorname{Im}((q_j-q_i)e^{i\varphi})^2>0$. Differentiation under the integral sign is therefore legitimate and $\Mv$ is real analytic off $\R$. For real $c$, Lemma \ref{lem:formula} gives $w_{ij}=y(q_i-q_j)/(\rho_i^2\rho_j^2)$ for $s=x+iy$, so the integrand is $|y|$ times the square root of a positive smooth function; in the polar charts $|y|=r|\sin\varphi|$ does not depend on $c$, and the same argument gives smoothness in the real variable $c$.
\end{proof}

\subsection{Collision of two vertices}

\begin{theorem}\label{thm:degen}
As $c\to0$, $c\to1$ or $c\to\infty$, $\Mv(c)\to(1+\sqrt2)\pi^2$. Hence $\Hm^2(\A(L_c))\to(1+\sqrt2)\pi^2$ along nonreal values and $\Hm^2(\A(L_c))\to\frac{1+\sqrt2}2\pi^2$ along real values.
\end{theorem}

\begin{proof}
By Proposition \ref{prop:sym12} it suffices to treat $c\to0$ with $q=(0,1,c)$. Let $w^0_{12}$ be the function $w_{12}$ for the points $0,1$. As $c\to0$, $v_3\to v_1$ pointwise, so $w_{23}\to-w^0_{12}$ and $w_{13}\to0$. For the upper bound, $(\sum w_{ij}^2)^{1/2}\le|w_{13}|+Q_c$ with $Q_c=(w_{12}^2+w_{23}^2)^{1/2}$; here $\int|w_{13}|=\pi^2$ by \eqref{eq:passare}, and $Q_c\to\sqrt2|w^0_{12}|$ pointwise with $Q_c\le|w_{12}|+|w_{23}|$, whose integral is $2\pi^2$ for every $c$ and whose pointwise limit $2|w^0_{12}|$ has the same integral; Pratt's lemma gives $\int Q_c\to\sqrt2\pi^2$, so $\limsup\Mv(c)\le(1+\sqrt2)\pi^2$. For the lower bound fix $\delta\in(0,\frac12)$ and let $|c|<\delta/2$; use $(\sum w_{ij}^2)^{1/2}\ge|w_{13}|$ on $B_\delta=\{|s|<\delta\}$ and $\ge Q_c$ outside. For $|s|\ge\delta$, the law of sines in the triangle $0,c,s$ gives $|\sin\theta_{13}|\le|c|/|s|$, and $|s-c|\ge|s|/2$, so $|w_{13}|\le2|c|/|s|^3$ and $\int_{\C\setminus B_\delta}|w_{13}|\,dA\le4\pi|c|/\delta\to0$; hence $\int_{B_\delta}|w_{13}|\to\pi^2$. Outside $B_\delta$ we have $|w_{12}|=|w^0_{12}|$ and, since $|s-c|\ge|s|/2$ and, by the law of sines in the triangle $1,c,s$, $|\sin\theta_{23}|\le|1-c|/|s-1|\le2/|s-1|$, also $|w_{23}|\le G(s)=2|s|^{-1}|s-1|^{-1}\min(1,2|s-1|^{-1})$; the function $G$ is integrable on $\C\setminus B_\delta$, being $O(|s-1|^{-1})$ near $1$ and $O(|s|^{-3})$ at infinity. Hence $Q_c\le|w^0_{12}|+G$ on $\C\setminus B_\delta$ for every $|c|<\delta/2$, and the dominated convergence theorem gives $\int_{\C\setminus B_\delta}Q_c\to\sqrt2\int_{\C\setminus B_\delta}|w^0_{12}|$, and $|w^0_{12}|\le2/|s|$ on $B_\delta$ gives $\int_{B_\delta}|w^0_{12}|\le4\pi\delta$. Thus $\liminf\Mv(c)\ge\pi^2+\sqrt2(\pi^2-4\pi\delta)$ for every $\delta$.
\end{proof}

When two vertices collide, the mass $\pi^2$ of the colliding pair concentrates at the scale $|c|$, where the amoeba looks like the amoeba of a line in a plane, while at scale one the two colliding coordinates become equal and the two remaining pairs contribute the length $\sqrt2$ of $(w^0_{12},-w^0_{12})$ instead of $2$. The limit $(1+\sqrt2)\pi^2$ is therefore strictly larger than the multivolume $\pi^2$ of the planar line obtained in the limit.

\begin{proposition}\label{prop:rate}
For $c=\varepsilon e^{i\tau}$ with $0<\varepsilon\le\frac18$ one has
$$\Mv(c)=(1+\sqrt2)\pi^2-\kappa(\tau)\sqrt\varepsilon+R(c),\qquad |R(c)|\le C\varepsilon^{2/3},$$
with an absolute constant $C$, where
$$\kappa(\tau)=\frac{2^{1/4}\,4\pi^{3/2}}{\Gamma(1/4)^2}\int_0^{2\pi}\sqrt{|\sin\varphi\,\sin(\varphi-\tau)|}\,d\varphi .$$
The function $\kappa$ is continuous and positive, and maximal exactly for $\tau\in\{0,\pi\}$, where $\kappa(0)=2^{1/4}\cdot16\pi^{3/2}/\Gamma(1/4)^2$. In particular there is $\varepsilon_0>0$ such that $\Mv(c)<(1+\sqrt2)\pi^2$ whenever $0<|c|<\varepsilon_0$, $0<|c-1|<\varepsilon_0$ or $|c|>1/\varepsilon_0$.
\end{proposition}

\begin{proof}
As in the proof of Theorem \ref{thm:degen} we use the points $q=(0,1,c)$ and write $s=\rho e^{i\varphi}$, $W=|w_{13}|$, $Q=(w_{12}^2+w_{23}^2)^{1/2}$ and $F=(W^2+Q^2)^{1/2}$, so that $\Mv(c)=\int_\C F\,dA$. Put $h(W,Q)=W+Q-(W^2+Q^2)^{1/2}$. This function is positively homogeneous of degree one, satisfies $0\le h(W,Q)\le\min(W,Q)$, and its partial derivatives $1-W/F$ and $1-Q/F$ lie in $[0,1]$, so $|h(W,Q)-h(W',Q')|\le|W-W'|+|Q-Q'|$. Since $\int W\,dA=\pi^2$ by \eqref{eq:passare},
\begin{equation}\label{eq:defect}
(1+\sqrt2)\pi^2-\Mv(c)=\Big(\sqrt2\pi^2-\int_\C Q\,dA\Big)+\int_\C h(W,Q)\,dA .
\end{equation}

We first collect pointwise bounds. Let $v_1=1/s$, $v_2=1/(s-1)$, $v_3=1/(s-c)$. Since $v_3-v_1=c/(s(s-c))$, we have $w_{13}=\operatorname{Im}(\bar v_1v_3)=\operatorname{Im}(c/(s-c))/\rho^2$, because $\operatorname{Im}(\bar v_1v_1)=0$. For $|s|\ge2\varepsilon$ we have $|s-c|\ge\rho/2$, hence
\begin{equation}\label{eq:Wbounds}
W\le\frac{2\varepsilon}{\rho^3},\qquad |W-A|\le\frac{|c|^2}{\rho^3|s-c|}\le\frac{2\varepsilon^2}{\rho^4},\qquad A=\frac{|\operatorname{Im}(c/s)|}{\rho^2}=\frac{\varepsilon|\sin(\varphi-\tau)|}{\rho^3},
\end{equation}
the second bound because $c/(s-c)-c/s=c^2/(s(s-c))$. Next, $w_{21}=\operatorname{Im}(\bar v_2v_1)=-w_{12}$, so $e=w_{23}+w_{12}=\operatorname{Im}(\bar v_2(v_3-v_1))$ satisfies $|e|\le\varepsilon\,|s-1|^{-1}|s|^{-1}|s-c|^{-1}$. For $|s|\le\frac12$ we have $|v_2|\le2$, and since $\bar v_1v_2+\bar v_1=\bar v_1\,s/(s-1)$ has modulus $|s-1|^{-1}\le2$, also $|w_{12}+\sin\varphi/\rho|\le2$. Hence, for $2\varepsilon\le|s|\le\frac12$, with $B=\sqrt2|\sin\varphi|/\rho$, the norm of the vector $(-\sin\varphi/\rho,\sin\varphi/\rho)$,
\begin{equation}\label{eq:Qbounds}
|e|\le\frac{4\varepsilon}{\rho^2},\qquad |Q-B|\le|w_{12}+\sin\varphi/\rho|+|w_{23}-\sin\varphi/\rho|\le4+\frac{4\varepsilon}{\rho^2},
\end{equation}
where we used $w_{23}-\sin\varphi/\rho=e-(w_{12}+\sin\varphi/\rho)$. Finally, for $|s|\le\frac12$ we have $Q\le|w_{12}|+|w_{23}|\le|v_1||v_2|+|v_2||v_3|\le2/\rho+2/|s-c|$.

The first term of \eqref{eq:defect} is small. Since $(w_{12},-w_{12})$ has norm $\sqrt2|w_{12}|$ and $\int\sqrt2|w_{12}|\,dA=\sqrt2\pi^2$ by \eqref{eq:passare}, we have $|\int Q\,dA-\sqrt2\pi^2|\le\int|Q-\sqrt2|w_{12}||\,dA\le\int\min(|e|,Q+\sqrt2|w_{12}|)\,dA$. On $\{|s|\le2\varepsilon\}$ we use $Q+\sqrt2|w_{12}|\le(1+\sqrt2)\,2/\rho+2/|s-c|$, whose integral over this disc is at most $8\pi(1+\sqrt2)\varepsilon+12\pi\varepsilon$, the disc being contained in $\{|s-c|\le3\varepsilon\}$. On $\{2\varepsilon\le|s|\le\frac12\}$, \eqref{eq:Qbounds} gives the bound $\int_{2\varepsilon}^{1/2}8\pi\varepsilon\rho^{-1}\,d\rho\le8\pi\varepsilon\log(1/\varepsilon)$. On $\{|s|\ge\frac12\}$ we have $|s-c|\ge|s|/2$, so $|e|\le2\varepsilon|s-1|^{-1}|s|^{-2}$, whose integral over this set is $C_0\varepsilon$ with $C_0=\int_{|s|\ge1/2}2|s-1|^{-1}|s|^{-2}\,dA<\infty$. Hence $|\sqrt2\pi^2-\int Q\,dA|\le C_1\varepsilon\log(1/\varepsilon)$ for an absolute constant $C_1$.

For the second term of \eqref{eq:defect} put $r=\varepsilon^{2/3}$ and $R=\varepsilon^{1/3}$; for $\varepsilon\le\frac18$ we have $2\varepsilon\le r<R\le\frac12$. On the disc $\{|s|\le r\}$ we use $h\le Q\le2/\rho+2/|s-c|$, whose integral is at most $4\pi r+4\pi(r+\varepsilon)\le12\pi\varepsilon^{2/3}$. On $\{|s|\ge R\}$ we use $h\le W\le2\varepsilon/\rho^3$, whose integral is at most $4\pi\varepsilon/R=4\pi\varepsilon^{2/3}$. On the annulus $\Omega=\{r\le|s|\le R\}$, the Lipschitz bound for $h$ together with \eqref{eq:Wbounds} and \eqref{eq:Qbounds} gives
$$\int_\Omega|h(W,Q)-h(A,B)|\,dA\le\int_\Omega\Big(\frac{2\varepsilon^2}{\rho^4}+4+\frac{4\varepsilon}{\rho^2}\Big)dA\le\frac{2\pi\varepsilon^2}{r^2}+4\pi R^2+8\pi\varepsilon\log\frac Rr\le6\pi\varepsilon^{2/3}+3\pi\varepsilon\log\frac1\varepsilon .$$
Finally, $h(A,B)\le B$ and $h(A,B)\le A$ give $\int_{|s|<r}h(A,B)\,dA\le\int_0^{2\pi}\sqrt2|\sin\varphi|\,d\varphi\cdot r=4\sqrt2\,\varepsilon^{2/3}$ and $\int_{|s|>R}h(A,B)\,dA\le\varepsilon\int_0^{2\pi}|\sin(\varphi-\tau)|\,d\varphi\int_R^\infty\rho^{-2}\,d\rho=4\varepsilon^{2/3}$. Collecting these estimates, and using $\varepsilon\log(1/\varepsilon)\le\frac3e\,\varepsilon^{2/3}$, which holds because $t^{1/3}\log(1/t)\le3/e$ for $0<t<1$, we obtain
$$(1+\sqrt2)\pi^2-\Mv(c)=\int_0^{2\pi}\!\!\int_0^\infty h(A,B)\,\rho\,d\rho\,d\varphi+R(c),\qquad|R(c)|\le C\varepsilon^{2/3}.$$

It remains to compute the integral of $h(A,B)$, which is finite by the bounds just used. Put $g(x)=1+x-(1+x^2)^{1/2}$, so that $h(A,B)=B\,g(A/B)$ wherever $\sin\varphi\ne0$. With $\rho=\sqrt\varepsilon\,u$ and $\lambda(\varphi)=|\sin(\varphi-\tau)|/(\sqrt2|\sin\varphi|)$ we have $A/B=\lambda/u^2$ and $B\rho\,d\rho=\sqrt2|\sin\varphi|\sqrt\varepsilon\,du$. By Tonelli's theorem and the substitution $u=\sqrt\lambda\,x$ (for $\lambda>0$; for $\lambda=0$ the inner integral vanishes),
$$\int_0^{2\pi}\!\!\int_0^\infty h(A,B)\,\rho\,d\rho\,d\varphi=\sqrt\varepsilon\int_0^{2\pi}\sqrt2|\sin\varphi|\int_0^\infty g(\lambda/u^2)\,du\,d\varphi=\sqrt\varepsilon\,K\int_0^{2\pi}\sqrt2|\sin\varphi|\sqrt{\lambda(\varphi)}\,d\varphi,$$
with $K=\int_0^\infty g(x^{-2})\,dx$, and $\sqrt2|\sin\varphi|\sqrt\lambda=2^{1/4}\sqrt{|\sin\varphi\sin(\varphi-\tau)|}$. The integral $K$ is finite, since $0\le g\le1$ and $g(t)\le t$.

To compute $K$, put $y=1/x$: $K=\int_0^\infty F(y)y^{-2}\,dy$ with $F(y)=1+y^2-(1+y^4)^{1/2}$, and since $F(y)/y\to0$ at $0$ and at $\infty$, integration by parts gives $K=2\int_0^\infty(1-y^2(1+y^4)^{-1/2})\,dy$. For $-1<a<1$, $\int_0^\infty y^a(1+y^4)^{-1/2}dy=\frac14B(\frac{a+1}4,\frac{1-a}4)$ \cite[3.241]{GR07}. Splitting at $y=1$ shows that this function and the function $\int_0^\infty y^a((1+y^4)^{-1/2}-y^{-2})dy$, defined for $1<a<5$, are both given by $\int_0^1y^a(1+y^4)^{-1/2}dy+\int_1^\infty y^a((1+y^4)^{-1/2}-y^{-2})dy+\frac1{1-a}$, real analytic on $(-1,5)\setminus\{1\}$; by analytic continuation the second equals $\frac14B(\frac{a+1}4,\frac{1-a}4)$ on $(1,5)$. At $a=2$, using $\Gamma(-\frac14)=-4\Gamma(\frac34)$, we get $\int_0^\infty(y^2(1+y^4)^{-1/2}-1)dy=-\Gamma(\frac34)^2/\sqrt\pi$, so $K=2\Gamma(\frac34)^2/\sqrt\pi=4\pi^{3/2}/\Gamma(\frac14)^2$ by $\Gamma(\frac14)\Gamma(\frac34)=\pi\sqrt2$ \cite[Section 12.14]{WW27}. This proves the expansion.

The function $\kappa$ is continuous by dominated convergence and positive, its integrand being positive off finitely many points. By the Cauchy--Schwarz inequality $\int_0^{2\pi}\sqrt{|\sin\varphi||\sin(\varphi-\tau)|}\,d\varphi\le(\int_0^{2\pi}|\sin\varphi|\,d\varphi)^{1/2}(\int_0^{2\pi}|\sin(\varphi-\tau)|\,d\varphi)^{1/2}=4$, with equality if and only if $|\sin\varphi|$ and $|\sin(\varphi-\tau)|$ are proportional, that is if and only if $\tau\in\{0,\pi\}$. Let $\kappa_0>0$ be the minimum of $\kappa$; then $\Mv(c)\le(1+\sqrt2)\pi^2-\kappa_0\sqrt\varepsilon+C\varepsilon^{2/3}<(1+\sqrt2)\pi^2$ for $0<\varepsilon<\min(\frac18,(\kappa_0/C)^6)$. The statements near $1$ and $\infty$ follow from the invariance of $\Mv$ under $c\mapsto1-c$ and $c\mapsto1/c$ (Proposition \ref{prop:sym12}).
\end{proof}

By Proposition \ref{prop:sym12} the same expansion holds as $c\to1$ with $\kappa(\arg(1-c))\sqrt{|c-1|}$ and as $c\to\infty$ with $\kappa(\arg c)/\sqrt{|c|}$, since $\kappa(-\tau)=\kappa(\tau)$. Thus $\Mv$ approaches $(1+\sqrt2)\pi^2$ strictly from below near the degenerate values.

\subsection{Critical points forced by symmetry}

\begin{lemma}\label{lem:equiv}
Let $F$ be a $C^1$ function on an open set $U\subset\R^m$ invariant under a finite group $H$ of diffeomorphisms fixing $x\in U$. If $H$ has no nonzero fixed vector in $T_xU$, then $dF(x)=0$. If moreover $m=2$, $F$ is $C^2$ and $H$ contains a rotation of order $3$ of $T_xU$, then the Hessian of $F$ at $x$ is a multiple of the identity.
\end{lemma}

\begin{proof}
Invariance gives $dF(x)\circ dh(x)=dF(x)$, so the vector dual to $dF(x)$ is $H$-fixed, hence zero. A quadratic form on $\R^2$ is a multiple of $|x|^2$ plus a traceless part on which a rotation by $\theta$ acts as a rotation by $2\theta$, and invariance under rotation by $2\pi/3$ kills the traceless part.
\end{proof}

\begin{proposition}\label{prop:crit}
The equilateral values $c=e^{\pm i\pi/3}$ are critical points of $\Mv$ with isotropic Hessian. Every point of the circles and lines $|c|=1$, $|c-1|=1$, $\operatorname{Re}c=\frac12$ off the real axis is critical in the direction normal to the curve. The harmonic values $c\in\{-1,\frac12,2\}$ are critical points of the restriction of $\Mv$ to the real axis.
\end{proposition}

\begin{proof}
The element $c\mapsto1/(1-c)$ of $\Sg$ fixes $e^{i\pi/3}$ and its derivative there is $1/(1-c)^2=e^{2\pi i/3}$, a rotation of order three; Lemma \ref{lem:equiv} and Proposition \ref{prop:reg} give the first assertion, and conjugation handles $e^{-i\pi/3}$. A point of an isosceles curve is fixed by a reflection of $\Sg$, for instance $c\mapsto1/\bar c$ for $|c|=1$, which reverses the normal direction. The restriction of $\Mv$ to $(-\infty,0)$ is smooth and invariant under $c\mapsto1/c$, which fixes $-1$ with derivative $-1$; the other harmonic values are images of $-1$ under $\Sg$.
\end{proof}

\subsection{Conclusions on the extremes}

\begin{theorem}\label{thm:linesummary}
The function $\Mv$ on $\C\setminus\{0,1\}$ is continuous, invariant under $\Sg$, real analytic on $\C\setminus\R$, satisfies $\sqrt3\,\pi^2<\Mv<3\pi^2$, and tends to $(1+\sqrt2)\pi^2$ as $c\to0,1,\infty$ at the rate of Proposition \ref{prop:rate}; in particular $\Mv(c)<(1+\sqrt2)\pi^2$ on a neighbourhood of $\{0,1,\infty\}$ in the Riemann sphere. Its critical points include $e^{\pm i\pi/3}$, with isotropic Hessian, and the harmonic values are critical along the real axis.
\end{theorem}

\begin{conjecture}\label{conj:lines}
For every $c\in\C\setminus\{0,1\}$ one has $\Mv(-1)\le\Mv(c)<(1+\sqrt2)\pi^2$, with equality on the left only for $c\in\{-1,\frac12,2\}$. The only other critical points of $\Mv$ off the real axis are the strict local minima $e^{\pm i\pi/3}$ and six saddles at the orbit of $e^{2\pi i/3}$, the isosceles triangle with angles $120^\circ,30^\circ,30^\circ$.
\end{conjecture}

Outside a large disc the function $\Mv$ is controlled by Theorem \ref{thm:degen} and Proposition \ref{prop:rate} via $c\mapsto1/c$, and the inequality $\Mv(c)<(1+\sqrt2)\pi^2$ of the conjecture is proved near $c=0,1,\infty$ by Proposition \ref{prop:rate}. The position of the saddles is not forced by symmetry: the only element of $\Sg$ fixing $e^{2\pi i/3}$ is the reflection $c\mapsto1/\bar c$ fixing the arc $|c|=1$, so Lemma \ref{lem:equiv} only gives the vanishing of the derivative in the direction normal to the arc. Granting the conjecture, the minimum of the multivolume is attained, exactly at the harmonic values, and the maximum is not attained. For the area of a nonreal line, which equals $\Mv(c)$, neither extreme would then be attained: the infimum $\Mv(-1)$ is approached only as $c$ tends to a harmonic real value and the supremum $(1+\sqrt2)\pi^2$ only as $c\to0,1,\infty$, while the most symmetric nonreal lines $c=e^{\pm i\pi/3}$ realise a local minimum. For the area of a real line, which equals $\frac12\Mv(c)$, the minimum would be attained at the harmonic values, for instance by $\{(s,1-s,1+s)\}$, and the supremum $\frac{1+\sqrt2}2\pi^2$ would not be attained. Unconditionally, the area is discontinuous across the real axis, where it drops by the factor $\frac12$ predicted by Corollary \ref{cor:half}, which for real lines is an equality at the level of $\MVol$ but whose last bound $\frac32\pi^2$ is not sharp, and the multivolume of lines in $(\C^*)^3$ is not constant, by Proposition \ref{prop:rate}, unlike the multivolume of half-dimensional linear spaces, which is always $\pi^{2k}$ \cite[Proposition 8]{Mik17}; its variation is governed by the shape of the triangle formed by the three finite points where the line meets the coordinate planes.

\section{\texorpdfstring{$k$-planes in $(\C^*)^{2k+1}$}{k-planes in (C*)\textasciicircum (2k+1)}}\label{sec:kplanes}

Let $V\subset(\C^*)^{2k+1}$ be the intersection with the torus of a generic affine $k$-plane of $\C^{2k+1}$, defined by $k+1$ affine linear equations, so that $\Delta_m\subset\Sigma_{2k+1}$. The $2k+1$ subsets $I$ are the complements of single indices $j$. For such $I$, the system \eqref{eq:doubled} consists of $k+1$ affine linear equations in $(r_I,z_j)$ and $k+1$ in $(r_I,w_j)$; eliminating $z_j$ and $w_j$ leaves $2k$ affine linear equations in the $2k$ unknowns $r_I$, with exactly one solution for generic coefficients. Hence $\alpha_I=1$ and Theorem \ref{thm:main} gives $\Hm^{2k}(\A(V))\le(2k+1)\pi^{2k}$. For generic $V$ the direction $e_j$ is not in the direction space of the $k$-plane, so $\pi_I|_V$ is injective and $\pi_I(V)$ is, up to a hyperplane section, a generic $k$-plane in $(\C^*)^{2k}$, whose multivolume is $\pi^{2k}$ by \cite[Proposition 8]{Mik17}, in accordance with the volume $\pi^{2k}$ of its coamoeba computed in \cite{NP17}. Hence $\MVol_I(V)=\pi^{2k}$ for every $I$ and $\sum_I\MVol_I(V)=(2k+1)\pi^{2k}$.

For the lower bound, the amoeba map of a generic $k$-plane $W\subset(\C^*)^{2k}$ covers its amoeba at most $\binom{2k}{k}$ times \cite[Corollary 3]{Mik17}, so $\Vol_{2k}(\A(W))\ge\pi^{2k}/\binom{2k}{k}$, and Theorem \ref{thm:main} gives $\Hm^{2k}(\A(V))\ge\sqrt{2k+1}\,(k!)^2\pi^{2k}/(2k)!$. For $k=1$ this is $\frac{\sqrt3}2\pi^2$, the bound of Theorem \ref{thm:line} for real lines. These general bounds are superseded by Corollary \ref{cor:kplanes2}: for almost every real $k$-plane each projected part equals $\pi^{2k}/2$ and $\frac{\sqrt{2k+1}}2\pi^{2k}\le\Hm^{2k}(\A(V))=\frac12\MVol(V)\le\frac{2k+1}2\pi^{2k}$, and for almost every complex $k$-plane each projected part equals $\pi^{2k}$ and $\sqrt{2k+1}\,\pi^{2k}\le\Hm^{2k}(\A(V))\le(2k+1)\pi^{2k}$. For $k=2$ the real statement holds for every real plane satisfying the condition of Theorem \ref{thm:conic}, and the real planes with six lines in general position which violate it are the squares of lines, for which each projected part equals $\pi^4/4$ and $\frac{\sqrt5}4\pi^4\le\Hm^4(\A(V))=\frac14\MVol(V)\le\frac54\pi^4$ (Corollary \ref{cor:squareproj}), with supremum at least $\frac{1+2\sqrt2}4\pi^4$ by Corollary \ref{cor:supsq}. For every $k$, the $k$-th symmetric powers of real lines are real $k$-planes of $(\C^*)^{2k+1}$ with projected parts $\pi^{2k}/2^k$ (Corollary \ref{cor:squareproj}), and those of generic nonreal lines have projected parts $\pi^{2k}$ (Corollary \ref{cor:nonrealproj}). In contrast with Corollary \ref{cor:kplanes2}, for real generic $k$-planes in $(\C^*)^{2k}$ with $k\ge3$ the amoeba map is not $2^k$-to-one (Corollary \ref{cor:fibre4}), contrary to \cite[Proposition 5.6]{NP17}; the value $\pi^{2k}/2^k$ of \cite[Theorem 1.2(ii)]{NP17} for the volume is not decided by Corollary \ref{cor:fibre4}, since fibres of more than $2^k$ points also occur, and it remains open (Section \ref{sec:open}).

\section{\texorpdfstring{The chain of inequalities for $n>2k$}{The chain of inequalities for n>2k}}\label{sec:chain}

From now on $n>2k$. We first make the third inequality of \eqref{eq:chain} explicit. Fix $|I|=2k$ and $\varphi\in\T^I_\pi$, and let $\hat I$ be the complement of $I$. Mikhalkin's complexification of the fibre $V\cap\Arg_{\pi,I}^{-1}(\varphi)$ is the square system in the unknowns $(r_I,z_{\hat I},w_{\hat I})\in(\C^*)^{2n-2k}$ of Section \ref{sec:counting}, namely the system \eqref{eq:doubled}, with $F_m(r_I,z_{\hat I})=f_m(e^{i\varphi_I}r_I,z_{\hat I})$; by Bernstein's theorem it has at most $\alpha_I$ isolated solutions in the torus (Proposition \ref{prop:fibre}). The antiholomorphic involution $\iota_\varphi(r_I,z_{\hat I},w_{\hat I})=(\bar r_I,\bar w_{\hat I},\bar z_{\hat I})$ preserves the solution set of \eqref{eq:doubled}, since $F_m(\bar r_I,\bar w_{\hat I})=\overline{G_m(r_I,w_{\hat I})}$ and $G_m(\bar r_I,\bar z_{\hat I})=\overline{F_m(r_I,z_{\hat I})}$.

\begin{proposition}\label{prop:BI}
The map $z\mapsto(r_I,z_{\hat I},\bar z_{\hat I})$, where $z_I=e^{i\varphi_I}r_I$, is a bijection from the fibre $V\cap\Arg_{\pi,I}^{-1}(\varphi)$ onto the set of solutions of \eqref{eq:doubled} fixed by $\iota_\varphi$. Consequently $\MVol_I(V)=\pi^{2k}\alpha_I$ if and only if the following condition holds. $(\mathrm B_I)$: for almost every $\varphi\in\T^I_\pi$, the system \eqref{eq:doubled} has exactly $\alpha_I$ isolated solutions in the torus, and all of them are fixed by $\iota_\varphi$. For $n=2k$ condition $(\mathrm B_I)$ is condition $(\mathrm B)$ of Theorem \ref{thm:char} together with $(\mathrm A)$.
\end{proposition}

\begin{proof}
A solution fixed by $\iota_\varphi$ has $r_I\in(\R^*)^{2k}$ and $w_{\hat I}=\bar z_{\hat I}$, and then $F_m(r_I,z_{\hat I})=0$ says that $z=(e^{i\varphi_I}r_I,z_{\hat I})$ lies in $V$; its rolled arguments in $I$ are $\varphi_I$ because $r_I$ is real. Conversely a point of the fibre has this form with $r_I$ real, and gives the fixed solution $(r_I,z_{\hat I},\bar z_{\hat I})$, because $G_m(r_I,\bar z_{\hat I})=\overline{F_m(r_I,z_{\hat I})}=0$. The two constructions are inverse to each other. By the proof of Proposition \ref{prop:fibre}, for almost every $\varphi$, namely for the regular values of $\Arg_{\pi,I}|_V$, these fixed solutions are nondegenerate, hence isolated. Thus the fibre has at most $\alpha_I$ points, and by \eqref{eq:rolled} $\MVol_I(V)=\pi^{2k}\alpha_I$ if and only if the fibre has exactly $\alpha_I$ points for almost every $\varphi$. If $(\mathrm B_I)$ holds, then for almost every $\varphi$ the $\alpha_I$ isolated solutions are fixed, so the fibre has at least $\alpha_I$ points, hence exactly $\alpha_I$. Conversely, if the fibre has $\alpha_I$ points for a regular value $\varphi$, the system has $\alpha_I$ isolated fixed solutions; since it has at most $\alpha_I$ isolated solutions in all, they are all of them, and $(\mathrm B_I)$ holds. For $n=2k$, \eqref{eq:doubled} reads $f(e^{i\varphi}r)=0$, $\bar f(e^{-i\varphi}r)=0$, whose solutions $u=e^{i\varphi}r$ are the points of $V\cap\sigma_\varepsilon(V)$, $\varepsilon=e^{2i\varphi}$, while $\iota_\varphi$ corresponds to $\sigma_\varepsilon$; the number of these points is $\alpha(V)$ for almost every $\varepsilon$ by Lemma \ref{lem:alpha}, so $(\mathrm B_I)$ amounts to $\alpha(V)=M$ and $(\mathrm B)$.
\end{proof}

\begin{lemma}\label{lem:CBeq}
Let $V_0$ be a connected component of $V$ such that $J\not\equiv0$ on $V_0$ and $\int_{V_0}J=\sum_I\int_{V_0}J_I$. Then there are $I$ and $t\in(\C^*)^n$ such that $V_0\subset T^I_t$.
\end{lemma}

\begin{proof}
By \eqref{eq:CBJ} the integrand $\sum_IJ_I-J$ is nonnegative and continuous, and its integral over $V_0$ vanishes, so it vanishes identically on $V_0$. Hence at every point of $V_0$ at most one $J_I$ is nonzero, that is, $J_I^2J_{I'}^2\equiv0$ on $V_0$ for $I\ne I'$. The ring of real analytic functions on the connected manifold $V_0$ has no zero divisors, since the zero set of a real analytic function which is not identically zero has empty interior. So for every pair $I\ne I'$ one of $J_I$, $J_{I'}$ vanishes identically on $V_0$. Since $J\not\equiv0$, some $J_I$ is not identically zero, and then $J_{I'}\equiv0$ for all $I'\ne I$. Let $U=\{z\in V_0:J(z)>0\}$, a nonempty open dense subset of $V_0$. At $z\in U$ the image $d\Log(T_zV_0)$ is a $2k$-dimensional subspace of $\R^n$ whose only nonzero Pl\"ucker coordinate is the one indexed by $I$, so it equals $\R^I$. Hence for $j\notin I$ the differential of $\log|z_j|$ vanishes on $T_zV_0$ for every $z\in U$, and by density for every $z\in V_0$. Thus $|z_j|$ is constant on the connected manifold $V_0$, and the holomorphic function $z_j|_{V_0}$, having constant modulus, is constant by the open mapping theorem. Calling $t_j$ this constant for $j\notin I$, and putting $t_j=1$ for $j\in I$, we get $V_0\subset T^I_t$.
\end{proof}

\begin{proof}[Proof of Theorem \ref{thm:strict}]
By \eqref{eq:area}, $\Hm^{2k}(\A(V))=\int_{\A(V)}1\,d\Hm^{2k}\le\int_{\A(V)}N\,d\Hm^{2k}=\MVol(V)$, since $N\ge1$ on $\A(V)$ and $N=0$ off it. The right-hand side is finite, so equality holds if and only if $N=1$ almost everywhere on $\A(V)$. The second inequality is the integral of \eqref{eq:CBJ}. If it is an equality, then it is an equality on each of the finitely many connected components of $V$, and Lemma \ref{lem:CBeq} shows that each component with $J\not\equiv0$ lies in some $T^I_t$. This proves the second assertion. The third assertion is Proposition \ref{prop:BI} applied to each $I$, the sum of the bounds $\MVol_I(V)\le\pi^{2k}\alpha_I$ of Proposition \ref{prop:fibre} being an equality if and only if each of them is.

Assume now that $\sum_I\alpha_I>0$ and that $\Hm^{2k}(\A(V))=\pi^{2k}\sum_I\alpha_I$. Then the three inequalities of \eqref{eq:chain} are equalities, and $\MVol(V)=\pi^{2k}\sum_I\alpha_I>0$. Let $V_0$ be a connected component of $V$ with $J\not\equiv0$; it exists, since $\MVol(V)>0$. By the second inequality and Lemma \ref{lem:CBeq}, $V_0\subset T^I_t$ for some $I$ and $t$. The projection $\pi_I$ maps $T^I_t$ isomorphically onto $(\C^*)^I$, and $p_I$ maps the affine space $\Log(T^I_t)=\{x:x_j=\log|t_j|\text{ for }j\notin I\}$ isometrically onto $\R^I$, compatibly with $\Log$. Hence the counting function of $\Log|_{V_0}$ at $x\in\Log(T^I_t)$ equals the counting function of $\Log|_W$ at $p_I(x)$, where $W=\pi_I(V_0)\subset(\C^*)^I$ is a closed complex submanifold of dimension $k$, closed because $V_0$ is closed in $(\C^*)^n$ and contained in the closed set $T^I_t$. The amoeba of $W$ is isometric to $\A(V_0)$, so its Lebesgue measure is at most $\MVol(V_0)<\infty$. By Lemma \ref{lem:parity} and Sard's theorem, the counting function of $\Log|_{V_0}$ is at least $2$ at $\Hm^{2k}$-almost every point of $\A(V_0)$. The counting function of $\Log|_V$ is at least that of $\Log|_{V_0}$, so $N\ge2$ almost everywhere on $\A(V_0)$. On the other hand, $\Hm^{2k}(\A(V_0))>0$, since $\Log|_{V_0}$ is an immersion on the nonempty open set $\{J>0\}\cap V_0$, hence an embedding on a small ball in it, whose image is an embedded $2k$-dimensional disc of positive measure. This contradicts the first equality, which says $N=1$ almost everywhere on $\A(V)$.
\end{proof}

\section{Injective amoeba maps and varieties of maximal volume}\label{sec:inj}

\begin{proof}[Proof of Theorem \ref{thm:emb}]
The equivalence of $(\mathrm I_1)$ with the statement that $N=1$ almost everywhere on $\A(V)$ is the first assertion of Theorem \ref{thm:strict}, and this statement is $(\mathrm I_2)$, since the points of $V$ fixed by $\sigma'_x$ are the points of $V\cap\Log^{-1}(x)$. If $(\mathrm I_3)$ holds, then $N=1$ on $\A(V)\setminus Z$, and $Z$ is null, so $(\mathrm I_2)$ holds.

Assume $(\mathrm I_2)$. We regard $V$ as a real algebraic subset of $(\C^*)^n\subset\R^{2n}$ and introduce the double point set $D=\{z\in V:\ \text{there is }w\in V\setminus\{z\}\text{ with }|w_j|=|z_j|\text{ for all }j\}.$ It is the projection to the first factor of the semialgebraic set $\{(z,w)\in V\times V:z\ne w,\ |z_j|^2=|w_j|^2\text{ for all }j\}$, so it is semialgebraic by the Tarski--Seidenberg theorem \cite[Proposition 2.2.7]{BCR98}. The set $V^+=\{z\in V:J(z)>0\}$ is semialgebraic and open, by Section \ref{sec:prelim}. We claim that $D\cap V^+$ has dimension smaller than $2k$. Otherwise it contains a nonempty open subset $O$ of the $2k$-dimensional manifold $V$ \cite[Section 2.8]{BCR98}. Shrinking $O$, we may assume that $\Log|_O$ is an embedding, since $J>0$ on $O$, so $\Hm^{2k}(\Log(O))>0$. But every point $x\in\Log(O)$ has, besides its preimage in $O$, a second preimage in $V$, so $N\ge2$ on $\Log(O)$, contradicting $(\mathrm I_2)$. Let $E=\{z\in V:J(z)=0\}\cup\overline{D\cap V^+}$. It is closed, and $\overline{D\cap V^+}$ has dimension smaller than $2k$ \cite[Proposition 2.8.13]{BCR98}, so it has volume zero in $V$, being a finite union of submanifolds of dimension at most $2k-1$. By \eqref{eq:area}, $\Hm^{2k}\big(\Log(E)\big)\le\int_{\{J=0\}}J\,d\mathrm{vol}+\int_{\overline{D\cap V^+}}J\,d\mathrm{vol}=0$. Put $Z=\Log(E)$. It is closed, because $\Log|_V$ is proper, hence closed, and it is contained in $\A(V)$. Let $U=V\setminus\Log^{-1}(Z)$, an open subset of $V$ with $\Log(U)=\A(V)\setminus Z$. Every $z\in U$ satisfies $J(z)>0$ and $z\notin D$, because $z\notin E$; thus $\Log|_U$ is an injective immersion, and no point of $V$ other than $z$ has the same image. The restriction $\Log|_U:U\to\R^n\setminus Z$ is proper, being the restriction of the proper map $\Log|_V$ to the full preimage of $\R^n\setminus Z$. A proper injective immersion is an embedding with closed image. Hence $\A(V)\setminus Z$ is a closed embedded submanifold of $\R^n\setminus Z$ of dimension $2k$, and $\Log:U\to\A(V)\setminus Z$ is a diffeomorphism. This is $(\mathrm I_3)$.

Assume now that the three conditions hold and that $\MVol(V)>0$. If $\G_V$ is infinite, then $\Hm^{2k}(\A(V))=0$ by Lemma \ref{lem:infinite}, while $\Hm^{2k}(\A(V))=\MVol(V)>0$ by $(\mathrm I_1)$. If $\G_V$ is finite and nontrivial, then by Lemma \ref{lem:count}, applied to $\Gamma=\G_V$, the fibre of $\Log|_V$ over $\Hm^{2k}$-almost every point of $\A(V)$ has at least $|\G_V|\ge2$ points, contradicting $(\mathrm I_2)$. So $\G_V$ is trivial, and $V$ is not real up to a multiplicative translation by Lemma \ref{lem:real}. Finally, if a component $V_0$ with $J\not\equiv0$ lies in a translate of a coordinate subtorus, the proof of Theorem \ref{thm:strict} shows that $N\ge2$ almost everywhere on $\A(V_0)$, a set of positive measure, which contradicts $(\mathrm I_2)$.
\end{proof}

\begin{proof}[Proof of Theorem \ref{thm:max}]
For every $V\in\Vc(\Delta)$, $\Hm^{2k}(\A(V))\le\MVol(V)\le m^*(\Delta)$ by \eqref{eq:chain}, which gives the first inequality. The equality $\Hm^{2k}(\A(V))=m^*(\Delta)$ holds if and only if both inequalities $\Hm^{2k}(\A(V))\le\MVol(V)$ and $\MVol(V)\le m^*(\Delta)$ are equalities, that is, by Theorem \ref{thm:emb}, if and only if the conditions of Theorem \ref{thm:emb} hold and $\MVol(V)=m^*(\Delta)$.

Assume the injectivity property, and let $V_0\in\Vc(\Delta)$. There is a sequence $V_l\in\Vc(\Delta)$ converging to $V_0$ in the space of coefficients, each $V_l$ satisfying the conditions of Theorem \ref{thm:emb}, so that $\Hm^{2k}(\A(V_l))=\MVol(V_l)$. By Lemma \ref{lem:lsc}, $\sup_{\Vc(\Delta)}\Hm^{2k}(\A(V))\ge\liminf_{l\to\infty}\MVol(V_l)\ge\MVol(V_0).$ Taking the supremum over $V_0$ gives $\sup\Hm^{2k}(\A(V))\ge m^*(\Delta)$, hence equality. A variety of maximal $2k$-volume then satisfies $\Hm^{2k}(\A(V))=m^*(\Delta)$, and the first part of the proof applies. If $m^*(\Delta)>0$, Theorem \ref{thm:emb} gives the triviality of $\G_V$ and the non-reality, and, since no component with $J\not\equiv0$ lies in a translate of a coordinate subtorus, the second inequality of \eqref{eq:chain} is strict by Theorem \ref{thm:strict}. If $m^*(\Delta)=0$, all amoebas have measure zero and there is nothing to prove.
\end{proof}

\section{Generic injectivity of the amoeba map}\label{sec:geninj}

Throughout this section $n>2k$ and $(\mathrm L_1)$, $(\mathrm L_2)$ hold, and we use the notation introduced before Theorem \ref{thm:geninj}. For $c\in\Vc(\Delta)$ we write $J_c(z)$ for the $2k$-dimensional Jacobian of $\Log|_{V_c}$ at $z$ and $N_c$ for the counting function of $\Log|_{V_c}$. For $u\in\T^n=(S^1)^n$ and $a\in\Z^n$ we write $u^a=\prod_ju_j^{a_j}$. For a subgroup $L\subset\Z^n$ we put $G(L)=\{u\in\T^n:u^v=1\text{ for all }v\in L\}$. This is a closed subgroup of $\T^n$, isomorphic to the group of characters of $\Z^n/L$, so it is a compact abelian Lie group of dimension $n-\rk L$. We write $G_m=G(L_m)$ and $G_S=G(L_S)=\bigcap_{m\in S}G_m$, with $G_\emptyset=\T^n$.

Put $F(c,z)=(f_{c,1}(z),\dots,f_{c,n-k}(z))$ on $\C^A\times(\C^*)^n$ and let $U$ be the set of $(c,z)$ with $F(c,z)=0$ and $\rk d_zF(c,z)=n-k$. By Lemma \ref{lem:U}, $U$ is a complex submanifold of dimension $|A|+k$, where $|A|=\sum_m|A_m|$, the projection $\pi:U\to\C^A$ is a submersion, $\pi^{-1}(c)=\{c\}\times V_c$ for $c\in\Vc(\Delta)$, and $\Vc(\Delta)$ contains a nonempty Zariski open subset of $\C^A$. For a subset $S\subsetneq\{1,\dots,n-k\}$ put $s=|S|$, $d_S=\rk L_S$, $G_S^\circ=G_S\setminus\bigcup_{m\notin S}G_m$, $\mathcal Y_S=U\times G_S^\circ$ and $\Psi_S(c,z,u)=(f_{c,m}(z/u))_{m\notin S}\in\C^{n-k-s}$. The set $G_S^\circ$ is open in the compact Lie group $G_S$, so $\mathcal Y_S$ is a real manifold of dimension $2|A|+2k+(n-d_S)$, and $\Psi_S$ is real analytic. For $c\in\C^A$ let $W_{S,c}=\{(z,u):(c,z,u)\in\Psi_S^{-1}(0)\}$ and let $Z_{S,c}$ be the set of $z$ with $(z,u)\in W_{S,c}$ for some $u$. The transversality argument is contained in the following statement, proved in Appendix \ref{app:transv} (Lemmas \ref{lem:sep} to \ref{lem:dim}): $\Psi_S$ is a submersion, and there is a Lebesgue null set $\mathcal N_S\subset\C^A$ such that for $c\notin\mathcal N_S$ the set $W_{S,c}$ is empty or a submanifold of dimension $e_S=4k-n-d_S+2s\le2k-1$, so that $Z_{S,c}$ has $2k$-dimensional volume zero in $V_c$. The inequality $e_S\le2k-1$ is exactly condition $(\mathrm L_2)$.

\begin{proof}[Proof of Theorem \ref{thm:geninj}]
Let $\mathcal N=\bigcup_S\mathcal N_S$, a null set, the union running over the finitely many $S\subsetneq\{1,\dots,n-k\}$. Let $c\in\Vc(\Delta)\setminus\mathcal N$ and consider the double point set $D_c=\{z\in V_c:\ \text{there is }w\in V_c\setminus\{z\}\text{ with }\Log w=\Log z\}.$ Let $z\in D_c$ with partner $w$, and put $u=z/w$. Then $|u_j|=|z_j|/|w_j|=1$, so $u\in\T^n$, and $u\ne1$. Let $S=\{m:u\in G_m\}$. If $S$ were all of $\{1,\dots,n-k\}$, then $u\in G(L_{\{1,\dots,n-k\}})=G(\Z^n)=\{1\}$ by $(\mathrm L_1)$, a contradiction. So $S$ is proper, $u\in G_S^\circ$, $(c,z)\in U$, and $\Psi_S(c,z,u)=(f_{c,m}(w))_{m\notin S}=0$. Hence $z\in Z_{S,c}$, and $D_c\subset\bigcup_SZ_{S,c}$ has volume zero in $V_c$ by Lemma \ref{lem:dim}. If $x\in\A(V_c)$ satisfies $N_c(x)\ge2$, every preimage of $x$ lies in $D_c$, so $\{N_c\ge2\}\subset\Log(D_c)$, and by \eqref{eq:area} $\Hm^{2k}(\{N_c\ge2\})\le\int_{D_c}J_c\,d\mathrm{vol}=0$. Thus the amoeba map of $V_c$ is injective almost everywhere for $c\in\Vc(\Delta)\setminus\mathcal N$.

We now replace $\mathcal N$ by a proper real algebraic subset. Let $\mathcal B_0$ be the set of $c\in\Vc(\Delta)$ with $\Hm^{2k}(\{N_c\ge2\})>0$, a null set since $\mathcal B_0\subset\mathcal N$. In real coordinates, $\Vc(\Delta)$ is semialgebraic by Tarski--Seidenberg \cite[Proposition 2.2.7]{BCR98}, $J_c(z)^2$ is a rational function of $(c,z)$, and $D^+=\{(c,z):z\in V_c,\ J_c(z)>0,\ \exists\,w\in V_c\setminus\{z\},\ |w_j|^2=|z_j|^2\ \forall j\}$ is semialgebraic. As in the proof of Theorem \ref{thm:emb}, $c\in\mathcal B_0$ if and only if $\dim D^+_c=2k$, the points with $J_c=0$ contributing a null set by \eqref{eq:area}. By Hardt's triviality theorem \cite[Theorem 9.3.2]{BCR98} this dimension is a semialgebraic function of $c$, so $\mathcal B_0$ is semialgebraic and null, and its Zariski closure $\mathcal B$ is a proper real algebraic subset \cite[Proposition 2.8.2]{BCR98}. This proves the first assertion.

A proper real algebraic subset of $\R^{2|A|}$ has empty interior. By Lemma \ref{lem:U}, $\Vc(\Delta)$ contains a nonempty Zariski open subset $W'$ of $\C^A$, and $W'\setminus\mathcal B$ is dense in $\C^A$, so every point of $\Vc(\Delta)$ is a limit of points of $\Vc(\Delta)\setminus\mathcal B$. This is the injectivity property. The remaining assertions follow from Theorem \ref{thm:max}, whose proof uses only the injectivity property and the lower semicontinuity of the multivolume Lemma \ref{lem:lsc}. That lemma is proved by perturbing compact pieces of $V_{c_0}$ by the implicit function theorem, which only requires $c$ to stay in $\Vc(\Delta)$.
\end{proof}

\begin{example}\label{ex:conditions}
If every $\Delta_m$ contains $0,e_1,\dots,e_n$, then $L_m=\Z^n$ for all $m$, so $(\mathrm L_1)$ holds, and $(\mathrm L_2)$ reads $n\ge2s+2k-n+1$ for $s\le n-k-1$, which holds because $2s+2k-n+1\le2(n-k-1)+2k-n+1=n-1$. For the lines of $(\C^*)^3$ with $\Delta_1=\operatorname{conv}(0,e_1,e_2)$ and $\Delta_2=\operatorname{conv}(0,e_1,e_3)$, $L_1=\Z e_1+\Z e_2$, $L_2=\Z e_1+\Z e_3$, $L_1+L_2=\Z^3$, and $(\mathrm L_2)$ for $S=\{1\}$ or $\{2\}$ reads $2\ge2+2-3+1=2$. In this case the exceptional set $\mathcal B$ contains the real lines, $c\in\R$, whose amoeba maps are two-to-one, see Section \ref{sec:injlines}.
\end{example}

\begin{corollary}\label{cor:planes}
Let $n>2k$ and let all $\Delta_m$ be the standard simplex $\operatorname{conv}(0,e_1,\dots,e_n)$, so that $\Vc(\Delta)$ consists of the affine $k$-planes of $\C^n$, intersected with $(\C^*)^n$, cut out by equations with all coefficients nonzero. Then $\Vc(\Delta)$ has the injectivity property, $\sup\Hm^{2k}(\A(V))=m^*(\Delta)$, and a plane has maximal $2k$-volume if and only if its amoeba map is an embedding off a null set and its multivolume is maximal. If $n=2k+1$, then $\sqrt{2k+1}\,\pi^{2k}\le\Hm^{2k}(\A(V))<(2k+1)\pi^{2k}$ for the planes outside a proper real algebraic subset.
\end{corollary}

\begin{proof}
The conditions $(\mathrm L_1)$ and $(\mathrm L_2)$ hold by Example \ref{ex:conditions}, so the first assertions follow from Theorem \ref{thm:geninj}. For $n=2k+1$, outside a proper real algebraic subset the amoeba map is injective almost everywhere, so $\Hm^{2k}(\A(V))=\MVol(V)$, and each projection $\pi_I|_V$ is injective onto a generic $k$-plane of $(\C^*)^{2k}$, whose multivolume is $\pi^{2k}$ \cite[Proposition 8]{Mik17}. The Cauchy--Binet formula \eqref{eq:CBJ} and the Cauchy--Schwarz inequality give $J\ge(2k+1)^{-1/2}\sum_IJ_I$, hence $\MVol(V)\ge(2k+1)^{-1/2}\sum_I\MVol_I(V)=\sqrt{2k+1}\,\pi^{2k}$, and the upper bound is Theorem \ref{thm:strict}.
\end{proof}

\section{Curves, eliminants and mixed fibre polytopes}\label{sec:curves}

For $k=1$ the multivolumes $\MVol_I$ can be bounded through the Newton polygons of the plane curves obtained by elimination. Let $C\subset(\C^*)^n$ be a smooth complete intersection curve defined by $f_1,\dots,f_{n-1}$ with generic coefficients. For $I=\{i,j\}$ the direct image $(\pi_I)_*[C]$ is the divisor of a Laurent polynomial $g_I$, and by Esterov and Khovanskii \cite{EK08} (see also Sturmfels, Tevelev and Yu \cite{STY07} and Billera and Sturmfels \cite{BS92}) the Newton polygon $P_I$ of $g_I$ is the mixed fibre polytope of $\Delta_1,\dots,\Delta_{n-1}$ with respect to $p_I$. If $g_I=h_I^{d_I}$ with $h_I$ irreducible, $d_I$ being the degree of $\pi_I|_C$ onto its image, then $\MVol_I(C)=d_I\MVol(\{h_I=0\})$, and Mikhalkin's bound $\MVol(\{h=0\})\le\pi^2\MV(\Delta_h,\Delta_h)=2\pi^2\Vol(\Delta_h)$ \cite[Theorem 1]{Mik17} gives $\MVol_I(C)\le2\pi^2d_I\Vol(\Delta_{h_I})\le2\pi^2\Vol(P_I)$. The first inequality of Theorem \ref{thm:main} therefore gives

$$\Hm^2(\A(C))\le2\pi^2\sum_{1\le i<j\le n}\Vol(P_{ij}),$$

the analogue for curves in $(\C^*)^n$ of the Passare--Rullg{\aa}rd bound \cite{PR04}, and, if $C$ is real up to translation, the same bound with $\pi^2$ in place of $2\pi^2$ by Theorem \ref{thm:sym}. For the lines $L_c$ the three polygons are unit triangles and the bound is $3\pi^2$.

\section{Open problems}\label{sec:open}

The estimates above depend on the Newton polytopes through the mixed volumes $\alpha_I$, which are defined through a chosen set of generators, so one may take the minimum over all presentations of $V$ as a smooth complete intersection. Two losses occur in the argument. The first is the inequality $\sum_IJ_I(T)\ge1$ of \eqref{eq:cb}, sharp only for coordinate planes; averaging over all orthogonal projections onto $2k$-planes, as in the Cauchy--Crofton formula \cite[Theorem 3.2.26]{Fed69}, would give an exact formula for $\Hm^{2k}(\A(V))$, but the corresponding projections are not monomial maps when the plane is not rational, and the fibre counting by Bernstein's theorem is then unavailable. The second is the multiplicity of the amoeba map, which by Theorem \ref{thm:sym} is controlled by the symmetries of $V$; for a variety with trivial stabiliser one expects $\Log|_V$ to be generically injective when $n>2k$, as it is for nonreal lines, and then $\Hm^{2k}(\A(V))=\MVol(V)$.

\begin{problem}
Prove Conjecture \ref{conj:lines}, in particular the strict inequality $\Mv(c)<(1+\sqrt2)\pi^2$, known by Proposition \ref{prop:rate} near $c=0,1,\infty$, on the whole of $\C\setminus\{0,1\}$, and the position of the saddles at apex angle $120^\circ$, and determine the extremal volumes of amoebas of generic $k$-planes in $(\C^*)^{2k+1}$, the analogue for $n=2k+1$ of \cite[Problem 3.1]{Mik17}.
\end{problem}

\begin{problem}
Corollary \ref{cor:fibre4} shows that for every real $k$-plane in general position in $(\C^*)^{2k}$ with $k\ge3$ which is not a symmetric power of a line, the amoeba map has fibres of fewer than $2^k$ points on a set of positive measure, of exactly $4$ points when $k=3$. Prove that fibres of exactly $4$ points occur also for $k\ge4$, decide whether the amoebas of such planes always have volume larger than $\pi^{2k}/2^k$, and determine which multiples of $4$ not exceeding $\binom{2k}k$ occur as fibre sizes (Corollary \ref{cor:sizes}).
\end{problem}

\begin{problem}
Find an explicit condition on the coefficients of a complex plane in $(\C^*)^5$ under which its amoeba map is injective almost everywhere, the analogue of the condition $c\notin\R$ for lines; by Theorem \ref{thm:twotoone} this holds outside a proper real algebraic subset, but the Hermitian analogue of \eqref{eq:xi} does not force the trivial solution, so that one cannot expect injectivity everywhere. Determine the extremal values of the volume of $\A(\ell)+\A(\ell)$ on the moduli space of real configurations in $(\C^*)^5$, and decide whether its supremum equals $\frac{1+2\sqrt2}4\pi^4$. For $k=2$ and $5\le n\le7$, find an explicit condition on nonreal points $\varpi_m$ under which the amoeba map of $\ell\odot\ell$ is injective everywhere, replacing the condition of Theorem \ref{thm:nonrealsq}(a), which needs $n\ge8$. Prove Conjecture \ref{conj:k3} by completing the case analysis of the degenerate configurations of seven points in $\mathbb P^3$ begun in Theorem \ref{thm:triple}, and classify the real $k$-planes of $(\C^*)^{2k}$ with $2^k$-to-one amoeba maps for $k\ge4$; by Theorem \ref{thm:fibred} they include the fibred extensions of products of symmetric powers, and one may ask whether these, up to monomial equivalence, exhaust them. Finally, for real $k$-planes in $(\C^*)^n$ with $n+1<\binom{k+2}2$, find the explicit condition replacing the quadric condition of Section \ref{sec:conic}, and decide whether the real $k$-planes with $n+1\ge\binom{k+2}2$ violating the quadric condition, with hyperplanes in general position, are images of symmetric products of real lines.
\end{problem}

\begin{problem}
Decide whether varieties of maximal volume exist in the half-dimensional case for $k\ge2$, the question of \cite[Problems 3.1 and 3.3]{Mik17}. By Corollary \ref{cor:nomaxreal} such a variety has trivial stabiliser and is not real, by Example \ref{ex:products} the multiHarnack property does not suffice, and by Theorem \ref{thm:planesgp} no plane in general position in $(\C^*)^4$ has maximal volume. Extend Theorem \ref{thm:planesgp} to $k$-planes in general position in $(\C^*)^{2k}$ for every $k\ge2$, replacing the Minkowski sums of two ellipses in $\R^2$ of Lemma \ref{lem:ellipses} by Minkowski sums of $k$ ellipses in $\R^k$, and determine the supremum of the volume of the amoebas of planes in $(\C^*)^4$.
\end{problem}

\appendix

\section{Transversality lemmas}\label{app:transv}

In this appendix we prove the statements used in Section \ref{sec:geninj}, with the notation introduced there.

\begin{lemma}\label{lem:sep}
Let $A\subset\Z^n$ be finite, let $L$ be the subgroup generated by the differences of the elements of $A$, let $z\in(\C^*)^n$ and $u\in\T^n$. The vectors $(z^a)_{a\in A}$ and $((z/u)^a)_{a\in A}$ of $\C^A$ are linearly dependent if and only if $u\in G(L)$. In that case $f(z/u)=u^{-a_0}f(z)$ for every $f$ with support in $A$ and any $a_0\in A$.
\end{lemma}

\begin{proof}
All the coordinates $z^a$ are nonzero, so the two vectors are dependent exactly when $(z/u)^a=\mu z^a$ for some $\mu$ and all $a\in A$, that is, when $u^{-a}$ is independent of $a\in A$. This means $u^{a-b}=1$ for all $a,b\in A$, that is, $u\in G(L)$. Then $u^{-a}=u^{-a_0}$ for all $a\in A$, which gives the formula.
\end{proof}

\begin{lemma}\label{lem:U}
The set $U$ is a complex submanifold of $\C^A\times(\C^*)^n$ of dimension $|A|+k$, where $|A|=\sum_m|A_m|$, and the projection $\pi:U\to\C^A$ is a submersion. For $c\in\Vc(\Delta)$ one has $\pi^{-1}(c)=\{c\}\times V_c$. The set $\Vc(\Delta)$ contains a nonempty Zariski open subset of $\C^A$, hence it is dense in $\C^A$.
\end{lemma}

\begin{proof}
For fixed $z$, the map $c\mapsto F(c,z)$ is linear, and it is onto $\C^{n-k}$, because the $m$-th component depends only on the coefficients $c_{m,\cdot}$ and has the nonzero coefficient $z^a$ in front of $c_{m,a}$. Hence $F$ is a submersion on $\C^A\times(\C^*)^n$, $F^{-1}(0)$ is a complex manifold of dimension $|A|+n-(n-k)=|A|+k$, and $U$ is an open subset of it. A tangent vector of $F^{-1}(0)$ at $(c,z)$ is a pair $(\dot c,\dot z)$ with $d_cF\,\dot c+d_zF\,\dot z=0$. If $\rk d_zF(c,z)=n-k$, then for every $\dot c$ this equation has a solution $\dot z$, so $d\pi$ is onto at the points of $U$. The equality $\pi^{-1}(c)=\{c\}\times V_c$ for $c\in\Vc(\Delta)$ is the definition of $\Vc(\Delta)$. Finally, $F^{-1}(0)$ is a smooth algebraic variety, and by the algebraic version of Sard's theorem (generic smoothness in characteristic zero \cite[Corollary III.10.7]{Har77}) there is a nonempty Zariski open subset $W\subset\C^A$ over which $\pi|_{F^{-1}(0)}$ is smooth. For $c\in W$ the differential of $F(c,\cdot)$ has rank $n-k$ along $V_c$. Intersecting $W$ with the Zariski open set where all vertex coefficients are nonzero, we obtain a nonempty Zariski open subset of $\C^A$ contained in $\Vc(\Delta)$.
\end{proof}

\begin{lemma}\label{lem:submersion}
For every $S\subsetneq\{1,\dots,n-k\}$, the map $\Psi_S:\mathcal Y_S\to\C^{n-k-s}$ is a submersion.
\end{lemma}

\begin{proof}
Let $(c,z,u)\in\mathcal Y_S$. For $m\notin S$ we have $u\notin G_m$, so by Lemma \ref{lem:sep} the vectors $(z^a)_{a\in A_m}$ and $((z/u)^a)_{a\in A_m}$ are linearly independent in $\C^{A_m}$. Given $\eta\in\C^{n-k-s}$, choose for every $m\notin S$ a vector $\dot c_m\in\C^{A_m}$ with $\sum_a\dot c_{m,a}z^a=0$ and $\sum_a\dot c_{m,a}(z/u)^a=\eta_m$, and put $\dot c_m=0$ for $m\in S$. The curve $\tau\mapsto(c+\tau\dot c,z,u)$, $\tau\in\R$ small, stays in $\mathcal Y_S$. Indeed, $F(c+\tau\dot c,z)=F(c,z)+\tau\big(\sum_a\dot c_{m,a}z^a\big)_m=0$, the differential $d_zF$ depends continuously on $\tau$, so its rank stays $n-k$ for small $\tau$, and $u$ does not move. Along this curve $\Psi_S$ changes by $\tau\eta$, because $\Psi_S$ is affine in $c$. Hence $d\Psi_S$ is onto.
\end{proof}

\begin{lemma}\label{lem:dim}
There is a Lebesgue null set $\mathcal N_S\subset\C^A$ such that for $c\notin\mathcal N_S$ the set $W_{S,c}=\{(z,u):(c,z,u)\in\Psi_S^{-1}(0)\}$ is empty or a submanifold of $V_c\times G_S^\circ$ of dimension $e_S=4k-n-d_S+2s$, and $e_S\le2k-1$. For such $c$ the projection $Z_{S,c}=\{z:(z,u)\in W_{S,c}\text{ for some }u\}$ has $2k$-dimensional volume zero in $V_c$.
\end{lemma}

\begin{proof}
By Lemma \ref{lem:submersion}, $W_S=\Psi_S^{-1}(0)$ is a submanifold of $\mathcal Y_S$ of dimension $2|A|+2k+n-d_S-2(n-k-s)$. Let $\mathcal N_S$ be the set of critical values of the smooth map $\pi\circ\mathrm{pr}:W_S\to\C^A$, a null set by Sard's theorem, since $W_S$ is second countable. For a regular value $c$, the fibre $W_S\cap\pi^{-1}(c)$ is empty or a submanifold of dimension $\dim W_S-2|A|=4k-n-d_S+2s=e_S$, which is $W_{S,c}$. By $(\mathrm L_2)$, $d_S\ge2s+2k-n+1$, hence $e_S\le4k-n+2s-2s-2k+n-1=2k-1$. The set $Z_{S,c}$ is the image of the manifold $W_{S,c}$ of dimension $e_S<2k$ under the smooth map $(z,u)\mapsto z$ into the $2k$-dimensional manifold $V_c$. Every point of $W_{S,c}$ is critical for this map, so $Z_{S,c}$ is null by Sard's theorem.
\end{proof}

\end{document}